\documentclass[10pt,leqno]{amsart}
\usepackage[T1]{fontenc}
\usepackage[utf8]{inputenc}
\usepackage{longtable}
\usepackage{mathtools}
\usepackage{url}
\usepackage{amsbsy}
\usepackage{amstext}
\usepackage{amsthm}
\usepackage{amssymb}
\usepackage{microtype}
\usepackage[bookmarks=false,
 breaklinks=false,pdfborder={0 0 1},backref=false,colorlinks=false]
 {hyperref}
\hypersetup{pdftitle={Holomorphic sectional curvature and the Calabi operator on Kahler C-spaces},
 pdfauthor={Xiaolong Li and Peter Petersen},
 hypertexnames=false}

\makeatletter

\providecommand{\tabularnewline}{\\}

\numberwithin{equation}{section}
\numberwithin{figure}{section}

\@ifundefined{date}{}{\date{}}
\usepackage{needspace}
\usepackage{amsfonts}\usepackage{array}
\usepackage{booktabs}
\newcommand{\CQB}{\operatorname{CQB}}
\newcommand{\dCQB}{{}^{d}\!\operatorname{CQB}}
\DeclareMathOperator{\Ric}{Ric}
\theoremstyle{plain}
\newtheorem{theorem}{Theorem}[section]
\newtheorem{proposition}[theorem]{Proposition}\newtheorem{lemma}[theorem]{Lemma}\newtheorem{corollary}[theorem]{Corollary}\theoremstyle{remark}
\newtheorem{remark}[theorem]{Remark}\newtheorem{example}[theorem]{Example}

\author[Li]{Xiaolong Li}
\address{Department of Mathematics and Statistics, Auburn University,
Auburn, AL 36849, USA}
\email{xil0005@auburn.edu}
\thanks{Xiaolong Li's research is partially supported by NSF-DMS
\#2553660 and a start-up grant at Auburn University.}

\author[Petersen]{Peter Petersen}
\address{University of California, Los Angeles, 520 Portola Plaza, CA, 90095}
\email{petersen@math.ucla.edu}

\let\calabimaketitle\maketitle
\renewcommand{\maketitle}{}
\AddToHook{env/abstract/after}{%
  \calabimaketitle}

\providecommand{\theoremname}{Theorem}
\theoremstyle{plain}

\makeatother

\providecommand\theoremname{Theorem}

\begin{document}
\title{Holomorphic sectional curvature and the Calabi operator on Kähler
C-spaces}

\maketitle
\subjclass[2020]{Primary 53C55; Secondary 53C30, 32M10, 32Q20} 
\keywords{Kähler C-space, generalized flag manifold, Calabi curvature operator,
holomorphic sectional curvature, Kähler-Einstein metrics.}
\begin{abstract}
We study holomorphic sectional curvature and the Calabi curvature
operator on compact Kähler manifolds. First we show that every homogeneous
Kähler metric on a generalized flag manifold has positive holomorphic
sectional curvature. For the standard homogeneous Kähler-Einstein
metric with Einstein constant 2, it is shown that the spectrum of
the Calabi operator lies in $(-2,2]$. We determine its sharp upper
bound in terms of the anticanonical degrees and characterize the 2-eigenspace.
The spectral bounds imply nonnegative cross quadratic bisectional
curvature and positive dual cross quadratic bisectional curvature,
proving conjectures of Ni and Zheng. We also show that these conditions
do not characterize homogeneity: we construct a nonhomogeneous Kähler-Einstein
13-fold whose Calabi spectrum lies strictly inside $(-2,2)$ after
the same normalization, as well as nonhomogeneous examples in dimensions
4 and 5. The proofs involve an integral identity for holomorphic sectional
curvature, projective embeddings, root-space estimates, and explicit
metrics on projective bundles. 
\end{abstract}

\section{Introduction}

\label{sec:1}

Curvature conditions that characterize complex homogeneous spaces
form a central theme in Kähler geometry. Siu and Yau \cite{SY80}
proved that a compact Kähler manifold with positive holomorphic bisectional
curvature is complex projective space. Mori's characterization by
ampleness of the tangent bundle \cite{Mor79} places this conclusion
in a broader algebraic setting. In the case of nonnegative holomorphic
bisectional curvature, Mok and Zhong showed in \cite{MZ86} that if
the metric is also Einstein then the space is locally symmetric. Mok
generalized this in \cite{Mok88} by dropping the Einstein condition
and showed that the universal cover is a product of a Euclidean factor,
a locally symmetric space, and factors that are biholomorphic to projective
spaces. In particular, nonnegative bisectional curvature is too restrictive
to include all compact homogeneous Kähler manifolds.

Compact simply connected homogeneous Kähler manifolds are also called
Kähler C-spaces. They are generalized flag manifolds $M=G/K$, where
$G$ is compact connected semisimple and $K$ is the centralizer of
a torus, or equivalently $G^{\mathbb{C}}/P$ for a parabolic subgroup
$P$. Kähler C-spaces form a hierarchy for each fixed simple compact
group $G$ as there is a submersion $G/K\rightarrow G/L$ as long
as $K\subset L$. The initial object is the flag variety $G/T^{r}$,
where $r$ is the rank of $G$ and $T^{r}$ is a maximal torus. All
other generalized flags have compact isotropy $K\supset T^{r}$. Up
to a finite quotient, $K$ is the product of its connected center
and its semisimple derived subgroup. When $K\subsetneqq G$ is maximal
$G/K$ becomes a terminal object in this hierarchy. These include,
but are not limited to, the irreducible Hermitian symmetric spaces
and all have $b_{2}=1$. Borel's formula \cite{Bor53} determines
the Poincaré polynomials, and hence the Betti numbers, of these homogeneous
spaces; see also \cite{PS25} for explicit computations.

Borel \cite{Bor54} described the invariant complex structures of
Kähler C-spaces, and Koszul's canonical Hermitian form \cite{Kos55}
describes the invariant Ricci form and standard (canonical) homogeneous
Kähler-Einstein metric. Matsushima \cite[Theorem~3]{Mat72} showed
that every Kähler-Einstein metric on such a manifold is obtained from
the standard metric by a holomorphic automorphism and a positive rescaling.
The Borel--Remmert decomposition \cite{BR62} expresses a general
compact homogeneous Kähler manifold as a product of a complex torus
and a Kähler C-space; see also \cite[Chapter~8]{Bes87}. We consider
the simply connected case. Kähler-Einstein Fano varieties that are
not Kähler C-spaces were discovered by Sakane \cite{Sak86}, followed
by \cite{KS88}. Tian's seminal paper \cite{Tian87} constructed examples
of Fano varieties with Kähler-Einstein metrics and finite automorphism
group.

A related algebraic problem is the Campana--Peternell conjecture
\cite{CP91}: a Fano manifold with nef tangent bundle should be rational
homogeneous. Homogeneous tangent bundles are generated by global vector
fields and hence are nef. This motivates the search for curvature
conditions that hold on all Kähler C-spaces and might characterize
them. We study two such conditions, expressed through holomorphic
sectional curvature and the spectrum of the Calabi curvature operator.

\subsection{Holomorphic sectional curvature}

For $0\ne X\in T^{1,0}M$, write 
\[
H_{g}(X)=\frac{R(X,\overline{X},X,\overline{X})}{|X|^{4}_{g}}.
\]
Our convention gives $H=2$ for the Fubini--Study metric with potential
$\log(1+|z|^{2})$. Positive holomorphic sectional curvature is weaker
than positive bisectional curvature, but still forces a compact Kähler
manifold to be projective and rationally connected \cite{Yan18}.
Metrics on projective bundles \cite{AHZ18,YZ19} and on Grassmannian
and flag bundles \cite{Mag22} illustrate the range of geometries
compatible with this condition.

For homogeneous metrics, Itoh \cite[Theorem~3.1]{Ito78} proved positivity
for Kähler C-spaces with $b_{2}=1$ and at most two irreducible isotropy
modules in $T^{1,0}M$. This includes all classical spaces with $b_{2}=1$
and several exceptional spaces. Lohove \cite[Theorems~22, 47, and 54]{Loh19}
generalized this to establish positivity for homogeneous Kähler metrics
with at most three irreducible isotropy summands, for the standard
metrics of type $G_{2}$, and for the standard metrics on classical
flag manifolds of rank at most four.

\begin{theorem}\label{theorem:1.1} Let $M=G/K$ be a positive-dimensional
compact generalized flag manifold with a fixed $G$-invariant complex
structure. Every $G$-invariant Kähler metric on $M$ has strictly
positive holomorphic sectional curvature. \end{theorem}

In particular, the standard homogeneous Kähler-Einstein metric has
$H>0$, as conjectured in \cite[Conjecture~1.8]{NZ25}. The proof
works for all homogeneous Kähler metrics. We describe the metric by
a positive grading of the tangent root spaces and dampen a tangent
vector $X$ exponentially according to this grading. The holomorphic
curvature of $X$ becomes an integral of squared norms of Hermitian
Lie algebra elements. Pairing with the grading element proves strict
positivity and gives a quantitative lower bound. Section~\ref{sec:3}
derives this identity, including the isotropy contribution.

\subsection{The Calabi curvature operator}

The Kähler curvature tensor defines a Hermitian endomorphism $\mathcal{C}$
of $\operatorname{Sym}^{2}T^{1,0}M$, which we call the \emph{Calabi
curvature operator}; compare \cite{CV60,Ito78}. In a unitary frame
it is characterized by 
\[
\langle\mathcal{C}S,S\rangle=\sum_{p,q,r,s}R_{p\bar{q}r\bar{s}}S_{pr}\overline{S_{qs}},\qquad S_{pr}=S_{rp}.
\]
The norm is the one induced from the tensor product, so $\|S\|^{2}=\sum_{p,r}|S_{pr}|^{2}$.
Holomorphic sectional curvature comes from $S=X\otimes X$. For $S=\frac{1}{2}(X\otimes Y+Y\otimes X)$,
the Kähler symmetries give $\langle\mathcal{C}S,S\rangle=R(X,\overline{X},Y,\overline{Y})$;
dividing by $|X|^{2}|Y|^{2}$ gives the holomorphic bisectional curvature
of nonzero $X,Y\in T^{1,0}M$. The spectrum $\operatorname{Spec}(\mathcal{C})$
is calculated for irreducible Hermitian symmetric spaces in \cite{CV60}
and \cite{Bor60}. Projective space has only one eigenvalue while
in all other cases there are two eigenvalues of opposite sign. These
spectral estimates enter into the Calabi--Vesentini vanishing argument
for tangent-bundle-valued cohomology. In the normalization $\Ric=2g$,
the relevant lower bound is $\mathcal{C}>-2\operatorname{Id}$. Ni
\cite[Theorem~2.6(ii)]{Ni21} proved that a compact Kähler manifold
with positive $\dCQB$ satisfies $H^{q}(M,T^{1,0}M)=0$ for $1\leq q\leq\dim_{\mathbb{C}}M$.

\begin{theorem}\label{theorem:1.2} Let $M$ be a positive-dimensional
compact Kähler C-space with its standard homogeneous Kähler-Einstein
metric, normalized by $\Ric=2g$. The spectrum of its Calabi curvature
operator satisfies 
\[
\operatorname{Spec}(\mathcal{C})\subset(-2,2].
\]
\end{theorem}

The upper bound has a precise geometric description. Choose a set
$\Phi$ of simple roots and let $\Phi_{K}\subset\Phi$ be the simple
roots for the Levi factor and $\Phi_{T}=\Phi-\Phi_{K}$ the painted
simple roots. For $\alpha\in\Phi_{T}$, let $k_{\alpha}$ be the degree
of the anticanonical bundle of $M$ on the simple Schubert curve corresponding
to $\alpha$, and set 
\[
k_{*}=\min_{\alpha\in\Phi_{T}}k_{\alpha}.
\]
These integers satisfy $k_{\alpha}\geq2$. We prove the sharper identity
\[
\lambda_{\max}(\mathcal{C})=\frac{4}{k_{*}}\leq2,
\]
and for unit vectors $e_{\alpha}$ in the root spaces, 
\[
\ker(2\operatorname{Id}-\mathcal{C})=\bigoplus_{\substack{\alpha\in\Phi_{T}\\
k_{\alpha}=2
}
}\mathbb{C}(e_{\alpha}\otimes e_{\alpha}).
\]
Thus $\lambda_{\max}(\mathcal{C})<2$ precisely when every simple
root $\alpha\in\Phi_{T}$ is adjacent to a root for the Levi factor,
or equivalently when $M$ has no homogeneous projection with fiber
$\mathbb{P}^{1}$.

The proof of the upper bound uses the fundamental representations
to embed the generalized flag isometrically into a product of projective
spaces, suitably scaled by the anticanonical degrees. The Gauss equation
then bounds $\mathcal{C}$ from above. For the lower bound, torus
invariance splits $\mathcal{C}$ into blocks indexed by sums of tangent
roots as observed by Lohove in \cite{Loh19}. A root-counting argument
gives a quantitative bound in the simply-laced types $A$, $D$, and
$E$. Scaled matrix estimates treat types $B$, $C$, and $F_{4}$,
while the $G_{2}$ blocks have size at most three and admit direct
formulas.

The background for the above theorem starts with Ni's work in \cite{Ni21}
where the cross quadratic bisectional curvature $\CQB$ and its dual
$\dCQB$ are introduced. When $\Ric=\lambda g$, these forms satisfy
\[
\CQB(A)=\lambda\|A\|^{2}-\langle\mathcal{C}S,S\rangle,\qquad\dCQB(A)=\lambda\|A\|^{2}+\langle\mathcal{C}S,S\rangle,
\]
where $S$ is the symmetric part of $A$; Section~\ref{sec:2} specifies
the two coefficient conventions. In particular, Theorem~\ref{theorem:1.2}
implies $\CQB\geq0$ and $\dCQB>0$, proving \cite[Conjecture~1.5]{NZ25}.
The sharp upper bound gives $\CQB(A)\geq(2-4/k_{*})\|A\|^{2}$ and
proves the strict-positivity assertion proposed on page~1719 of \cite{NZ25}.
Earlier cases include the classical spaces with $b_{2}=1$ \cite[Theorems~5.4 and 6.7]{Ni21},
flags of type $A$, and $\mathrm{SU}(6)/\mathrm{S}(\mathrm{U}(2)^{3})$
\cite[Theorems~1.6 and 1.7]{NZ25}. These conditions differ from quadratic
bisectional curvature, whose sign was studied by Chau and Tam \cite{CT13},
and from dual-Nakano semipositivity of the generally non-Kähler normal
homogeneous metric \cite{GO26}.

\subsection{Counterexamples to homogeneity}

The spectral condition in Theorem~\ref{theorem:1.2} does not characterize
homogeneous Kähler-Einstein manifolds. We give an explicit 13-dimensional
example, together with two lower-dimensional non-Einstein examples
for the corresponding curvature conditions involving the Ricci eigenvalues.

Throughout, $\mathbb{P}_{\mathrm{lines}}(E)$ denotes the bundle of
complex lines in a holomorphic vector bundle $E$.

\begin{theorem}\label{theorem:1.3} The following compact complex
manifolds are not homogeneous. 
\begin{enumerate}
\item Let $B=\mathbb{P}^{5}\times\mathbb{P}^{5}$, $L_{1}=\mathcal{O}_{B}(-1,0)$,
and $L_{2}=\mathcal{O}_{B}(0,-1)$. The 13-fold 
\[
M=\mathbb{P}_{\mathrm{lines}}(L^{\oplus2}_{1}\oplus L^{\oplus2}_{2})
\]
admits an explicit smooth Kähler-Einstein metric $g$ with $\Ric(g)=g$
and 
\[
-\frac{1}{2}\operatorname{Id}\leq\mathcal{C}_{g}\leq\frac{11}{12}\operatorname{Id}.
\]
Both spectral bounds are attained. 
\item The 4-fold $X_{2}=\mathbb{P}_{\mathrm{lines}}(\mathcal{O}^{\oplus2}_{\mathbb{P}^{2}}\oplus\mathcal{O}_{\mathbb{P}^{2}}(-1))$
admits an explicit Kähler metric satisfying 
\[
\CQB(A)\geq0,\qquad\dCQB(A)\geq\frac{1}{8}\|A\|^{2}.
\]
\item The 5-fold $X_{3}=\mathbb{P}_{\mathrm{lines}}(\mathcal{O}^{\oplus2}_{\mathbb{P}^{3}}\oplus\mathcal{O}_{\mathbb{P}^{3}}(-1))$
admits an explicit Kähler metric satisfying 
\[
\CQB(A)\geq\frac{1}{4}\|A\|^{2},\qquad\dCQB(A)\geq\frac{1}{4}\|A\|^{2}.
\]
The metrics in parts~(2) and~(3) are not Einstein. 
\end{enumerate}
\end{theorem}

Rescaling the metric in part~(1) to $g/2$ gives Einstein constant
$2$ and Calabi spectrum in $[-1,11/6]\subset(-2,2)$. For the metric
$g$ as stated in part~(1), the spectral bounds give the sharp inequalities
$\CQB(A)\geq\frac{1}{12}\|A\|^{2}$ and $\dCQB(A)\geq\frac{1}{2}\|A\|^{2}$.
Thus part~(1) disproves the Kähler-Einstein assertion of \cite[Conjecture~1.9]{NZ25},
and parts~(2) and~(3) disprove its general Kähler assertion in smaller
dimensions. In the non-Einstein case, a bound on $\mathcal{C}$ alone
does not establish these inequalities.

The Einstein metric belongs to the projective-bundle framework of
Koiso and Sakane \cite{KS88} and the admissible metrics of Apostolov,
Calderbank, Gauduchon, and Tønnesen-Friedman \cite{ACGTF08}. We describe
a direct potential, verify the Einstein equation, and check smoothness
at both ends of the fiber coordinate. The curvature estimate controls
the entire Calabi operator. The fivefold $X_{3}$ already appears
in \cite[Example~1.4 and p.~837]{NWZ21} as a nonhomogeneous manifold
with positive orthogonal Ricci curvature.

Nonhomogeneity has a unified explanation: on each example, the tangent
bundle restricted to a suitable rational curve has a quotient line
bundle of degree $-1$. A homogeneous tangent bundle is globally generated
and cannot have such a quotient. These tangent bundles are therefore
not nef, so the examples do not contradict the Campana--Peternell
conjecture.

Section~\ref{sec:2} fixes the Calabi operator, $\CQB$ and $\dCQB$
conventions and recalls the invariant curvature formulas from Lohove
\cite{Loh19}. Sections~\ref{sec:3}, \ref{sec:4}, and~\ref{sec:5}
prove Theorems~\ref{theorem:1.1}, \ref{theorem:1.2}, and~\ref{theorem:1.3},
respectively.

\medskip{}

\noindent\textbf{Note added during preparation.} While we were preparing
this manuscript, the paper of Jian and Zhang \cite{JZ26} appeared.
It proves that every Kähler C-space of classical type, with its canonical
Kähler-Einstein metric, has $\CQB\geq0$ and $\dCQB>0$. In our normalization
$\Ric=2g$, their conclusion is equivalent to $\operatorname{Spec}(\mathcal{C})\subset(-2,2]$
and thus overlaps with the classical-type part of Theorem~\ref{theorem:1.2}.

\medskip{}

\textbf{AI Disclosure:} The authors used OpenAI’s ChatGPT as an assistive
tool in the development of the paper. The authors take full responsibility
for all mathematical arguments.

\section{Preliminaries}

\label{sec:2}

\subsection{The Calabi curvature operator}

Let $(M,g)$ be a Kähler manifold of complex dimension $n$. We use
Latin letters for indices since Greek letters are used for roots.
In the actual calculations we will use the roots as indices. The curvature
is given by 
\begin{equation}
R_{a\bar{b}c\bar{d}}=-\partial_{a}\partial_{\bar{b}}g_{c\bar{d}}+g^{p\bar{q}}(\partial_{a}g_{c\bar{q}})(\partial_{\bar{b}}g_{p\bar{d}}),\label{eq:2.1}
\end{equation}
so that $\Ric_{a\bar{b}}=\sum_{c}R_{a\bar{b}c\bar{c}}$ in a unitary
frame. The Kähler symmetries make 
\[
\langle\mathcal{C}S,S\rangle=\sum_{a,b,c,d}R_{a\bar{b}c\bar{d}}S_{ac}\overline{S_{bd}},\qquad S\in\operatorname{Sym}^{2}T^{1,0}M,
\]
the quadratic form of a Hermitian operator. We call $\mathcal{C}$
the \emph{Calabi curvature operator}. Symmetric tensors carry the
inner product induced from $T^{1,0}M\otimes T^{1,0}M$; thus $\|S\|^{2}=\sum_{a,c}|S_{ac}|^{2}$.
In particular, an off-diagonal tensor $(e_{a}\otimes e_{c}+e_{c}\otimes e_{a})/\sqrt{2}$, where $a\neq c$,
has unit length, as does the square of a unit vector $e_{a}\otimes e_{a}$.

The Fubini--Study metric with potential $\log(1+|z|^{2})$ has 
\[
R_{a\bar{b}c\bar{d}}=\delta_{ab}\delta_{cd}+\delta_{ad}\delta_{cb},\qquad\mathcal{C}=2\operatorname{Id},\qquad\Ric=(n+1)g.
\]
Holomorphic sectional curvature is the restriction of its quadratic
form to rank-one symmetric tensors: 
\[
H(X)=\frac{\langle\mathcal{C}(X\otimes X),X\otimes X\rangle}{|X|^{4}}.
\]
Positivity of $H$ does not imply that $\mathcal{C}$ is positive.

\subsubsection{Relation to the cross quadratic forms}

We record the definitions of \cite{Ni21,NZ25} to relate the spectral
results to those conjectures and to treat the non-Einstein examples
in Section~\ref{sec:5}. In a unitary frame of eigenvectors for $\mathrm{Ric}$,
write the eigenvalues as $r_{1},\ldots,r_{n}$. For $A:T^{0,1}M\to T^{1,0}M$,
let $A_{pq}$ be the coefficient of $e_{p}$ in $A(\bar{e}_{q})$.
For $A:T^{1,0}M\to T^{0,1}M$, use the complex conjugate of the coefficient
of $\bar{e}_{p}$ in $A(e_{q})$. In either case set 
\[
\mathcal{M}(A)=\sum_{p,q}r_{p}|A_{pq}|^{2},\qquad\mathcal{R}(A)=\sum_{p,q,r,s}R_{p\bar{q}r\bar{s}}A_{pr}\overline{A_{qs}}.
\]
Then 
\begin{align}
\CQB(A) & =\mathcal{M}(A)-\mathcal{R}(A),\nonumber \\
\dCQB(A) & =\mathcal{M}(A)+\mathcal{R}(A).\label{eq:2.2}
\end{align}
These contractions are independent of the frame. If we decompose $A=S+K$
into its orthogonal symmetric and skew-symmetric parts, then the Kähler
symmetry implies that $\mathcal{R}(A)=\langle\mathcal{C}S,S\rangle$.
In case $\Ric=\lambda g$ 
\begin{align}
\CQB(A) & =\lambda\|A\|^{2}-\langle\mathcal{C}S,S\rangle,\nonumber \\
\dCQB(A) & =\lambda\|A\|^{2}+\langle\mathcal{C}S,S\rangle.\label{eq:2.3}
\end{align}
For $\lambda>0$, the conditions $\CQB\geq0$ and $\dCQB>0$ are therefore
equivalent to 
\[
-\lambda\operatorname{Id}<\mathcal{C}\leq\lambda\operatorname{Id}.
\]
The skew-symmetric part contributes $\lambda\|K\|^{2}$ to both forms.

\subsection{Flag manifolds, root spaces, and invariant metrics}

We use the standard description of homogeneous Kähler manifolds in
\cite[Chapter~8]{Bes87} and \cite[Sections~1.1--1.4]{Loh19}. Let
$M=G/K$, with $G$ compact connected semisimple and $K$ the centralizer
of a torus. We may replace $G$ by its simply connected cover. Let
$\mathfrak{g}$ be the Lie algebra and $\mathfrak{g}^{\mathbb{C}}$
the complexification. Fix a maximal torus of $K$ and let $\mathfrak{h}\subset\mathfrak{k}^{\mathbb{C}}$
be the complexification of its Lie algebra, so that $\mathfrak{h}$
is preserved by compact conjugation. Choose simple roots $\Phi$ and
positive roots $\Delta^{+}$ compatible with the complex structure.
Let $\Phi_{K}$ index the simple roots of the Levi factor and let
$\Phi_{T}=\Phi\setminus\Phi_{K}$ be the painted simple roots or simple
tangent roots. Let $\Delta_{\Phi_{K}}$ be the root subsystem generated
by $\Phi_{K}$, and set 
\[
\Delta^{+}_{\Phi_{K}}=\Delta_{\Phi_{K}}\cap\Delta^{+},\qquad\Delta^{+}_{\Phi_{T}}=\Delta^{+}\setminus\Delta^{+}_{\Phi_{K}}.
\]
Thus $\Delta^{+}_{\Phi_{T}}$ is the set of positive tangent roots.
The corresponding decompositions are 
\[
\mathfrak{k}^{\mathbb{C}}=\mathfrak{h}\oplus\bigoplus_{\alpha\in\Delta_{\Phi_{K}}}\mathfrak{g}_{\alpha},\qquad\mathfrak{m}_{\pm}=\bigoplus_{\alpha\in\Delta^{+}_{\Phi_{T}}}\mathfrak{g}_{\pm\alpha}.
\]
Thus $T^{1,0}_{o}M=\mathfrak{m}_{+}$ and $T^{0,1}_{o}M=\mathfrak{m}_{-}$,
where $o=eK$, and $\mathfrak{m}^{\mathbb{C}}=\mathfrak{m}_{-}\oplus\mathfrak{m}_{+}$
is the tangent complement of $\mathfrak{k}^{\mathbb{C}}$. In the
realization $M=G^{\mathbb{C}}/P$, our parabolic $P$ contains the
Borel algebra that uses the negative root spaces in $\mathfrak{m}_{-}$.
Nodes indexed by $\Phi_{K}$ are called white, and the remaining nodes
painted. Adjacent nodes have nonzero root inner product. Simply-laced
root systems have equal root lengths and correspond to types $A,D,E$;
the other irreducible types are $B,C,F_{4},G_{2}$.

Let $B$ be an invariant symmetric complex bilinear form, negative
definite on $\mathfrak{g}$. Complex conjugation $\sigma$ defines
$U^{*}=-\sigma(U)$, so 
\[
[U,V]^{*}=[V^{*},U^{*}],\qquad\langle U,V\rangle_{B}=B(U,V^{*})>0\quad(U=V\ne0).
\]
For $\alpha\in\mathfrak{h}^{*}$, define $H_{\alpha}$ by $B(H_{\alpha},H)=\alpha(H)$
for all $H\in\mathfrak{h}$, and put $(\alpha,\beta)=B(H_{\alpha},H_{\beta})$.
Thus $H_{\alpha}$ is the $B$-dual of $\alpha$; for a root $\alpha$,
the Cartan element of its standard $\mathfrak{sl}_{2}$ triple is
$h_{\alpha}=2H_{\alpha}/|\alpha|^{2}$. Choose root vectors with real
structure constants satisfying 
\begin{gather*}
E^{*}_{\alpha}=E_{-\alpha},\qquad B(E_{\alpha},E_{-\alpha})=z_{\alpha}>0,\\{}
[E_{\alpha},E_{\beta}]=N_{\alpha\beta}E_{\alpha+\beta}\quad\text{if }\alpha+\beta\text{ is a root},\\{}
[E_{\alpha},E_{-\alpha}]=z_{\alpha}H_{\alpha}.
\end{gather*}
When $\alpha+\beta$ is neither a root nor zero, the bracket $[E_{\alpha},E_{\beta}]=0$
and $N_{\alpha\beta}=0$. Retaining $z_{\alpha}$ allows the use of
convenient matrix root vectors in Section~\ref{sec:4}. This is also
the convention in \cite[Proposition~11]{Loh19}.

An invariant Kähler metric is determined by additive metric coefficients
$g_{\alpha}$, positive on $\Delta^{+}_{\Phi_{T}}$, zero on the Levi
roots, and satisfying $g_{-\alpha}=-g_{\alpha}$. For $\alpha\in\Delta^{+}_{\Phi_{T}}$,
write $w_{\alpha}=\|E_{\alpha}\|^{2}_{g}=g_{\alpha}z_{\alpha}$. Additivity
means $g_{\alpha+\beta}=g_{\alpha}+g_{\beta}$ whenever $\alpha+\beta$
is a root. Equivalently, there is a Hermitian central element $\xi\in\mathfrak{z}(\mathfrak{k}^{\mathbb{C}})$
with $g_{\alpha}=\alpha(\xi)$. If $D=\operatorname{ad}\xi$, then
\[
|X|^{2}_{g}=B(DX,X^{*}),\qquad h(U,V)=-B(|D|U,V),
\]
for $X\in\mathfrak{m}_{+}$ and $U,V\in\mathfrak{m}^{\mathbb{C}}$,
where $h$ is the complex bilinear extension of the real metric. Here
$D$ is self-adjoint for $\langle\cdot,\cdot\rangle_{B}$, and $|D|$
acts by $|s|$ on the $s$-eigenspace of $D$. The positive, zero,
and negative subspaces of $D$ are $\mathfrak{m}_{+}$, $\mathfrak{k}^{\mathbb{C}}$,
and $\mathfrak{m}_{-}$, respectively. Section~\ref{sec:3} uses this description
for arbitrary invariant Kähler metrics.

For the canonical metric, let 
\[
\rho=\frac{1}{2}\sum_{\alpha\in\Delta^{+}}\alpha,\qquad\rho_{\Phi_{K}}=\frac{1}{2}\sum_{\alpha\in\Delta^{+}_{\Phi_{K}}}\alpha,\qquad\delta=2\rho-2\rho_{\Phi_{K}}=\sum_{\alpha\in\Delta^{+}_{\Phi_{T}}}\alpha.
\]
The sum $\delta$ is the weight of the anticanonical bundle. The anticanonical
degrees are 
\[
k_{\alpha}=\frac{2(\delta,\alpha)}{|\alpha|^{2}},\qquad\alpha\in\Phi_{T}.
\]
The invariant Ricci form is determined by $\delta$, independently
of the metric \cite{Kos55}. The normalization $\Ric=2g$ corresponds
to 
\begin{equation}
g_{\alpha}=\frac{1}{2}(\delta,\alpha),\qquad w_{\alpha}=\|E_{\alpha}\|^{2}=g_{\alpha}z_{\alpha},\qquad\Ric=2g.\label{eq:2.4}
\end{equation}
Thus $\xi=H_{\delta}/2$ and $[\omega_{g}]/(2\pi)=c_{1}(M)/2$. Formula
\eqref{eq:2.4} is the formula for the metric coefficients in \cite[Section~3.2, (3.3)]{Loh19},
scaled to the stated Einstein constant. It specifies the canonical
choice of $g_{\alpha}$; we use $g_{\alpha}$ also for arbitrary invariant
Kähler metrics. On the full flag $\mathrm{SU}(3)/T^{2}$, with simple
root lengths of square norm 2, the canonical metric coefficients of
$\alpha_{1},\alpha_{2},\alpha_{1}+\alpha_{2}$ are $1,1,2$.

\subsection{Projective realizations and the Gauss equation}

The fundamental weights $\varpi_{\beta}$ are characterized by 
\[
\frac{2(\varpi_{\beta},\alpha)}{|\alpha|^{2}}=\delta_{\alpha\beta},\qquad\alpha,\beta\in\Phi.
\]
They are dual to the simple coroots $\frac{2\alpha}{|\alpha|^{2}}$,
rather than to the simple roots. This is the normalized pairing of
\cite[Section~13.1]{Hum72}, written using the ordinary inner product;
angle brackets here denote Hermitian inner products only. For each
painted node $\alpha$, let $V_{\alpha}$ be the irreducible $G^{\mathbb{C}}$-representation
of lowest weight $-\varpi_{\alpha}$, and choose a nonzero lowest-weight
vector $v_{\alpha}\in V_{\alpha}$. The stabilizer, $P_{\alpha}$,
of $\mathbb{C}v_{\alpha}$ is the maximal parabolic obtained by deleting
node $\alpha$. Since $P=\bigcap_{\alpha\in\Phi_{T}}P_{\alpha}$,
the maps 
\[
F_{\alpha}(gP)=[gv_{\alpha}],\qquad F=(F_{\alpha})_{\alpha\in\Phi_{T}}:M\longrightarrow\prod_{\alpha\in\Phi_{T}}\mathbb{P}(V_{\alpha})
\]
yield a holomorphic embedding. Indeed, the common stabilizer is $P$,
and the differential has kernel $(\bigcap_{\alpha}\mathfrak{p}_{\alpha})/\mathfrak{p}=0$;
compactness then makes the injective immersion an embedding. This
is the standard fundamental-representation embedding; see \cite[Remark~3.7]{Kad23}
for this product realization of a partial flag variety. For type $A$,
it is the product of the Plücker embeddings of the subspaces defining
a flag.

Let $C_{\alpha}\cong\mathbb{P}^{1}$ be the $\mathrm{SL}_{2}\subset G^{\mathbb{C}}$
orbit through $o$ of the root associated with a painted simple root
$\alpha$. The cyclic $\mathfrak{sl}_{2}$-submodule generated by
$v_{\beta}$ has lowest weight $-\delta_{\alpha\beta}$, by the defining
formula for the fundamental weights. Hence $F_{\alpha}$ maps $C_{\alpha}$
to a projective line, and every $F_{\beta}$ with $\beta\ne\alpha$
is constant on $C_{\alpha}$. With $G$-invariant Hermitian inner
products on the spaces $V_{\alpha}$, the Fubini--Study forms satisfy
\[
\frac{1}{2\pi}\int_{C_{\alpha}}F^{*}_{\beta}\omega_{\mathrm{FS},\beta}=\delta_{\alpha\beta}.
\]
For $\alpha\in\Phi_{T}$, the root formula and the Einstein normalization
give 
\[
\int_{C_{\alpha}}c_{1}(M)=k_{\alpha}=\frac{2(\delta,\alpha)}{|\alpha|^{2}}=\frac{4g_{\alpha}}{|\alpha|^{2}},\qquad\frac{1}{2\pi}\int_{C_{\alpha}}\omega_{g}=\frac{k_{\alpha}}{2}.
\]
The invariant form induced by the product metric with factors $(k_{\alpha}/2)g_{\mathrm{FS},\alpha}$
therefore has the same periods as the canonical form: 
\[
\omega_{g}=\sum_{\alpha\in\Phi_{T}}\frac{k_{\alpha}}{2}F^{*}_{\alpha}\omega_{\mathrm{FS},\alpha}.
\]
Here we use the fact that invariant closed $(1,1)$-forms are determined
by their periods on the simple Schubert curves; see \cite[Proposition~2.2 and Theorem~3.1]{AB03}.
Thus $F$ is also an isometry with these scales.

\subsection{Root-sum blocks and curvature formulas}

We use the curvature formulas already calculated in \cite[Section~1.4.2, Propositions~12 and~16]{Loh19}.
We translate them into our convention in a symmetric form convenient
for calculations. Lohove's curvature endomorphism has the opposite
sign to the one in \eqref{eq:2.7} below. With the same root vectors
and metric, his component $R_{\alpha\beta\gamma\epsilon}$ is therefore
$-R_{\alpha\bar{\beta}\gamma\bar{\epsilon}}$ in our notation. His
formulas include the isotropy contributions, which must be retained
on partial flags.

In the unit frame $e_{\alpha}=E_{\alpha}/\sqrt{w_{\alpha}}$, the
torus acts on $e_{\alpha}\otimes e_{\gamma}$ with weight $\alpha+\gamma$.
Thus $\mathcal{C}$ splits into orthogonal blocks indexed by root
sums. Use the unit pair $e_{\alpha}\otimes e_{\alpha}$ for a repeated
root and $(e_{\alpha}\otimes e_{\gamma}+e_{\gamma}\otimes e_{\alpha})/\sqrt{2}$
otherwise, and write $p_{\alpha\gamma}=1$ or $2$, respectively.
The block entries are 
\begin{equation}
\mathcal{C}_{\alpha\gamma,\beta\epsilon}=\frac{\sqrt{p_{\alpha\gamma}p_{\beta\epsilon}}\,R_{\alpha\bar{\beta}\gamma\bar{\epsilon}}}{\sqrt{w_{\alpha}w_{\gamma}w_{\beta}w_{\epsilon}}},\qquad\alpha+\gamma=\beta+\epsilon=\tau.\label{eq:2.5}
\end{equation}
The numerator is evaluated on the original root vectors $E_{\alpha}$.

\begin{proposition}[Invariant curvature formula]\label{proposition:2.1}
Suppose $\alpha,\beta,\gamma,\epsilon\in\Delta^{+}_{\Phi_{T}}$ and
$\alpha+\gamma=\beta+\epsilon$. Set 
\[
s=g_{\alpha}+g_{\gamma}=g_{\beta}+g_{\epsilon},
\]
and 
\[
N_{1}=B([E_{\alpha},E_{-\beta}],[E_{\gamma},E_{-\epsilon}]),\qquad N_{2}=B([E_{\alpha},E_{-\epsilon}],[E_{\gamma},E_{-\beta}]).
\]
Then 
\begin{equation}
R_{\alpha\bar{\beta}\gamma\bar{\epsilon}}=\frac{\min(g_{\alpha},g_{\beta})\min(g_{\gamma},g_{\epsilon})N_{1}+\min(g_{\alpha},g_{\epsilon})\min(g_{\gamma},g_{\beta})N_{2}}{s}.\label{eq:2.6}
\end{equation}
Components with unequal root sums vanish. The same formula holds for
any additive positive metric coefficients in place of $g_{\alpha}$.
\end{proposition} 
\begin{proof}
This is a symmetric reformulation of \cite[Proposition~16]{Loh19}.
To make the sign and the case of vanishing metric coefficients explicit,
his connection formulas \cite[Proposition~12]{Loh19} give, for signed
tangent roots $r,t$ with tangent sum, 
\[
\Lambda_{E_{r}}E_{t}=\frac{1}{2}\left(1+\frac{|g_{t}|-|g_{r}|}{|g_{r}+g_{t}|}\right)[E_{r},E_{t}].
\]
An isotropy output is projected out of $\Lambda$. In our sign convention
the invariant curvature formula is 
\begin{equation}
\mathcal{R}(X,Y)=[\Lambda_{X},\Lambda_{Y}]-\Lambda_{[X,Y]_{\mathfrak{m}}}-\operatorname{ad}([X,Y]_{\mathfrak{k}})|_{\mathfrak{m}}.\label{eq:2.7}
\end{equation}
Here $\mathcal{R}(X,Y)$ denotes the curvature endomorphism, whereas
$\mathcal{R}(A)$ above denotes its scalar contraction. Since $h(E_{\epsilon},E_{-\epsilon})=-g_{\epsilon}z_{\epsilon}$,
its component is 
\[
-g_{\epsilon}B(\mathcal{R}(E_{\alpha},E_{-\beta})E_{\gamma},E_{-\epsilon}).
\]
The calculation in \cite[(1.5)--(1.10)]{Loh19}, with this sign change
and the Jacobi identity, gives 
\[
R_{\alpha\bar{\beta}\gamma\bar{\epsilon}}=\left(\min(g_{\gamma},g_{\epsilon})-\frac{g_{\gamma}g_{\epsilon}}{s}\right)N_{1}+\left(\frac{g_{\gamma}g_{\epsilon}}{s}-\max\{g_{\gamma}-g_{\beta},0\}\right)N_{2}.
\]
The isotropy term supplies the coefficient $g_{\gamma}=g_{\epsilon}$
when $g_{\alpha}=g_{\beta}$, including the Cartan case $\alpha=\beta$.
Finally, 
\[
s\min(g_{\gamma},g_{\epsilon})-g_{\gamma}g_{\epsilon}=\min(g_{\alpha},g_{\beta})\min(g_{\gamma},g_{\epsilon}),
\]
\[
g_{\gamma}g_{\epsilon}-s\max\{g_{\gamma}-g_{\beta},0\}=\min(g_{\alpha},g_{\epsilon})\min(g_{\gamma},g_{\beta}),
\]
which proves \eqref{eq:2.6}.

The torus invariance shows that curvatures where $\alpha+\gamma\neq\beta+\epsilon$
vanish. 
\end{proof}

For the diagonal entries, the same source \cite[Proposition~16, cases Ia and Ib]{Loh19}
gives 
\begin{equation}
R_{\alpha\bar{\alpha}\gamma\bar{\gamma}}=\begin{cases}
g_{\alpha}z_{\alpha}z_{\gamma}(\alpha,\gamma)+\dfrac{g^{2}_{\alpha}}{g_{\alpha+\gamma}}z_{\alpha+\gamma}N^{2}_{\alpha\gamma}, & \gamma-\alpha\in\Delta^{+}_{\Phi_{T}},\\[4pt]
g_{\gamma}z_{\alpha}z_{\gamma}(\alpha,\gamma)+\dfrac{g^{2}_{\gamma}}{g_{\alpha+\gamma}}z_{\alpha+\gamma}N^{2}_{\alpha\gamma}, & \gamma-\alpha\notin\Delta^{+}_{\Phi_{T}}.
\end{cases}\label{eq:2.8}
\end{equation}
The sum-root term is omitted if $\alpha+\gamma$ is not a root. In
particular, $R_{\alpha\bar{\alpha}\alpha\bar{\alpha}}=g_{\alpha}z^{2}_{\alpha}|\alpha|^{2}$.
The numerator in the top line is $g^{2}_{\alpha}$, as in Lohove's
formula; this corrects $g_{\alpha}g_{\gamma}$ in the top line of
\cite[(4.2)]{NZ25}.

For example, on the flag $\mathrm{SU}(3)/T^{2}$, take $z_{\alpha}=1$
and $\left|\alpha\right|^{2}=2$ for simple roots. The two simple
roots have inner product $-1$, $N^{2}_{\alpha_{1}\alpha_{2}}=1$,
and $g_{\alpha_{1}}=g_{\alpha_{2}}=1$ . Therefore 
\[
R_{\alpha_{1}\bar{\alpha}_{1}\alpha_{2}\bar{\alpha}_{2}}=-1+\frac{1}{2}=-\frac{1}{2}.
\]
Thus the Calabi operator has eigenvalue $-1$ with eigenvector 
\[
\frac{1}{\sqrt{2}}\left(e_{\alpha_{1}}\otimes e_{\alpha_{2}}+e_{\alpha_{2}}\otimes e_{\alpha_{1}}\right).
\]
The eigenvalue $2$ is attained on both $e_{\alpha_{1}}\otimes e_{\alpha_{1}}$
and $e_{\alpha_{2}}\otimes e_{\alpha_{2}}$. This small example exhibits
both the normalization of the pair basis and the need to control the
lower end of the spectrum.

\subsection{Projective-bundle conventions}

For a holomorphic vector bundle $E\to B$, $\mathbb{P}_{\mathrm{lines}}(E)$
parametrizes one-dimensional subspaces of its fibers. It has tautological
subbundle $\mathcal{O}_{\mathbb{P}(E)}(-1)\subset\pi^{*}E$, with
dual $\mathcal{O}_{\mathbb{P}(E)}(1)$. On a product of projective
spaces, $\mathcal{O}(a,b)$ denotes the tensor product of the pullbacks
from the two factors. A line subbundle $L\subset E$ defines a section,
along which 
\[
T_{\mathbb{P}(E)/B}|_{B}=\operatorname{Hom}(L,E/L),\qquad T\mathbb{P}(E)|_{B}=TB\oplus\operatorname{Hom}(L,E/L).
\]
The second identity uses the differential of the section to split
the tangent projection.

\begin{lemma}\label{lemma:2.2} If an embedded curve $C\cong\mathbb{P}^{1}$
in a compact complex manifold $Y$ admits a surjective holomorphic
bundle map $TY|_{C}\to\mathcal{O}_{C}(-1)$, then $Y$ is not homogeneous.
\end{lemma} 
\begin{proof}
A transitive holomorphic group action supplies global vector fields
spanning every tangent space, so $TY$ is globally generated. Restriction
and passage to a quotient preserve global generation, whereas $\mathcal{O}_{\mathbb{P}^{1}}(-1)$
has no nonzero holomorphic section. 
\end{proof}

In Section~\ref{sec:5} we use local potentials $g_{a\bar{b}}=\partial_{a}\partial_{\bar{b}}\Phi$;
potentials differing by the real part of a holomorphic function define
the same metric. The potentials are shown to extend smoothly with
a positive metric across the complement of the coordinate charts.
A uniform curvature bound established on a dense open set then extends
by continuity of the eigenvalues of Hermitian matrices.

\section{Positive holomorphic sectional curvature}

\label{sec:3}

We prove Theorem~\ref{theorem:1.1} for an arbitrary invariant Kähler
metric. The proof does not require the Einstein condition. The main
step is to express the holomorphic sectional curvature as an integral
of squared norms. The expression retains the isotropy part of every
bracket, which is essential for strict positivity.

\subsection{The invariant connection and the curvature numerator}

Recall the operator $D=\operatorname{ad}\xi$, the adjoint $U^{*}=-\sigma(U)$,
and the inner product $\langle U,V\rangle_{B}=B(U,V^{*})$ from Section~\ref{sec:2}.
Thus $D$ is positive on $\mathfrak{m}_{+}$, zero on $\mathfrak{k}^{\mathbb{C}}$,
and negative on $\mathfrak{m}_{-}$. The metric is given by $h(U,V)=-B(|D|U,V)$,
and $|X|^{2}_{g}=B(DX,X^{*})$ for $X\in\mathfrak{m}_{+}$. Inverse
powers of $D$ below are taken on $\mathfrak{m}_{+}$; the inverse
of $|D|$ is taken on $\mathfrak{m}^{\mathbb{C}}$. The eigenvalues
need not be integers.

For clarity, the adjoint relations used in the calculation are 
\begin{align}
B(DU,V) & =-B(U,DV), & B(|D|U,V) & =B(U,|D|V),\label{eq:3.1}\\
(DU)^{*} & =-D(U^{*}), & (|D|U)^{*} & =|D|(U^{*}),\nonumber \\
|D|^{2} & =D^{2}.\nonumber 
\end{align}
For example, the first identity follows from invariance of $B$. If
$DU=sU$, then $D(U^{*})=-sU^{*}$; the last two identities follow
by applying this observation to each eigenspace.

We shall use the invariant connection in an operator form. The Levi-Civita
connection at the origin is encoded by its Nomizu map $\Lambda$,
with 
\begin{equation}
\Lambda_{U}V=\frac{1}{2}[U,V]_{\mathfrak{m}}+\frac{1}{2}|D|^{-1}\bigl([U,|D|V]+[V,|D|U]\bigr)_{\mathfrak{m}}.\label{eq:3.2}
\end{equation}
The inverse is taken only on $\mathfrak{m}^{\mathbb{C}}$; the subscripts
$\mathfrak{m}$ in \eqref{eq:3.2} denote the complex-linear extension
of the projection onto the real tangent complement. To check the formula
directly, insert $h=-B(|D|\cdot,\cdot)$ into the homogeneous Koszul
identity 
\[
2h(\Lambda_{U}V,W)=h([U,V]_{\mathfrak{m}},W)+h([W,U]_{\mathfrak{m}},V)+h(U,[W,V]_{\mathfrak{m}})
\]
and use $B([U,V],W)=B(U,[V,W])$ in the last two terms. This gives
\eqref{eq:3.2}; see also \cite{Nom54}. Two consequences of this
formula will be used repeatedly: 
\begin{equation}
\Lambda_{U}V=D^{-1}[U,DV],\quad\Lambda_{U^{*}}V=[U^{*},V]_{\mathfrak{m}_{+}}\quad\text{for }U,V\in\mathfrak{m}_{+}.\label{eq:3.3}
\end{equation}
It suffices to verify these identities when $U,V$ lie in the $s$-
and $t$-eigenspaces of $D$, with $s,t>0$. The bracket $[U,V]$
lies in the $(s+t)$-eigenspace, so its coefficient in \eqref{eq:3.2}
is 
\[
\frac{1}{2}\left(1+\frac{t-s}{s+t}\right)=\frac{t}{s+t}.
\]
For $[U^{*},V]$, which lies in the $(t-s)$-eigenspace, the coefficient
is 
\[
\frac{1}{2}\left(1+\frac{t-s}{|t-s|}\right)\quad\text{for }s\ne t.
\]
It is one when $t>s$ and zero when $t<s$. If $s=t$, the bracket
belongs to the isotropy algebra and its tangent projection vanishes.
This proves \eqref{eq:3.3}, including the case of equal eigenvalues.

The curvature endomorphism is given by \eqref{eq:2.7}. Since $\overline{X}=\sigma(X)=-X^{*}$,
the curvature convention fixed in Section~\ref{sec:2} gives 
\begin{equation}
R(X,\overline{X},X,\overline{X})=-B(D\mathcal{R}(X,X^{*})X,X^{*}).\label{eq:3.4}
\end{equation}

\begin{lemma}\label{lemma:3.1} For $X\in\mathfrak{m}_{+}$, set
\[
C=[X,X^{*}],\quad K=[DX,X^{*}],\quad F=[X,DX].
\]
Then 
\[
R(X,\overline{X},X,\overline{X})=\|D^{-1/2}F\|^{2}_{B}+B(C,K)-B(C,|D|C).
\]
\end{lemma} 
\begin{proof}
Decompose $C=C_{-}+C_{0}+C_{+}$ according to the spaces $\mathfrak{m}_{-}$
, $\mathfrak{k}^{\mathbb{C}}$, and $\mathfrak{m}_{+}$ corresponding
to negative, zero, and positive eigenspaces of $D$. The identity
$C^{*}=C$ gives $C_{-}=C^{*}_{+}$. Moreover, the derivation rule
for $D$ gives 
\begin{equation}
DC=[DX,X^{*}]-[X,(DX)^{*}]=K-K^{*}.\label{eq:3.5}
\end{equation}
Inserting \eqref{eq:3.3} into \eqref{eq:2.7} yields 
\begin{align*}
\mathcal{R}(X,X^{*})X={} & -D^{-1}[X,DC_{+}]-[X^{*},D^{-1}F]_{\mathfrak{m}_{+}}\\
 & -D^{-1}[C_{+},DX]-[C_{-},X]_{\mathfrak{m}_{+}}-[C_{0},X].
\end{align*}
Indeed, $\Lambda_{X^{*}}X=-C_{+}$ and $\Lambda_{X}X=D^{-1}F$ give
the first two terms. The tangent and isotropy components of $[X,X^{*}]$
give the remaining three terms.

After applying \eqref{eq:3.4}, invariance of $B$ gives 
\begin{align*}
R(X,\overline{X},X,\overline{X})={} & \|D^{-1/2}F\|^{2}_{B}-B(DC_{+},C_{-})\\
 & +B(C_{+},K)+B(C_{-},K^{*})+B(C_{0},K^{*}).
\end{align*}
For example, the first connection term contributes 
\[
B([X,DC_{+}],X^{*})=-B(DC_{+},C_{-}),
\]
whereas the second contributes 
\[
B([X^{*},D^{-1}F],(DX)^{*})=B(D^{-1}F,[(DX)^{*},X^{*}])=\|D^{-1/2}F\|^{2}_{B}.
\]
$B$ pairs the $s$- and $t$-eigenspaces of $D$ only when $s+t=0$.
Consequently, 
\begin{align*}
B(C_{0},DC) & =0,\\
B(C_{-},DC) & =B(C_{-},DC_{+}),\\
B(C,|D|C) & =2B(C_{-},DC_{+}).
\end{align*}
Replacing $K^{*}$ by $K-DC$ using \eqref{eq:3.5} therefore combines
the last three terms with $-B(DC_{+},C_{-})$ into $B(C,K)-B(C,|D|C)$.
This proves the lemma. 
\end{proof}

\subsection{An integral of squared norms}

For $X\in\mathfrak{m}_{+}$ and $t\geq0$, define 
\begin{equation}
X_{t}=e^{-\frac{tD}{2}}X,\quad C_{t}=[X_{t},X^{*}_{t}],\quad M_{t}=-C_{t}'-\frac{1}{2}|D|C_{t}.\label{eq:3.6}
\end{equation}
The exponential multiplies the component of $X$ in the $s$-eigenspace
by $e^{-\frac{ts}{2}}$. Since $D$ is positive on $\mathfrak{m}_{+}$,
$X_{t}$ and all the bracket expressions below tend exponentially
to zero as $t\longrightarrow\infty$.

\begin{proposition}\label{proposition:3.2} For every $X\in\mathfrak{m}_{+}$,
\begin{equation}
R(X,\overline{X},X,\overline{X})=2\int^{\infty}_{0}\|M_{t}\|^{2}_{B}\,dt.\label{eq:3.7}
\end{equation}
The term $M_{t}$ is Hermitian, and it is nonzero for every finite
$t\geq0$ whenever $X\ne0$. \end{proposition} 
\begin{proof}
Put 
\[
K_{t}=[DX_{t},X^{*}_{t}],\quad F_{t}=[X_{t},DX_{t}],\quad L_{t}=[DX_{t},(DX_{t})^{*}].
\]
Differentiating $X_{t}'=-\frac{1}{2}DX_{t}$ gives 
\begin{align}
C_{t}' & =-\frac{1}{2}(K_{t}+K^{*}_{t}),\label{eq:3.8}\\
DC_{t} & =K_{t}-K^{*}_{t},\nonumber \\
C_{t}'' & =\frac{1}{4}D^{2}C_{t}+L_{t}.\nonumber 
\end{align}
For the third identity, direct differentiation first gives 
\[
C_{t}''=\frac{1}{4}[D^{2}X_{t},X^{*}_{t}]+\frac{1}{2}[DX_{t},(DX_{t})^{*}]+\frac{1}{4}[X_{t},(D^{2}X_{t})^{*}].
\]
On the other hand, 
\[
D^{2}C_{t}=[D^{2}X_{t},X^{*}_{t}]-2[DX_{t},(DX_{t})^{*}]+[X_{t},(D^{2}X_{t})^{*}]
\]
showing that 
\[
C_{t}''-\frac{1}{4}D^{2}C_{t}=[DX_{t},(DX_{t})^{*}]=L_{t}
\]

The crucial Jacobi identity responsible for the cancellation in the
proof is 
\begin{equation}
\|F_{t}\|^{2}_{B}=B(K_{t},K^{*}_{t})-B(C_{t},L_{t}).\label{eq:3.9}
\end{equation}
Temporarily write $Y=DX_{t}$ and $X=X_{t}$, so $F=[X,Y]$ and $K=[Y,X^{*}]$
and note that 
\begin{align*}
B(K,K^{*}) & =B([[Y,X^{*}],X],Y^{*})\\
 & =B(-[Y,[X,X^{*}]]+[X^{*},[X,Y]],Y^{*})\\
 & =B([X,X^{*}],[Y,Y^{*}])+B([X,Y],[Y^{*},X^{*}])\\
 & =B(C,L)+\|F\|^{2}_{B}.
\end{align*}
This proves \eqref{eq:3.9}.

The first two identities in \eqref{eq:3.8} express $K_{t}=-C_{t}'+\frac{1}{2}DC_{t}$
and $K^{*}_{t}=-C_{t}'-\frac{1}{2}DC_{t}$. Using symmetry of $B$
and \eqref{eq:3.1}, we obtain 
\[
B(K_{t},K^{*}_{t})=B(C_{t}',C_{t}')+\frac{1}{4}B(C_{t},D^{2}C_{t}).
\]
The last identity in \eqref{eq:3.8} similarly gives 
\[
B(C_{t},L_{t})=B(C_{t},C_{t}'')-\frac{1}{4}B(C_{t},D^{2}C_{t}).
\]
Substitution in \eqref{eq:3.9} yields 
\begin{equation}
\|F_{t}\|^{2}_{B}=B(C_{t}',C_{t}')-B(C_{t},C_{t}'')+\frac{1}{2}B(C_{t},D^{2}C_{t}).\label{eq:3.10}
\end{equation}
Because $D$ is a derivation, damping commutes with brackets in $\mathfrak{m}_{+}$.
Thus $F_{t}=e^{-\frac{tD}{2}}F(0)$. Integration on each eigenspace
gives 
\[
\int^{\infty}_{0}\|F_{t}\|^{2}_{B}\,dt=\|D^{-1/2}F(0)\|^{2}_{B}.
\]
Integration by parts in \eqref{eq:3.10} therefore gives 
\begin{align}
\|D^{-1/2}F(0)\|^{2}_{B}={} & 2\int^{\infty}_{0}B(C_{t}',C_{t}')\,dt+B(C(0),C'(0))\nonumber \\
 & +\frac{1}{2}\int^{\infty}_{0}B(C_{t},D^{2}C_{t})\,dt.\label{eq:3.11}
\end{align}
Here $C(0)$, $F(0)$, and $K(0)$ denote values at time zero. The
boundary term at infinity is zero by exponential decay.

Since $C_{t}$ is Hermitian and $|D|$ commutes with $*$, both $C_{t}'$
and $|D|C_{t}$ are Hermitian. We may therefore expand their squared
norm using $B$. The relations $|D|^{2}=D^{2}$ and $2B(C_{t}',|D|C_{t})=\frac{d}{dt}B(C_{t},|D|C_{t})$
give 
\begin{align*}
2\int^{\infty}_{0}\left\Vert C_{t}'+\frac{1}{2}|D|C_{t}\right\Vert ^{2}_{B}\,dt & =2\int^{\infty}_{0}B(C_{t}',C_{t}')\,dt-B(C(0),|D|C(0))\\
 & \quad+\frac{1}{2}\int^{\infty}_{0}B(C_{t},D^{2}C_{t})\,dt\\
 & =\|D^{-1/2}F(0)\|^{2}_{B}-B(C(0),C'(0))-B(C(0),|D|C(0)),
\end{align*}
where the last line uses \eqref{eq:3.11}. Finally, $B(C(0),DC(0))=0$
by \eqref{eq:3.1}. Hence \eqref{eq:3.8} implies $B(C(0),C'(0))=-B(C(0),K(0))$.
Lemma~\ref{lemma:3.1} proves \eqref{eq:3.7}.

It remains to prove strict positivity. Since $|D|\xi=0$, invariance
of $B$ gives 
\begin{align*}
B(\xi,M_{t}) & =-\frac{d}{dt}B(\xi,[X_{t},X^{*}_{t}])\\
 & =-\frac{d}{dt}B(DX_{t},X^{*}_{t})\\
 & =B(DX_{t},(DX_{t})^{*})>0\quad\text{for }X\ne0.
\end{align*}
In the last step, $D$ is self-adjoint for the Hermitian inner product
and is invertible on $\mathfrak{m}_{+}$. Thus $M_{t}\ne0$. Because
$M_{t}$ is Hermitian, its squared norm is strictly positive. 
\end{proof}

\begin{proof}[Proof of Theorem~\ref{theorem:1.1}]
Proposition~\ref{proposition:3.2} proves $R(X,\overline{X},X,\overline{X})>0$
for every nonzero $X\in T^{1,0}_{o}M$. Dividing by $|X|^{4}_{g}$
and using homogeneity proves the theorem at every point of $M$. The
argument applies to semisimple groups, so it also includes products
of flag manifolds. In particular, it applies to the canonical Kähler-Einstein
metric in Conjecture 1.8 of \cite{NZ25}. 
\end{proof}

\subsection{A lower bound and the root interpretation}

The same identity gives a uniform positive lower bound for each fixed
invariant metric.

\begin{corollary}\label{corollary:3.3} Let $m$ be the smallest
eigenvalue of $D$ on $\mathfrak{m}_{+}$. Then 
\[
H(X)\geq\frac{m}{B(\xi,\xi)}>0\quad\text{for every }X\ne0.
\]
\end{corollary} 
\begin{proof}
Set $f(t)=|X_{t}|^{2}_{g}=B(DX_{t},X^{*}_{t})$. The preceding proof
and the inequality $D\geq m\operatorname{Id}$ on $\mathfrak{m}_{+}$
give 
\[
B(\xi,M_{t})=-f'(t)=\|DX_{t}\|^{2}_{B}\geq mf(t).
\]
Cauchy-Schwarz on the real space of Hermitian elements and \eqref{eq:3.7}
now imply 
\begin{align*}
R(X,\overline{X},X,\overline{X}) & \geq\frac{2}{B(\xi,\xi)}\int^{\infty}_{0}(f'(t))^{2}\,dt\\
 & \geq-\frac{2m}{B(\xi,\xi)}\int^{\infty}_{0}f(t)f'(t)\,dt\\
 & =\frac{m}{B(\xi,\xi)}|X|^{4}_{g}.
\end{align*}
The second inequality uses $-f'\geq mf\geq0$, and the last equality
uses $f(0)=|X|^{2}_{g}$ and $f(t)\longrightarrow0$. 
\end{proof}

The integral has a useful expression in the root vectors already fixed
in Section~\ref{sec:2}. Write 
\[
X=\sum_{\alpha\in\Delta^{+}_{\Phi_{T}}}x_{\alpha}E_{\alpha},\quad g_{\alpha}=\alpha(\xi)>0.
\]
The bracket $[E_{\alpha},E_{-\beta}]$ lies in the $(g_{\alpha}-g_{\beta})$-eigenspace
of $D$. Formula \eqref{eq:3.6} therefore becomes 
\begin{equation}
M_{t}=\sum_{\alpha,\beta\in\Delta^{+}_{\Phi_{T}}}\min(g_{\alpha},g_{\beta})e^{-\frac{t(g_{\alpha}+g_{\beta})}{2}}x_{\alpha}\overline{x_{\beta}}[E_{\alpha},E_{-\beta}].\label{eq:3.12}
\end{equation}
Indeed, the coefficient obtained by differentiating is $\frac{1}{2}(g_{\alpha}+g_{\beta})$,
and the term $\frac{1}{2}|D|C_{t}$ subtracts $\frac{1}{2}|g_{\alpha}-g_{\beta}|$.
Their difference is the minimum of the two metric coefficients. Formula
\eqref{eq:3.12} includes brackets in the Cartan subalgebra and in
the semisimple part of the isotropy algebra. No term is discarded
when two metric coefficients coincide. 
\begin{example}
\label{example:3.4} Consider $\mathrm{SU}(3)/T^{2}$, with positive
root vectors $E_{\alpha_{1}}=E_{12}$, $E_{\alpha_{2}}=E_{23}$, and
$E_{\alpha_{1}+\alpha_{2}}=E_{13}=[E_{12},E_{23}]$, where $B(U,V)=\operatorname{tr}(UV)$.
Their metric coefficients are $g_{\alpha_{1}},g_{\alpha_{2}},g_{\alpha_{1}}+g_{\alpha_{2}}$,
with $g_{\alpha_{1}},g_{\alpha_{2}}>0$. For $X=xE_{12}+yE_{23}$,
the mixed brackets $[E_{12},E_{32}]$ and $[E_{23},E_{21}]$ vanish.
Thus 
\[
M_{t}=g_{\alpha_{1}}e^{-g_{\alpha_{1}}t}|x|^{2}H_{12}+g_{\alpha_{2}}e^{-g_{\alpha_{2}}t}|y|^{2}H_{23},
\]
where $H_{12}=\operatorname{diag}(1,-1,0)$ and $H_{23}=\operatorname{diag}(0,1,-1)$.
Since $B(H_{12},H_{12})=B(H_{23},H_{23})=2$ and $B(H_{12},H_{23})=-1$,
the integral gives 
\[
R(X,\overline{X},X,\overline{X})=2g_{\alpha_{1}}|x|^{4}+2g_{\alpha_{2}}|y|^{4}-\frac{4g_{\alpha_{1}}g_{\alpha_{2}}}{g_{\alpha_{1}}+g_{\alpha_{2}}}|x|^{2}|y|^{2}.
\]
In particular, a negative mixed coefficient is compatible with strictly
positive holomorphic sectional curvature. If $g_{\alpha_{1}}=g_{\alpha_{2}}=1$,
the canonical normalization is $\Ric=2g$. In this case $\xi=\operatorname{diag}(1,0,-1)$,
$m=1$, and $B(\xi,\xi)=2$. Corollary~\ref{corollary:3.3} shows
that $H\geq\frac{1}{2}$ with the lower bound attained for $\frac{1}{\sqrt{2}}(E_{12}+E_{23})$. 
\end{example}

\section{The spectrum of the Calabi operator}

\label{sec:4}

We prove Theorem~\ref{theorem:1.2} and its sharper upper bound.
The Calabi operator of the canonical metric, normalized by $\Ric=2g$,
satisfies 
\[
-2\operatorname{Id}<\mathcal{C}\leq\frac{4}{k_{*}}\operatorname{Id}\leq2\operatorname{Id}.
\]
The upper bound follows from an isometric embedding into a product
of projective spaces and is sharp in every Lie type. The lower bound
is proved directly on the root-sum blocks of $\mathcal{C}$: a root
count suffices in the simply-laced types, scaled matrix estimates
treat types $B$, $C$, and $F_{4}$, and explicit blocks settle $G_{2}$.
The curvature formulas recalled from \cite[Proposition~16]{Loh19}
in Section~\ref{sec:2} include roots in the Levi factor, so all
painted diagrams are covered.

We also identify the eigenspace of eigenvalue $2$. It is nontrivial
precisely when the flag admits a homogeneous projection with fiber
$\mathbb{P}^{1}$. The quadratic-curvature consequences follow immediately
from \eqref{eq:2.3}.

\subsection{The sharp upper bound in every Lie type}

Recall the anticanonical degrees $k_{\alpha}$ and their minimum $k_{*}$
introduced after Theorem~\ref{theorem:1.2}. The $k_{\alpha}$ are
also called Koszul coefficients.

\begin{theorem}\label{theorem:4.1} For every generalized flag manifold
of positive dimension with its canonical metric normalized by $\Ric=2g$,
\[
\lambda_{\max}(\mathcal{C})=\frac{4}{k_{*}}\leq2.
\]
In particular, the upper bound is strict precisely when every $k_{\alpha}>2$.
\end{theorem} 
\begin{proof}
For $\alpha\in\Phi_{T}$, since $2(\rho,\alpha)/|\alpha|^{2}=1$,
\[
k_{\alpha}=\frac{2(\delta,\alpha)}{|\alpha|^{2}}=\frac{2(2\rho-2\rho_{\Phi_{K}},\alpha)}{|\alpha|^{2}}=2-\frac{4(\rho_{\Phi_{K}},\alpha)}{|\alpha|^{2}}\geq2.
\]
Indeed, every positive Levi root is a nonnegative integral combination
of simple roots other than $\alpha$, whose inner products with $\alpha$
are nonpositive.

Use the fundamental projective embedding $F$ from Section~\ref{sec:2},
with the factor indexed by $\alpha$ carrying the metric $\frac{k_{\alpha}}{2}g_{\mathrm{FS},\alpha}$.
Its induced metric is the canonical metric with $\operatorname{Ric}=2g$.
The invariant metrics and their representation-theoretic potentials
are described in \cite[Sections~2--5]{AB03}.

For a holomorphic isometric immersion, the Kähler-Gauss equation is
\[
R_{a\bar{b}c\bar{d}}=\overline{R}_{a\bar{b}c\bar{d}}-\langle\mathrm{II}(e_{a},e_{c}),\mathrm{II}(e_{b},e_{d})\rangle,
\]
where $\overline{R}$ denotes ambient curvature and $\mathrm{II}$
the complex second fundamental form. On the symmetric tensors this
becomes 
\[
\langle\mathcal{C}S,S\rangle=\langle\overline{\mathcal{C}}S,S\rangle-\left\Vert \sum_{a,c}\mathrm{II}(e_{a},e_{c})S_{ac}\right\Vert ^{2}.
\]
An ambient upper spectral bound for the Calabi curvature operator
therefore passes to the submanifold. The ambient projective factors
have Calabi eigenvalues $4/k_{\alpha}$, and mixed-factor symmetric
tensors have eigenvalue zero. This accounts for the upper bound $4/k_{*}$
in Section~\ref{sec:4}; the factor 4 corresponds to $\Ric=2g$.
In particular, we have 
\[
\langle\mathcal{C}S,S\rangle=\langle\overline{\mathcal{C}}S,S\rangle-\left\Vert \sum_{a,c}\mathrm{II}(e_{a},e_{c})S_{ac}\right\Vert ^{2}\leq\frac{4}{k_{*}}\|S\|^{2}.
\]
Thus $\mathcal{C}\leq(4/k_{*})\operatorname{Id}$.

To see sharpness, take a painted simple root $\alpha\in\Phi_{T}$
with $k_{\alpha}=k_{*}$. Its root $\mathrm{SL}_{2}\subset G^{\mathbb{C}}$
orbit is a simple Schubert curve. It maps to a projective line in
the factor indexed by $\alpha$ and to a point in every other factor,
as noted in Section~\ref{sec:2}. The curve is totally geodesic in
the ambient product. Its ambient second fundamental form is the orthogonal
sum of its second fundamental form in $M$ and the restriction of
$\mathrm{II}$. Both therefore vanish along the curve. Its induced
metric is $\frac{k_{\alpha}}{2}g_{\mathrm{FS},\alpha}$, so its unit
tangent $e_{\alpha}$ satisfies $H(e_{\alpha})=\frac{4}{k_{\alpha}}$
in $M$. Since this attains the upper bound, it follows that $\mathcal{C}(e_{\alpha}\otimes e_{\alpha})=\frac{4}{k_{\alpha}}e_{\alpha}\otimes e_{\alpha}$.
In particular, the eigenvalue is $2$ when $k_{\alpha}=2$. 
\end{proof}

\begin{proposition}[Equality space]\label{proposition:4.2} The
2-eigenvalue space of the Calabi operator consists precisely of the
symmetric root squares 
\[
\ker(2\operatorname{Id}-\mathcal{C})=\bigoplus_{\substack{\alpha\in\Phi_{T}\\
k_{\alpha}=2
}
}\mathbb{C}(e_{\alpha}\otimes e_{\alpha}).
\]
In particular, a root-sum block of $\mathcal{C}$ has largest eigenvalue
strictly less than two unless it is indexed by $2\alpha$ for such
a painted simple root. \end{proposition} 
\begin{proof}
We already saw in the proof of Theorem~\ref{theorem:4.1} that 
\[
\bigoplus_{\substack{\alpha\in\Phi_{T}\\
k_{\alpha}=2
}
}\mathbb{C}(e_{\alpha}\otimes e_{\alpha})\subset\ker(2\operatorname{Id}-\mathcal{C}).
\]
Conversely, suppose $\mathcal{C}S=2S$. Identify $T^{1,0}_{o}M$ with
its image under $DF$. Equality in the ambient upper bound requires
$S$ to have components only in the pure-factor symmetric squares
for which $k_{\alpha}=2$. For such an $\alpha$, the root-space description
of the maps $F_{\beta}$ gives 
\[
DF(T^{1,0}_{o}M)\cap T^{1,0}_{F_{\alpha}(o)}\mathbb{P}(V_{\alpha})=\mathbb{C}e_{\alpha}.
\]
Indeed, this intersection is tangent to the fiber of the projection
obtained by forgetting $F_{\alpha}$. The equality $k_{\alpha}=2$
means that $\alpha$ is not adjacent to any Levi node; adding it to
$\Phi_{K}$ therefore adds only an isolated $A_{1}$ component, so
that fiber has tangent space $\mathbb{C}e_{\alpha}$. Contracting
the pure-factor component of $S$ with a covector on that factor,
extended by zero on the other factors, yields a vector in this intersection.
Hence the component lies in $\mathbb{C}(e_{\alpha}\otimes e_{\alpha})$.
This proves the opposite inclusion. 
\end{proof}

\begin{corollary}\label{corollary:4.3} The canonical metric satisfies
$\lambda_{\max}(\mathcal{C})<2$ exactly when every painted node is
adjacent to at least one white node. Equivalently, there is no homogeneous
projection from $M$ with fiber $\mathbb{P}^{1}$. By \eqref{eq:2.3},
these are also the equivalent conditions for $\CQB>0$. \end{corollary} 
\begin{proof}
The coefficient $k_{\alpha}$ equals two precisely when $\alpha$
is not adjacent to $\Phi_{K}$. This follows because $\rho_{\Phi_{K}}$
has strictly positive simple-root coefficients on every nonempty Levi
component, and all cross-Cartan pairings have the same nonpositive
sign. Adding such an $\alpha$ to $\Phi_{K}$ adds an isolated $A_{1}$
Levi component, giving a homogeneous projection with fiber $\mathrm{SL}_{2}/B\cong\mathbb{P}^{1}$.
Conversely a one-dimensional homogeneous parabolic fiber has precisely
this form. Apply Theorem~\ref{theorem:4.1}. 
\end{proof}

\subsection{The lower spectral bound in simply-laced types}

\label{subsec:4.2}

Normalize all roots to squared length 2 and put 
\[
G_{*}=\max_{\alpha\in\Delta^{+}_{\Phi_{T}}}g_{\alpha}.
\]
In a root basis with $B(E_{\alpha},E_{-\alpha})=1$, the squared norm
of $E_{\alpha}$ is $g_{\alpha}$. Each positive root has a nonnegative
metric coefficient, which vanishes precisely on the Levi roots. On
a painted simple root $\alpha\in\Phi_{T}$, $g_{\alpha}=\frac{k_{\alpha}}{2}\geq1$,
so $g_{\alpha}\geq1$ for every $\alpha\in\Delta^{+}_{\Phi_{T}}$.

\begin{theorem}\label{theorem:4.4} For every positive-dimensional
generalized flag of simply-laced type with $\Ric=2g$, 
\[
\lambda_{\min}(\mathcal{C})\geq-2+\frac{2}{G_{*}}>-2.
\]
Together with Theorem~\ref{theorem:4.1}, this proves Theorem~\ref{theorem:1.2}
in types $A$, $D$, and $E$. \end{theorem}

\begin{lemma}[Root count]\label{lemma:4.5} For $\tau\in\Delta^{+}_{\Phi_{T}}$,
let $m_{\tau}$ count unordered decompositions $\tau=\alpha+\gamma$
with $\alpha,\gamma\in\Delta^{+}_{\Phi_{T}}$, and let $h_{\tau}$
count $\beta\in\Delta^{+}_{\Phi_{T}}$ for which $(\tau,\beta)=1$
and $\tau-\beta$ is a Levi root, of either sign. Then 
\[
g_{\tau}=1+m_{\tau}+\frac{1}{2}h_{\tau}.
\]
\end{lemma} 
\begin{proof}
In $2g_{\tau}=\sum_{\beta\in\Delta^{+}_{\Phi_{T}}}(\tau,\beta)$ the
term $\beta=\tau$ contributes 2. Every root with pairing $-1$ pairs
with $\beta+\tau\in\Delta^{+}_{\Phi_{T}}$, whose pairing is $+1$,
and these cancel. Each remaining $+1$ root belongs either to a tangent
decomposition of $\tau$, contributing 2 per unordered pair, or to
the indicated Levi-difference set, contributing 1. A root with positive
metric coefficient is positive and tangent, which justifies this exhaustive
pairing. 
\end{proof}

\begin{proof}[Proof of Theorem~\ref{theorem:4.4}]
We use the root-sum blocks and unit pair basis of \eqref{eq:2.5}.
The possible squared norms of $\tau$ are $8,6,4,2$; these exhaust
the cases because the roots are simply laced and the summands are
positive.

If $|\tau|^{2}=8$, the block consists of the repeated pair $\tau=\alpha+\alpha$
and $\mathcal{C}_{\tau}=\frac{2}{g_{\alpha}}>0$.

If $|\tau|^{2}=6$, the two roots $\alpha,\gamma$ have pairing 1,
the unordered pair is unique, and $\mathcal{C}_{\tau}=\frac{2}{\max(g_{\alpha},g_{\gamma})}>0$.
Uniqueness follows since another summand $\beta$ would satisfy $(\tau,\beta)=3$,
forcing one pairing with a summand to be 2 and hence forcing equality
with that summand. \smallskip{}
 When $|\tau|^{2}=2$, the sum $\tau$ is itself a tangent root. Only
sums that index a nonempty block need to be considered, so $m_{\tau}>0$.
Every diagonal block entry is $-\frac{2}{g_{\tau}}$, and every off-diagonal
entry has absolute value at most $\frac{2}{g_{\tau}}$.

To verify the off-diagonal bound, consider distinct decompositions
$\alpha+\gamma=\beta+\epsilon=\tau$. Since the decompositions are
distinct, all cross-pairings are integers $\leq1$. Moreover, $(\tau,\beta)=(\tau,\epsilon)=1$,
so each column sum of the $2\times2$ matrix of cross-pairings is
1; its row sums are also 1. Its entries must therefore be 0 or 1,
with one of each in every row and column. Exactly one of the bracket
pairings $N_{1},N_{2}$ in Proposition~\ref{proposition:2.1} is
nonzero, and its absolute value is 1: the corresponding differences
are roots, while the other two differences are not. Since $z_{\alpha}=1$
in this normalization, \eqref{eq:2.5} and \eqref{eq:2.6} give, in
the first case, 
\[
|\mathcal{C}_{\alpha\gamma,\beta\epsilon}|=\frac{2}{g_{\tau}}\frac{\min(g_{\alpha},g_{\beta})}{\sqrt{g_{\alpha}g_{\beta}}}\frac{\min(g_{\gamma},g_{\epsilon})}{\sqrt{g_{\gamma}g_{\epsilon}}}\leq\frac{2}{g_{\tau}}.
\]
The second case is identical after interchanging $\beta$ and $\epsilon$.
The formula also covers root differences in the Levi system. Consequently
the sum of the absolute values in a row is at most $\frac{2m_{\tau}}{g_{\tau}}$.
Lemma~\ref{lemma:4.5} gives 
\[
\lambda_{\min}(\mathcal{C}_{\tau})\geq-\frac{2m_{\tau}}{g_{\tau}}\geq-2+\frac{2}{g_{\tau}}\geq-2+\frac{2}{G_{*}}.
\]
\smallskip{}
 When $|\tau|^{2}=4$ every pair $\alpha_{j}+\gamma_{j}=\tau$ is
orthogonal. Let $m_{\tau}$ be the number of such unordered decompositions.
When $m_{\tau}=1$ the block vanishes. For $m_{\tau}\geq2$, orient
each pair so that $g_{\alpha_{j}}\leq g_{\gamma_{j}}$; the sum $s=g_{\alpha_{j}}+g_{\gamma_{j}}$
is fixed.

All four cross-pairings between different pairs equal 1: they are
each at most 1, and each pair of them sums to 2. Consider $\frac{\alpha_{j}+\gamma_{j}}{2}$,
$\frac{\alpha_{j}-\gamma_{j}}{2}$. These are mutually orthogonal
unit vectors, and the roots are $\alpha_{j}$, $\gamma_{j}$, together
with $\frac{\alpha_{j}-\gamma_{j}}{2}\pm\frac{\alpha_{k}-\gamma_{k}}{2}$.
They generate a subsystem of type $D_{m_{\tau}+1}$ (with $D_{3}=A_{3}$).
In its usual orthogonal matrix realization, the difference-bracket
products are $-1$. Choosing these normalized root vectors simultaneously
for the whole subsystem gives one compatible phase choice for all
pairs. In the corresponding unit symmetric-pair basis, the block has
zero diagonal and off-diagonal entries 
\[
(\mathcal{C}_{\tau})_{jk}=-\frac{2\min(g_{\alpha_{j}},g_{\alpha_{k}})}{\sqrt{g_{\alpha_{j}}g_{\gamma_{j}}g_{\alpha_{k}}g_{\gamma_{k}}}},\quad\text{for }j\ne k.
\]
This computation is intrinsic to the root subsystem and applies also
when it lies in an exceptional root system. To verify the sign, identify
the roots with $e_{0}-e_{j},e_{0}+e_{j}$ and write $F_{ab}$ for
elementary matrices, with distinct formal indices $0,1,\ldots,m_{\tau},\bar{0},\bar{1},\ldots,\bar{m_{\tau}}$.
Use 
\[
E_{e_{0}-e_{j}}=F_{0j}-F_{\bar{j},\bar{0}},\qquad E_{e_{0}+e_{j}}=F_{0,\bar{j}}-F_{j,\bar{0}};
\]
their difference-bracket products are negative. The root $2e_{0}$
does not exist, so no sum-root term enters the curvature formula.

For this fixed root sum, the matrix in the same ordered pair basis
\[
K_{jk}=\frac{\min(g_{\alpha_{j}},g_{\alpha_{k}})}{\sqrt{g_{\alpha_{j}}g_{\gamma_{j}}g_{\alpha_{k}}g_{\gamma_{k}}}}
\]
is positive semidefinite. To see this without computing its eigenvalues,
use 
\[
\min(g_{\alpha_{j}},g_{\alpha_{k}})=\int^{\infty}_{0}\mathbf{1}_{[0,g_{\alpha_{j}}]}(t)\mathbf{1}_{[0,g_{\alpha_{k}}]}(t)\,dt.
\]
For every vector $(c_{j})$, its quadratic form is therefore the integral
of 
\[
\left|\sum_{j}\frac{c_{j}}{\sqrt{g_{\alpha_{j}}g_{\gamma_{j}}}}\mathbf{1}_{[0,g_{\alpha_{j}}]}(t)\right|^{2},
\]
which is nonnegative. Hence 
\begin{equation}
\mathcal{C}_{\tau}=2\operatorname{diag}(\frac{1}{g_{\gamma_{j}}})-2K.\label{eq:4.1}
\end{equation}
We claim $g_{\gamma_{j}}\geq m_{\tau}$ for every $j$. Apply Lemma~\ref{lemma:4.5}
to $\gamma_{j}$. For each $k\ne j$, $(\gamma_{j},\alpha_{k})=1$
and $g_{\alpha_{k}}\leq\frac{s}{2}\leq g_{\gamma_{j}}$. If $g_{\alpha_{k}}<g_{\gamma_{j}}$,
then $\alpha_{k}+(\gamma_{j}-\alpha_{k})=\gamma_{j}$ is a tangent
decomposition. Such decompositions are distinct: the chosen roots
have pairing 2 with $\tau$, whereas their complements have pairing
0. If $g_{\alpha_{k}}=g_{\gamma_{j}}$, all four relevant metric coefficients
equal $\frac{s}{2}$, and both $\alpha_{k}$ and $\gamma_{k}$ contribute
to $h_{\gamma_{j}}$, for a total of 1 in the root-count formula.
These contributions are distinct and do not overlap the earlier ones.
Thus the $m_{\tau}-1$ other indices contribute at least $m_{\tau}-1$,
in addition to the original 1.

It follows that $\operatorname{tr}K=\sum_{j}\frac{1}{g_{\gamma_{j}}}\leq1$,
so $K\leq\operatorname{Id}$. Equation \eqref{eq:4.1} now gives 
\[
\mathcal{C}_{\tau}\geq-2\operatorname{Id}+2\operatorname{diag}\left(\frac{1}{g_{\gamma_{j}}}\right)\geq\left(-2+\frac{2}{G_{*}}\right)\operatorname{Id}.
\]
Since $G_{*}\geq1$, all blocks satisfy the same bound. 
\end{proof}

\begin{remark}\label{remark:4.6} For the flag of type $SU(3)/T^{2}$,
the canonical metric coefficients of $\alpha,\beta,\alpha+\beta$
are $1,1,2$. The pair $(\alpha,\beta)$ gives the scalar block $\mathcal{C}_{\alpha+\beta}=-1$,
the simple-root squares give the value $2$, and the other three pairs
give the value $1$. Thus 
\[
\operatorname{Spec}(\mathcal{C})=\{-1,1,2\},\qquad\lambda_{\min}(\mathcal{C})=-2+\frac{2}{G_{*}}=-1.
\]
The lower estimate in Theorem~\ref{theorem:4.4} is therefore attained,
while the eigenvalue $2$ occurs precisely on the two simple-root
squares. \end{remark}

\subsection{The lower spectral bound in types $B$ and $C$}

We use matrix root vectors to compute the blocks in types $B$ and
$C$. Throughout this subsection, $i,j,k,l$ are coordinate indices,
while $a,b$ label the unitary blocks $I_{a},I_{b}$, of sizes $n_{a},n_{b}$.
The coordinates $t_{i}$ describe the canonical metric and the ranges
of indices; curvature entries and estimates use the metric coefficients
$g_{\alpha}$. Each block $\mathcal{C}_{\tau}$ is written in its
specified unit symmetric-pair basis from \eqref{eq:2.5}.

\begin{lemma}[Matrix conventions]\label{lemma:4.7} Let $F_{ij}$
denote the elementary matrices, with indices $1,\ldots,r,\bar{1},\ldots,\bar{r}$,
and with an additional index $0$ in type $B_{r}$. Put $H_{i}=F_{ii}-F_{\bar{i}\bar{i}}$
and use $B(U,V)=\frac{1}{2}\operatorname{tr}(UV)$. In both types
set 
\[
E_{e_{i}-e_{j}}=F_{ij}-F_{\bar{j}\bar{i}}\quad(i<j).
\]
Choose the remaining positive root vectors as 
\[
\begin{array}{c|cc}
 & E_{e_{i}+e_{j}}\ (i<j) & \text{\itshape one-coordinate roots}\\
\hline B_{r} & F_{i\bar{j}}-F_{j\bar{i}} & E_{e_{i}}=F_{i0}-F_{0\bar{i}}\\
C_{r} & F_{i\bar{j}}+F_{j\bar{i}} & E_{2e_{i}}=F_{i\bar{i}}.
\end{array}
\]
For every positive root define $E_{-\alpha}=E^{\mathsf{T}}_{\alpha}$.
Then conjugation is $\sigma(U)=-U^{*}$, the coordinate vectors satisfy
$(e_{i},e_{j})=\delta_{ij}$, and 
\[
z_{\alpha}=B(E_{\alpha},E_{-\alpha})=\begin{cases}
\frac{1}{2}, & \text{\itshape type }C_{r},\ \alpha=\pm2e_{i},\\
1, & \text{\itshape otherwise.}
\end{cases}
\]
Moreover, $[E_{\alpha},E_{-\alpha}]=z_{\alpha}H_{\alpha}$, where
$H_{\alpha}=\sum_{i}\alpha_{i}H_{i}$ for $\alpha=\sum_{i}\alpha_{i}e_{i}$.
\end{lemma} 
\begin{proof}
The displayed matrices belong respectively to the orthogonal algebra
for the symmetric form pairing $i$ with $\bar{i}$ and assigning
norm 1 to the $0$ coordinate, and to the symplectic algebra for the
standard alternating form pairing $i$ with $\bar{i}$. Direct multiplication
gives $B(H_{i},H_{j})=\delta_{ij}$ and the stated norms and opposite-root
brackets. These algebras are stable under conjugate transpose, giving
the asserted conjugation. 
\end{proof}

The resulting curvature blocks are real symmetric. Changing the phases
of the root vectors conjugates each block by a diagonal unitary matrix
and preserves its spectrum. For $i<j<k$, the pairs $(e_{j}-e_{k},e_{i}+e_{k})$
and $(e_{j}+e_{k},e_{i}-e_{k})$ have the bracket pairings $N_{1}=0$,
$N_{2}=1$ in type $B$, and $N_{1}=2$, $N_{2}=1$ in type $C$,
in the notation of Proposition~\ref{proposition:2.1}. The additional
term in type $C$ comes from the root $2e_{k}$: the two brackets
are $-2E_{-2e_{k}}$ and $-2E_{2e_{k}}$, and their pairing is $4z_{2e_{k}}=2$.

\begin{theorem}\label{theorem:4.8} For the canonical metric with
$\Ric=2g$ on every generalized flag manifold of type $B_{r}$ or
$C_{r}$, the Calabi operator satisfies $\mathcal{C}>-2\operatorname{Id}$.
\end{theorem}

\subsubsection{Canonical coordinates and the block list}

A parabolic subgroup is specified by unitary group blocks of sizes
$n_{1},\ldots,n_{\ell}$ and a residual orthogonal or symplectic block
indexed by a zero tail of size $n_{0}$, with $\sum_{a}n_{a}+n_{0}=r$.
The nonnegative coordinates $t_{i}$ of $\delta/2$ are constant on
these blocks, strictly decreasing between positive blocks, and zero
on the tail. For each coordinate index $i$ in the positive block
$I_{a}$, 
\begin{equation}
t_{i}=n_{0}+\sum_{b>a}n_{b}+\frac{n_{a}+\varepsilon}{2},\quad\varepsilon=0\ (B),\quad\varepsilon=1\ (C).\label{eq:4.2}
\end{equation}
This follows by subtracting the Levi half-sum from the usual vectors
$\rho_{B}=(r-\frac{1}{2},r-\frac{3}{2},\ldots,\frac{1}{2})$ and $\rho_{C}=(r,r-1,\ldots,1)$.
The metric coefficients are 
\[
g_{e_{i}-e_{j}}=t_{i}-t_{j},\quad g_{e_{i}+e_{j}}=t_{i}+t_{j},\quad g_{e_{i}}=t_{i}\ (B),\quad g_{2e_{i}}=2t_{i}\ (C).
\]
A positive root belongs to $\Delta^{+}_{\Phi_{T}}$ exactly when its
metric coefficient is positive. In particular, $t_{i}$ is greater
than the number of indices in all blocks after $I_{a}$, including
the zero tail.

As in Section~\ref{subsec:4.2}, write $\tau=\alpha+\gamma$ for
a root sum and $s=g_{\alpha}+g_{\gamma}$ for the sum of its metric
coefficients. The absolute values of the nonzero coefficients of $\tau$
have one of the following patterns: 
\[
\begin{array}{c|l}
B & (1),(2),(1,1),(2,1),(2,2),(1,1,1),(2,1,1),(1,1,1,1),\\
C & (2),(4),(1,1),(3,1),(2,2),(2,1,1),(1,1,1,1).
\end{array}
\]
We consider only sign choices for which $\tau$ is a sum of two positive
tangent roots. The number of decompositions of $e_{i}-e_{j}$, $e_{i}+e_{j}$,
and $2e_{i}$, and of $e_{i}$ in type $B$, can grow with the rank.
We treat these blocks first. All other blocks have dimension at most
three. By \eqref{eq:2.8}, a repeated pair, or a pair with positive
inner product, has positive diagonal curvature. An orthogonal pair
whose sum is not a root has zero diagonal curvature. These observations
already settle the one-dimensional blocks of patterns $(2,1)$, $(3,1)$,
$(2,1,1)$, and $(2,2)$ whenever they occur.

We shall use the following elementary criterion. Let $P,Q$ index
the listed unit root-pair basis of a block. If positive numbers $u_{P}$
satisfy 
\begin{equation}
2+(\mathcal{C}_{\tau})_{P,P}-\sum_{Q\ne P}|(\mathcal{C}_{\tau})_{P,Q}|\frac{u_{Q}}{u_{P}}>0,\label{eq:4.3}
\end{equation}
then $2\operatorname{Id}+\mathcal{C}_{\tau}$ is positive definite.
Indeed, 
\[
2|c_{P}||c_{Q}|\leq\frac{u_{Q}}{u_{P}}|c_{P}|^{2}+\frac{u_{P}}{u_{Q}}|c_{Q}|^{2}
\]
follows by expanding a square. Applying it to each off-diagonal term
gives 
\[
c^{*}(2\operatorname{Id}+\mathcal{C}_{\tau})c\geq\sum_{P}\left(2+(\mathcal{C}_{\tau})_{P,P}-\sum_{Q\ne P}|(\mathcal{C}_{\tau})_{P,Q}|\frac{u_{Q}}{u_{P}}\right)|c_{P}|^{2}>0\quad(c\ne0).
\]
For the tables below, write 
\[
W_{P,Q}=|(\mathcal{C}_{\tau})_{P,Q}|\frac{u_{Q}}{u_{P}}.
\]
These quantities depend on the direction from row $P$ to column $Q$;
in general $W_{P,Q}\ne W_{Q,P}$. All block entries are computed from
\eqref{eq:2.6} in the unit pair basis of \eqref{eq:2.5}; the diagonal
entries are given by \eqref{eq:2.8}.

\subsubsection{Difference sums and the short-root sums in type $B$}

Fix $i<j$ with $\tau=e_{i}-e_{j}\in\Delta^{+}_{\Phi_{T}}$. Its decompositions
are precisely $(e_{i}-e_{k},e_{k}-e_{j})$ with both roots tangent.
Index the corresponding unit symmetric pairs by $k$, in the coordinate
order, and let $m_{\tau}$ be their number. If $m_{\tau}=0$, this
block is absent; assume $m_{\tau}>0$. With $s=g_{e_{i}-e_{j}}$,
the entries are 
\[
\begin{aligned}(\mathcal{C}_{e_{i}-e_{j}})_{kk} & =-\frac{2}{s},\\
(\mathcal{C}_{e_{i}-e_{j}})_{kl} & =-\frac{2}{s}\sqrt{\frac{g_{e_{i}-e_{k}}g_{e_{l}-e_{j}}}{g_{e_{i}-e_{l}}g_{e_{k}-e_{j}}}}\qquad(k<l).
\end{aligned}
\]
Thus $(\mathcal{C}_{e_{i}-e_{j}})_{kl}$ denotes the entry with root-pair
indices $(e_{i}-e_{k},e_{k}-e_{j})$ and $(e_{i}-e_{l},e_{l}-e_{j})$
in \eqref{eq:2.5}. Since $k<l$ implies $t_{k}\geq t_{l}$, both
ratios under the square root are at most $1$. Every absolute row
sum is therefore at most $2m_{\tau}/s$, so 
\[
\|\mathcal{C}_{e_{i}-e_{j}}\|\leq\frac{2m_{\tau}}{s},\qquad\lambda_{\min}(\mathcal{C}_{e_{i}-e_{j}})\geq-\frac{2m_{\tau}}{s}.
\]
If $i\in I_{a}$ and $j\in I_{b}$ belong to positive blocks, then
$s=m_{\tau}+(n_{a}+n_{b})/2>m_{\tau}$. If $j$ belongs to the zero
tail, then $s=m_{\tau}+n_{0}+(n_{a}+\varepsilon)/2>m_{\tau}$. Hence
the smallest eigenvalue is greater than $-2$ in both cases.

For the type-$B$ sum $\tau=e_{i}\in\Delta^{+}_{\Phi_{T}}$, put $s=g_{e_{i}}$.
The rows are $(e_{i}-e_{k},e_{k})$ with $0<g_{e_{k}}<g_{e_{i}}$.
If there are no such indices, the block is absent. Otherwise, its
entries are 
\[
(\mathcal{C}_{e_{i}})_{kk}=-\frac{2}{s},\qquad(\mathcal{C}_{e_{i}})_{kl}=-\frac{2}{s}\sqrt{\frac{g_{e_{i}-e_{k}}g_{e_{l}}}{g_{e_{i}-e_{l}}g_{e_{k}}}}\quad(k<l).
\]
Again, each absolute entry is at most $2/s$. The number $m_{\tau}$
of rows is smaller than $s$ by \eqref{eq:4.2}, so $\lambda_{\min}(\mathcal{C}_{e_{i}})\geq-2m_{\tau}/s>-2$.

\subsubsection{The blocks with sum $2e_{i}$}

Fix $i\in I_{a}$, where $I_{a}$ is a positive coordinate block.
Let $k,l$ range over the indices in strictly later blocks, including
the zero tail, and put $N=n_{0}+\sum_{b>a}n_{b}$. The corresponding
rows of $\mathcal{C}_{2e_{i}}$ are the unit pairs $(e_{i}-e_{k},e_{i}+e_{k})$.
Equation~\eqref{eq:4.2} gives 
\[
g_{e_{i}}=N+\frac{n_{a}}{2}>N\quad(B),\qquad g_{2e_{i}}=2N+n_{a}+1>2N\quad(C).
\]
The common coefficient sum is $s=2g_{e_{i}}$ in type $B$ and $s=g_{2e_{i}}$
in type $C$.

In type $B$ there is also the repeated pair $(e_{i},e_{i})$, which
we use as its row label. After choosing the signs of the root vectors,
the entries are 
\[
\begin{aligned}(\mathcal{C}_{2e_{i}})_{(e_{i},e_{i}),(e_{i},e_{i})} & =\frac{1}{g_{e_{i}}}, & (\mathcal{C}_{2e_{i}})_{kk} & =0,\\
(\mathcal{C}_{2e_{i}})_{kl} & =-\frac{2\min(g_{e_{i}-e_{k}},g_{e_{i}-e_{l}})}{\sqrt{g_{e_{i}-e_{k}}g_{e_{i}+e_{k}}g_{e_{i}-e_{l}}g_{e_{i}+e_{l}}}} &  & (k\ne l),\\
(\mathcal{C}_{2e_{i}})_{k,(e_{i},e_{i})} & =-\frac{\sqrt{2}\,g_{e_{i}-e_{k}}}{g_{e_{i}}\sqrt{g_{e_{i}-e_{k}}g_{e_{i}+e_{k}}}}.
\end{aligned}
\]
Take $u_{(e_{i},e_{i})}=g_{e_{i}}$ and $u_{k}=\sqrt{g_{e_{i}-e_{k}}g_{e_{i}+e_{k}}/2}$,
where $u_{k}$ is attached to the pair $(e_{i}-e_{k},e_{i}+e_{k})$.
In row $k$, the weighted off-diagonal sum is 
\[
\frac{2}{g_{e_{i}-e_{k}}g_{e_{i}+e_{k}}}\left(\sum_{l\ne k}\min(g_{e_{i}-e_{k}},g_{e_{i}-e_{l}})+g_{e_{i}-e_{k}}\right)\leq\frac{2N}{g_{e_{i}+e_{k}}}\leq\frac{2N}{g_{e_{i}}}<2.
\]
In the repeated-pair row it is 
\[
\sum_{k}\frac{g_{e_{i}-e_{k}}}{g^{2}_{e_{i}}}\leq\frac{N}{g_{e_{i}}}<1.
\]
The diagonal entries are nonnegative, so \eqref{eq:4.3} applies.
If $N=0$, only the positive repeated-pair row remains.

In type $C$ there is no repeated pair, since $e_{i}$ is not a root.
For $N\geq1$, the entries are 
\[
\begin{aligned}(\mathcal{C}_{2e_{i}})_{kk} & =\frac{4g_{e_{i}-e_{k}}}{g_{2e_{i}}g_{e_{i}+e_{k}}},\\
(\mathcal{C}_{2e_{i}})_{kl} & =-\frac{2g_{e_{i}-e_{k}}g_{2e_{l}}}{g_{2e_{i}}\sqrt{g_{e_{i}-e_{k}}g_{e_{i}+e_{k}}g_{e_{i}-e_{l}}g_{e_{i}+e_{l}}}}\qquad(k<l).
\end{aligned}
\]
Use $u_{k}=\sqrt{g_{e_{i}-e_{k}}g_{e_{i}+e_{k}}}$. If $l>k$, the
weighted entry in row $k$ is 
\[
\frac{2g_{2e_{l}}}{g_{2e_{i}}g_{e_{i}+e_{k}}}\leq\frac{2g_{2e_{k}}}{g_{2e_{i}}g_{e_{i}+e_{k}}}.
\]
If $l<k$, it is 
\[
\frac{2g_{e_{i}-e_{l}}g_{2e_{k}}}{g_{2e_{i}}g_{e_{i}-e_{k}}g_{e_{i}+e_{k}}}\leq\frac{2g_{2e_{k}}}{g_{2e_{i}}g_{e_{i}+e_{k}}},
\]
since $g_{e_{i}-e_{l}}\leq g_{e_{i}-e_{k}}$. Thus the row sum is
at most 
\[
\frac{2(N-1)g_{2e_{k}}}{g_{2e_{i}}g_{e_{i}+e_{k}}}\leq\frac{2(N-1)}{g_{2e_{i}}}<1,
\]
where $g_{e_{i}+e_{k}}-g_{2e_{k}}=g_{e_{i}-e_{k}}>0$. For a zero-tail
index, $g_{2e_{k}}=0$ and every off-diagonal entry in that row vanishes.
The diagonal is positive, so \eqref{eq:4.3} proves $\mathcal{C}_{2e_{i}}>-2\operatorname{Id}$.
If $N=0$, the block is absent.

\subsubsection{The type-$C$ sum $2e_{i}+2e_{j}$}

Fix $i<j$ with $g_{2e_{i}}\geq g_{2e_{j}}\geq0$ and $s=g_{2e_{i}}+g_{2e_{j}}>0$.
When $g_{2e_{j}}>0$, the ordered unit pair basis is 
\[
(e_{i}+e_{j},e_{i}+e_{j}),\qquad(2e_{i},2e_{j}).
\]
After choosing signs for this basis, the block is 
\[
\mathcal{C}_{2e_{i}+2e_{j}}=\begin{pmatrix}\dfrac{4}{s} & -\dfrac{4}{s}\sqrt{\dfrac{2g_{2e_{j}}}{g_{2e_{i}}}}\\[4pt]
-\dfrac{4}{s}\sqrt{\dfrac{2g_{2e_{j}}}{g_{2e_{i}}}} & 0
\end{pmatrix}.
\]
Its trace is $4/s>0$. The upper bound from Theorem~\ref{theorem:4.1}
therefore gives 
\[
\lambda_{\min}(\mathcal{C}_{2e_{i}+2e_{j}})=\operatorname{tr}\mathcal{C}_{2e_{i}+2e_{j}}-\lambda_{\max}(\mathcal{C}_{2e_{i}+2e_{j}})\geq\frac{4}{s}-2>-2.
\]
If $g_{2e_{j}}=0$, the second pair is absent and only the positive
one-dimensional block remains. The sum $2e_{i}-2e_{j}$ has only a
repeated-root pair.

\subsubsection{The type-$B$ sums $e_{i}+e_{j}$}

Fix $i<j$ with $g_{e_{i}+e_{j}}>0$, write $\tau=e_{i}+e_{j}$ and
$s=g_{e_{i}+e_{j}}$, and use $u_{(\alpha,\gamma)}=\sqrt{g_{\alpha}g_{\gamma}/2}$
in \eqref{eq:4.3}. The row labels used here and in type $C$ are
\[
\begin{aligned}A_{ij;k} & =(e_{j}-e_{k},e_{i}+e_{k}),\\
B_{ij;k} & =(e_{j}+e_{k},e_{i}-e_{k}),\\
L_{ij;k} & =(e_{i}-e_{k},e_{j}+e_{k}).
\end{aligned}
\]
For $A_{ij;k},B_{ij;k}$, the index $k$ lies in a block strictly
after the block containing $j$. For $L_{ij;k}$, require $k\ne j$,
$t_{j}\leq t_{k}<t_{i}$, and $t_{j}+t_{k}>0$. These conditions ensure
that both roots of every retained pair are tangent. Different indices
give different rows even when their metric coordinates agree. The
additional type-$B$ row is $T^{B}_{ij}=(e_{i},e_{j})$, present when
both roots are tangent. All these labels denote the unit symmetric
pairs of \eqref{eq:2.5}.

First suppose $i\in I_{a}$, $j\in I_{b}$ with $a<b$ and both blocks
positive. Let $p$ count the indices in blocks between $I_{a}$ and
$I_{b}$, and let $q$ count all indices after $I_{b}$, including
the zero tail. There are $q$ rows of each of the families $A,B$
and $h=p+n_{b}-1$ rows of the family $L$. Equation~\eqref{eq:4.2}
gives 
\[
g_{e_{i}}=q+n_{b}+p+\frac{n_{a}}{2},\qquad g_{e_{j}}=q+\frac{n_{b}}{2},\qquad s=2q+h+1+\frac{n_{a}+n_{b}}{2}.
\]
The $A,B,L$ diagonal entries of $\mathcal{C}_{e_{i}+e_{j}}$ are
$-2/s$, and its $T^{B}_{ij}$ diagonal entry is $2g_{e_{j}}/(sg_{e_{i}})>0$.
The following table gives the weighted absolute entries $W_{P,Q}$.
When two columns are listed, the formula is their sum. \begingroup
\setlength{\LTpre}{6pt} \setlength{\LTpost}{6pt} 
\global\long\def\arraystretch{1.6}%
\begin{longtable}[c]{c|c|c}
row $P$  & column(s) $Q$  & weighted entry or sum\tabularnewline
\endhead
\hline 
$A_{ij;k}$  & $A_{ij;l},B_{ij;l}$, $l\ne k$  & ${\displaystyle \text{sum}\leq\frac{4g_{e_{i}}}{sg_{e_{i}+e_{k}}}\leq\frac{4}{s}}$\tabularnewline
$B_{ij;k}$  & $A_{ij;l},B_{ij;l}$, $l\ne k$  & ${\displaystyle \text{sum}\leq\frac{4g_{e_{j}}}{sg_{e_{j}+e_{k}}}\leq\frac{4}{s}}$\tabularnewline
$A_{ij;k}$  & $B_{ij;k}$  & ${\displaystyle \frac{2g_{e_{j}+e_{k}}}{sg_{e_{i}+e_{k}}}\leq\frac{2}{s}}$\tabularnewline
$B_{ij;k}$  & $A_{ij;k}$  & ${\displaystyle \frac{2g_{e_{j}-e_{k}}}{sg_{e_{i}-e_{k}}}\leq\frac{2}{s}}$\tabularnewline
$A_{ij;k}$  & $L_{ij;l}$  & ${\displaystyle \frac{2g_{e_{i}-e_{l}}}{sg_{e_{i}+e_{k}}}\leq\frac{2}{s}}$\tabularnewline
$B_{ij;k}$  & $L_{ij;l}$  & ${\displaystyle \frac{2g_{e_{i}-e_{l}}}{sg_{e_{i}-e_{k}}}\leq\frac{2}{s}}$\tabularnewline
$A_{ij;k}$  & $T^{B}_{ij}$  & ${\displaystyle \frac{2g_{e_{i}}}{sg_{e_{i}+e_{k}}}\leq\frac{2}{s}}$\tabularnewline
$B_{ij;k}$  & $T^{B}_{ij}$  & ${\displaystyle \frac{2g_{e_{j}}}{sg_{e_{j}+e_{k}}}\leq\frac{2}{s}}$\tabularnewline
$L_{ij;k}$  & $A_{ij;l},B_{ij;l}$  & ${\displaystyle \text{sum}=\frac{4g_{e_{j}}}{sg_{e_{j}+e_{k}}}\leq\frac{2}{s}}$\tabularnewline
$L_{ij;k}$  & $L_{ij;l}$, $l\ne k$  & ${\displaystyle \leq\frac{2}{s}}$\tabularnewline
$L_{ij;k}$  & $T^{B}_{ij}$  & ${\displaystyle \frac{2g_{e_{j}}}{sg_{e_{j}+e_{k}}}\leq\frac{1}{s}}$\tabularnewline
$T^{B}_{ij}$  & $A_{ij;k},B_{ij;k}$  & ${\displaystyle \text{sum}=\frac{2}{s}\left(\frac{g_{e_{j}-e_{k}}}{g_{e_{j}}}+\frac{g_{e_{i}-e_{k}}}{g_{e_{i}}}\right)\leq\frac{4}{s}}$\tabularnewline
$T^{B}_{ij}$  & $L_{ij;k}$  & ${\displaystyle \frac{2g_{e_{i}-e_{k}}}{sg_{e_{i}}}\leq\frac{2}{s}}$ \tabularnewline
\end{longtable}\endgroup For example, the unnormalized curvature component between
the first two families with the same index is 
\[
\left|R_{e_{j}-e_{k},\overline{e_{j}+e_{k}},e_{i}+e_{k},\overline{e_{i}-e_{k}}}\right|=\frac{g_{e_{j}-e_{k}}g_{e_{j}+e_{k}}}{s},
\]
since $2e_{k}$ is not a root in type $B$. The components involving
$T^{B}_{ij}$ satisfy 
\[
\begin{aligned}\left|R_{e_{i},\overline{e_{j}-e_{k}},e_{j},\overline{e_{i}+e_{k}}}\right| & =\frac{g_{e_{i}}g_{e_{j}-e_{k}}}{s},\\
\left|R_{e_{i},\overline{e_{j}+e_{k}},e_{j},\overline{e_{i}-e_{k}}}\right| & =\frac{g_{e_{j}}g_{e_{i}-e_{k}}}{s},\\
\left|R_{e_{i},\overline{e_{i}-e_{k}},e_{j},\overline{e_{j}+e_{k}}}\right| & =\frac{g_{e_{j}}g_{e_{i}-e_{k}}}{s}.
\end{aligned}
\]
The sums over $A_{ij;l},B_{ij;l}$ are estimated as in the type-$C$
calculation below.

For an $A$ or $B$ row, the absolute value of the diagonal entry
plus the weighted off-diagonal sum is at most 
\[
\frac{2+4(q-1)+2+2h+2}{s}=\frac{2(2q+h+1)}{s}=2-\frac{n_{a}+n_{b}}{s}<2.
\]
For an $L$ row it is at most $\{2+2(h-1)+2q+1\}/s<2$, since the
numerator is less than $2s$ by $2q+1+n_{a}+n_{b}$. The off-diagonal
sum in the $T^{B}_{ij}$ row is at most $(4q+2h)/s<2$, and its diagonal
is positive. Thus \eqref{eq:4.3} holds.

If $i,j$ belong to the same positive block $I_{a}$, then $n_{a}\geq2$
and $s=2q+n_{a}$. The rows are $A_{ij;k},B_{ij;k}$ for the $q$
later indices, together with $T^{B}_{ij}$. The same bounds give $(4q+2)/s<2$
for each $A,B$ row, including its absolute diagonal entry, and $4q/s<2$
for the off-diagonal sum in the $T^{B}_{ij}$ row. If $g_{e_{j}}=0$,
only the rows $L_{ij;k}$ with $0<t_{k}<t_{i}$ remain. Let $N$ be
their number. If $N=0$, the block is absent. Otherwise, as in the
difference-sum calculation, every absolute entry is at most $2/s$,
and 
\[
s=g_{e_{i}+e_{j}}=N+n_{0}+\frac{n_{a}}{2}>N.
\]
Hence $\lambda_{\min}(\mathcal{C}_{e_{i}+e_{j}})\geq-2N/s>-2$ also
in this case.

\subsubsection{The type-$C$ blocks with sum $e_{i}+e_{j}$}

Fix $i<j$ with $g_{e_{i}+e_{j}}>0$ and use the same pair labels
$A_{ij;k},B_{ij;k},L_{ij;k}$ and index ranges as above, with $\tau=e_{i}+e_{j}$
and $s=g_{e_{i}+e_{j}}$. The additional pair is $T^{C}_{ij}=(e_{i}-e_{j},2e_{j})$,
which occurs precisely when $t_{i}>t_{j}>0$. We again use $u_{(\alpha,\gamma)}=\sqrt{g_{\alpha}g_{\gamma}/2}$
and the directional entries $W_{P,Q}$ defined above. Every pair in
this block has distinct roots.

First suppose $i\in I_{a}$, $j\in I_{b}$ with $a<b$ and both blocks
positive. Let $p$ count the indices in blocks between them and $q$
count the indices after $I_{b}$, including the zero tail. Then 
\[
\begin{aligned}\frac{g_{2e_{i}}}{2} & =q+n_{b}+p+\frac{n_{a}+1}{2}, & \frac{g_{2e_{j}}}{2} & =q+\frac{n_{b}+1}{2},\\
g_{e_{i}-e_{j}} & =p+\frac{n_{a}+n_{b}}{2}, & s & =2q+g_{e_{i}-e_{j}}+n_{b}+1.
\end{aligned}
\]
There are $q$ rows of each of the families $A,B$, $h=p+n_{b}-1$
rows of the family $L$, and the row $T^{C}_{ij}$. The $A,B,L$ diagonal
entries are $-2/s$ and the $T^{C}_{ij}$ diagonal entry is $-4/s$.
Equation~\eqref{eq:2.6} gives the following weighted entries. As
before, two listed columns mean the sum of their entries. \Needspace{30\baselineskip}
\begingroup \setlength{\LTpre}{6pt} \setlength{\LTpost}{6pt} 
\global\long\def\arraystretch{1.6}%
\begin{longtable}[c]{c|c|c}
row $P$  & column(s) $Q$  & weighted entry or sum\tabularnewline
\endhead
\hline 
$A_{ij;k}$  & $B_{ij;k}$  & ${\displaystyle \frac{2(g_{2e_{i}}+g_{e_{j}-e_{k}})}{sg_{e_{i}+e_{k}}}}$\tabularnewline
$B_{ij;k}$  & $A_{ij;k}$  & ${\displaystyle \frac{2g_{e_{j}-e_{k}}(g_{2e_{i}}+g_{e_{j}-e_{k}})}{sg_{e_{j}+e_{k}}g_{e_{i}-e_{k}}}}$\tabularnewline
$A_{ij;k}$  & $A_{ij;l},B_{ij;l}$, $l\ne k$  & ${\displaystyle \text{sum}\leq\frac{2g_{2e_{i}}}{sg_{e_{i}+e_{k}}}}$\tabularnewline
$B_{ij;k}$  & $A_{ij;l},B_{ij;l}$, $l\ne k$  & ${\displaystyle \text{sum}\leq\frac{2g_{2e_{j}}}{sg_{e_{j}+e_{k}}}}$\tabularnewline
$A_{ij;k}$  & $L_{ij;l}$  & ${\displaystyle \frac{2g_{e_{i}-e_{l}}}{sg_{e_{i}+e_{k}}}}$\tabularnewline
$B_{ij;k}$  & $L_{ij;l}$  & ${\displaystyle \frac{2g_{e_{i}-e_{l}}}{sg_{e_{i}-e_{k}}}}$\tabularnewline
$A_{ij;k}$  & $T^{C}_{ij}$  & ${\displaystyle \frac{2\sqrt{2}\,g_{e_{i}-e_{j}}}{sg_{e_{i}+e_{k}}}}$\tabularnewline
$B_{ij;k}$  & $T^{C}_{ij}$  & ${\displaystyle \frac{2\sqrt{2}\,g_{e_{i}-e_{j}}}{sg_{e_{i}-e_{k}}}}$\tabularnewline
$L_{ij;k}$  & $A_{ij;l},B_{ij;l}$  & ${\displaystyle \text{sum}=\frac{2g_{2e_{j}}}{sg_{e_{j}+e_{k}}}}$\tabularnewline
$L_{ij;k}$  & $T^{C}_{ij}$  & ${\displaystyle \frac{2\sqrt{2}\,g_{2e_{j}}}{sg_{e_{j}+e_{k}}}}$\tabularnewline
$L_{ij;k}$  & $L_{ij;l}$, $l\ne k$  & ${\displaystyle \leq\frac{2}{s}}$\tabularnewline
$T^{C}_{ij}$  & $A_{ij;k},B_{ij;k}$  & ${\displaystyle \text{sum}=\frac{2\sqrt{2}}{s}}$\tabularnewline
$T^{C}_{ij}$  & $L_{ij;k}$  & ${\displaystyle \frac{2\sqrt{2}\,g_{e_{i}-e_{k}}}{sg_{e_{i}-e_{j}}}}$ \tabularnewline
\end{longtable}\endgroup To verify the first entry, take $\alpha=e_{j}-e_{k}$,
$\gamma=e_{i}+e_{k}$, $\beta=e_{j}+e_{k}$, and $\epsilon=e_{i}-e_{k}$
in Proposition~\ref{proposition:2.1}. The matrix brackets give $N_{1}=2$
and $N_{2}=1$. The order $t_{i}>t_{j}>t_{k}\geq0$ therefore gives
\[
\begin{aligned}R_{e_{j}-e_{k},\overline{e_{j}+e_{k}},e_{i}+e_{k},\overline{e_{i}-e_{k}}} & =\frac{2g_{e_{j}-e_{k}}g_{e_{i}-e_{k}}+g_{e_{j}-e_{k}}g_{e_{j}+e_{k}}}{s}\\
 & =\frac{g_{e_{j}-e_{k}}(g_{2e_{i}}+g_{e_{j}-e_{k}})}{s}.
\end{aligned}
\]
The pair factor in \eqref{eq:2.5} is $2$, and all four root-vector
norms $z_{\alpha}$ here equal $1$. Normalizing and multiplying by
$u_{B_{ij;k}}/u_{A_{ij;k}}$ gives 
\[
W_{A_{ij;k},B_{ij;k}}=\frac{2R_{e_{j}-e_{k},\overline{e_{j}+e_{k}},e_{i}+e_{k},\overline{e_{i}-e_{k}}}}{g_{e_{j}-e_{k}}g_{e_{i}+e_{k}}}=\frac{2(g_{2e_{i}}+g_{e_{j}-e_{k}})}{sg_{e_{i}+e_{k}}}.
\]
The other entries follow in the same way, retaining $z_{2e_{j}}=1/2$
when $T^{C}_{ij}$ occurs. In particular, 
\[
\begin{aligned}\left|R_{e_{j}-e_{k},\overline{e_{i}-e_{l}},e_{i}+e_{k},\overline{e_{j}+e_{l}}}\right| & =\frac{g_{e_{i}-e_{l}}g_{e_{j}-e_{k}}}{s},\\
\left|R_{e_{i}-e_{j},\overline{e_{i}-e_{k}},2e_{j},\overline{e_{j}+e_{k}}}\right| & =\frac{g_{2e_{j}}g_{e_{i}-e_{k}}}{s}.
\end{aligned}
\]
For the sum over $A_{ij;l},B_{ij;l}$ with $l\ne k$ in row $A_{ij;k}$,
when $t_{l}>t_{k}$ we obtain 
\[
W_{A_{ij;k},A_{ij;l}}+W_{A_{ij;k},B_{ij;l}}=\frac{2}{s}\left(\frac{g_{e_{j}-e_{l}}}{g_{e_{j}-e_{k}}}+\frac{g_{e_{i}-e_{l}}}{g_{e_{i}+e_{k}}}\right)\leq\frac{2g_{2e_{i}}}{sg_{e_{i}+e_{k}}}.
\]
When $t_{l}\leq t_{k}$, the sum equals the right-hand side. The other
sums are estimated by the same comparison. Distinct indices with equal
metric coordinates remain separate terms.

We now verify \eqref{eq:4.3}. Since $n_{a},n_{b}\geq1$, we have
$(n_{a}+n_{b})/2\geq1$. In row $A_{ij;k}$, multiply its left-hand
side by $sg_{e_{i}+e_{k}}/2$ and use $g_{e_{i}-e_{l}}\leq g_{e_{i}-e_{j}}$
for $L_{ij;l}$. The result is at least 
\[
\frac{g_{2e_{i}}}{2}+\frac{hg_{2e_{j}}}{2}+\frac{sg_{2e_{k}}}{2}-(\sqrt{2}-1)g_{e_{i}-e_{j}}>0,
\]
because $g_{e_{i}-e_{j}}<g_{2e_{i}}/2$ and $\sqrt{2}-1<1$. For row
$B_{ij;k}$, multiplication by $sg_{e_{j}+e_{k}}g_{e_{i}-e_{k}}/2$
similarly gives the lower bound 
\[
\begin{aligned} & \frac{g_{2e_{j}}}{2}\left\{ \frac{g_{2e_{i}}}{2}-(\sqrt{2}-1)g_{e_{i}-e_{j}}\right\} +\frac{hg^{2}_{2e_{j}}}{4}\\
 & \qquad+\frac{g_{2e_{k}}}{2}\left\{ (2-\sqrt{2})g_{e_{i}-e_{j}}+g_{2e_{i}}(q+1)-\frac{sg_{2e_{k}}}{2}\right\} >0.
\end{aligned}
\]
Indeed, every index after $I_{b}$ satisfies $g_{2e_{k}}\leq2q$,
so 
\[
g_{2e_{i}}(q+1)-\frac{sg_{2e_{k}}}{2}\geq qg_{e_{i}-e_{j}}+g_{2e_{i}}>0.
\]
All three terms in the preceding lower bound are nonnegative, and
the first is strictly positive.

In row $L_{ij;k}$, the left-hand side of \eqref{eq:4.3} is at least
\[
\begin{aligned}\frac{2}{s}\left\{ s-h-\frac{g_{2e_{j}}(q+\sqrt{2})}{g_{e_{j}+e_{k}}}\right\}  & \geq\frac{2}{s}\{s-h-q-\sqrt{2}\}\\
 & =\frac{2}{s}\left\{ q+\frac{n_{a}+n_{b}+4}{2}-\sqrt{2}\right\} >0.
\end{aligned}
\]
For row $T^{C}_{ij}$, the finite sum of the coefficients contributed
by the $L$ rows is 
\[
U_{ij}=\sum_{\substack{k\ne j\\
t_{j}\leq t_{k}<t_{i}
}
}g_{e_{i}-e_{k}}=\frac{p(n_{a}+p)}{2}+(n_{b}-1)g_{e_{i}-e_{j}}.
\]
The equality follows from \eqref{eq:4.2}: the sum of the $p$ coordinates
$t_{k}$ in the intermediate blocks is $p(q+n_{b})+p(p+1)/2$. The
left-hand side of \eqref{eq:4.3} in row $T^{C}_{ij}$ is 
\[
2-\frac{4}{s}-\frac{2\sqrt{2}}{s}\left(q+\frac{U_{ij}}{g_{e_{i}-e_{j}}}\right).
\]
The identity 
\[
\begin{aligned}2g_{e_{i}-e_{j}}(g_{e_{i}-e_{j}}+n_{b}-1)-3U_{ij} & =\frac{p^{2}}{2}+p\left(\frac{n_{a}}{2}+n_{b}+1\right)\\
 & \qquad+\frac{(n_{a}+n_{b})(n_{a}+1)}{2}>0
\end{aligned}
\]
and $s-2=2q+g_{e_{i}-e_{j}}+n_{b}-1$ give $3(q+U_{ij}/g_{e_{i}-e_{j}})<2(s-2)$.
The row expression is therefore greater than 
\[
\left(2-\frac{4\sqrt{2}}{3}\right)\frac{s-2}{s}>0.
\]
Here $s\geq3$ because $n_{a},n_{b}\geq1$. This proves \eqref{eq:4.3}
when the two positive blocks are distinct.

If $i,j$ belong to the same positive block $I_{a}$, then $n_{a}\geq2$,
$g_{2e_{i}}=g_{2e_{j}}=2q+n_{a}+1$, and $s=g_{2e_{i}}$. Only $A_{ij;k},B_{ij;k}$
for the $q$ later indices occur. If $q=0$, the block is absent.
Otherwise, the same estimates, after multiplying the left-hand side
of \eqref{eq:4.3} by $sg_{e_{i}+e_{k}}/2$, give 
\[
\frac{g_{2e_{i}}}{2}(n_{a}-1)+\frac{sg_{2e_{k}}}{2}>0.
\]
Finally, if $g_{2e_{j}}=0$, only the rows $L_{ij;k}$ with $0<t_{k}<t_{i}$
occur. Let $N$ be their number. If $N=0$, the block is absent. For
$N>0$, the zero tail has size $n_{0}\geq1$ and 
\[
s=g_{e_{i}+e_{j}}=N+n_{0}+\frac{n_{a}+1}{2}>N.
\]
As in the difference-sum calculation, each absolute entry is at most
$2/s$, so $\lambda_{\min}(\mathcal{C}_{e_{i}+e_{j}})\geq-2N/s>-2$.
These cases exhaust the type-$C$ blocks with sum $e_{i}+e_{j}$.

\subsubsection{The remaining small root-sum blocks}

It remains to consider sums involving at least three coordinate indices
and the pattern $(3,1)$ in type $C$.

By Theorem~\ref{theorem:4.1}, $\mathcal{C}\leq2\operatorname{Id}$.
Proposition~\ref{proposition:4.2} identifies the eigenspace of eigenvalue
$2$ as the span of $e_{\alpha}\otimes e_{\alpha}$ for painted simple
roots $\alpha\in\Phi_{T}$ with anticanonical degree $k_{\alpha}=2$.
Consequently, on a root-sum block indexed by $\tau$ different from
all these $2\alpha$, 
\begin{equation}
\lambda_{\max}(\mathcal{C}_{\tau})<2.\label{eq:4.4}
\end{equation}
In particular, a block of size two with nonnegative trace satisfies
\[
\lambda_{\min}(\mathcal{C}_{\tau})=\operatorname{tr}\mathcal{C}_{\tau}-\lambda_{\max}(\mathcal{C}_{\tau})>-2
\]
whenever \eqref{eq:4.4} applies. With strictly positive trace, the
weak upper bound $\mathcal{C}_{\tau}\leq2\operatorname{Id}$ already
suffices. All zero-diagonal blocks considered below have root sums
distinct from these sums $2\alpha$.

\smallskip{}
 \noindent\emph{Patterns $(3,1)$ and $(2,1,1)$.} A type-$C$ sum
of pattern $(3,1)$ has at most the single decomposition $(2e_{i},e_{i}\pm e_{j})$,
when both summands are positive tangent roots. Their inner product
is positive, and the corresponding scalar block has positive diagonal
curvature. In type $B$, a sum of pattern $(2,1,1)$ has just one
decomposition, whose two roots have positive inner product and hence
positive diagonal curvature. In type $C$ there can additionally be
the decomposition $(2e_{i},e_{j}\pm e_{k})$, whose diagonal curvature
is zero. Thus the block has size at most two, and if it has size two
its trace is strictly positive. The preceding trace estimate gives
$2\operatorname{Id}+\mathcal{C}_{\tau}>0$. If only one of these decompositions
consists of positive tangent roots, the block is one-dimensional and
nonnegative.

\smallskip{}
 \noindent\emph{Type $B$, pattern $(1,1,1)$.} Three decompositions
occur only when $\tau=e_{i}+e_{j}+e_{k}$ and all three short roots
are tangent. Denote these three short roots, after permuting them,
by $\alpha_{1},\alpha_{2},\alpha_{3}$, with 
\[
0<g_{\alpha_{1}}\leq g_{\alpha_{2}}\leq g_{\alpha_{3}},\qquad\gamma_{m}=\tau-\alpha_{m}\quad(m=1,2,3).
\]
Thus row $m$ is the unit pair $(\alpha_{m},\gamma_{m})$ and $s=g_{\alpha_{1}}+g_{\alpha_{2}}+g_{\alpha_{3}}$.
The matrix $\mathcal{C}_{\tau}$ has zero diagonal and 
\[
|(\mathcal{C}_{\tau})_{mn}|=\frac{2\min(g_{\alpha_{m}},g_{\alpha_{n}})}{\sqrt{g_{\alpha_{m}}g_{\gamma_{m}}g_{\alpha_{n}}g_{\gamma_{n}}}}\qquad(m\ne n).
\]
Take $u_{m}=\sqrt{g_{\alpha_{m}}g_{\gamma_{m}}}$. The three weighted
absolute row sums satisfy 
\[
\begin{aligned}\sum_{n\ne1}|(\mathcal{C}_{\tau})_{1n}|\frac{u_{n}}{u_{1}} & =\frac{4}{g_{\gamma_{1}}}\leq\frac{6}{s},\\
\sum_{n\ne2}|(\mathcal{C}_{\tau})_{2n}|\frac{u_{n}}{u_{2}} & =\frac{2g_{\gamma_{3}}}{g_{\alpha_{2}}g_{\gamma_{2}}}\leq\frac{6}{s},\\
\sum_{n\ne3}|(\mathcal{C}_{\tau})_{3n}|\frac{u_{n}}{u_{3}} & =\frac{2}{g_{\alpha_{3}}}\leq\frac{6}{s}.
\end{aligned}
\]
The middle inequality follows from 
\[
\begin{aligned}3g_{\alpha_{2}}g_{\gamma_{2}}-sg_{\gamma_{3}} & =(g_{\alpha_{2}}-g_{\alpha_{1}})g_{\gamma_{3}}\\
 & \quad+(g_{\alpha_{3}}-g_{\alpha_{2}})(2g_{\alpha_{2}}-g_{\alpha_{1}})\geq0.
\end{aligned}
\]
The quadratic-form estimate used to prove \eqref{eq:4.3} gives $\mathcal{C}_{\tau}\geq-(6/s)\operatorname{Id}$.

We have $s\geq9/2$. Indeed, the nonzero type-$B$ coordinates are
averages over consecutive blocks of $\rho_{B}=(r-1/2,\ldots,1/2)$,
with the entries in a terminal $B_{n_{0}}$ Levi block replaced by
zero. Above this tail the averaged entries are those of $n_{0}+(r-n_{0}-1/2,\ldots,1/2)$.
For any $q$ selected nonzero coordinates, their sum is at least $qn_{0}+q^{2}/2\geq q^{2}/2$:
it suffices to select the $q$ smallest entries, and averaging a consecutive
block cannot lower the sum of its smallest partially selected entries.
Taking $q=3$ gives $s\geq9/2$. Hence 
\[
\mathcal{C}_{\tau}\geq-\frac{4}{3}\operatorname{Id},\qquad2\operatorname{Id}+\mathcal{C}_{\tau}\geq\frac{2}{3}\operatorname{Id}>0.
\]
If the sum has at most two decompositions, the block has zero diagonal
and \eqref{eq:4.4} gives $2\operatorname{Id}+\mathcal{C}_{\tau}>0$.

\smallskip{}
 \noindent\emph{Pattern $(1,1,1,1)$.} There are at most three decompositions,
corresponding to the three partitions of four distinct indices into
two unordered pairs. Every diagonal entry is zero. Equation~\eqref{eq:4.4}
settles blocks of dimension at most two, so suppose that all three
decompositions occur. The sum is not a root, and only differences
of roots contribute to the curvature formula.

In type $C$, the product of the three off-diagonal entries is positive.
To see this, first take the sum $e_{i}+e_{j}+e_{k}+e_{l}$ and its
ordered decompositions 
\[
(e_{i}+e_{j},e_{k}+e_{l}),\quad(e_{i}+e_{k},e_{j}+e_{l}),\quad(e_{i}+e_{l},e_{j}+e_{k}).
\]
Using the symplectic root vectors in Lemma~\ref{lemma:4.7}, one
has, for distinct indices, 
\[
[E_{e_{i}+e_{j}},E_{-(e_{i}+e_{k})}]=E_{e_{j}-e_{k}}.
\]
Applying this identity to the brackets between the three pairs gives
positive bracket products. The metric factors in \eqref{eq:2.6} are
also positive. Signed coordinate permutations carry any other choice
of four signs to this one, while a change of basis signs preserves
the product of the three off-diagonal entries. Therefore 
\[
\det\mathcal{C}_{\tau}=2(\mathcal{C}_{\tau})_{12}(\mathcal{C}_{\tau})_{13}(\mathcal{C}_{\tau})_{23}>0,\qquad\operatorname{tr}\mathcal{C}_{\tau}=0.
\]
The three real eigenvalues hence consist of two negative eigenvalues
$\lambda_{1}\leq\lambda_{2}<0$ and one positive eigenvalue $\lambda_{3}$.
It follows that 
\[
\lambda_{1}=-\lambda_{2}-\lambda_{3}>-\lambda_{3}\geq-2.
\]
If an off-diagonal entry is zero, both the determinant and the trace
are zero. The spectrum then has the form $(-\lambda,0,\lambda)$,
with $\lambda\geq0$; by \eqref{eq:4.4}, $\lambda<2$. Thus this
case also has $2\operatorname{Id}+\mathcal{C}_{\tau}>0$.

In type $B$, the long roots $\{\pm e_{i}\pm e_{j}\}$ form a closed
subsystem of type $D_{r}$. Its intersection with $\Delta^{+}_{\Phi_{T}}$
is the positive tangent root set of a type-$D_{r}$ parabolic. The
Levi subsystem has the same unitary blocks and the subsystem $D_{n_{0}}$
on the zero tail. Denote its canonical metric coefficients by $g^{D}$.
The omitted tangent roots are precisely $e_{i}$ with $g_{e_{i}}>0$,
so for every long root $\alpha$ the given metric coefficients satisfy
\[
g_{\alpha}-g^{D}_{\alpha}=\frac{1}{2}\left(\alpha,\sum_{i:\,g_{e_{i}}>0}e_{i}\right).
\]
This difference is nonnegative on every positive tangent long root.
For $e_{i}+e_{j}$ this is immediate; for $e_{i}-e_{j}$ with $i<j$,
it is zero when $g_{e_{i}},g_{e_{j}}>0$, and is $\frac{1}{2}$ when
$g_{e_{i}}>0=g_{e_{j}}$. The opposite pattern is impossible for a
positive tangent difference root.

The three decompositions use only long roots, and every root arising
as a difference of their summands lies in this closed subsystem. Write
the three pairs as $(\alpha_{j},\gamma_{j})$, $j=1,2,3$, and orient
each so that $g_{\alpha_{j}}\leq g_{\gamma_{j}}$. The calculation
leading to \eqref{eq:4.1} then gives 
\[
\mathcal{C}_{\tau}=2\operatorname{diag}\left(\frac{1}{g_{\gamma_{j}}}\right)-2K,\qquad K_{jk}=\frac{\min(g_{\alpha_{j}},g_{\alpha_{k}})}{\sqrt{g_{\alpha_{j}}g_{\gamma_{j}}g_{\alpha_{k}}g_{\gamma_{k}}}}.
\]
The matrix $K$ is positive semidefinite by the same integral formula
as in Section~\ref{subsec:4.2}. The type-$D$ root count in the
proof of Theorem~\ref{theorem:4.4} gives $\max(g^{D}_{\alpha_{j}},g^{D}_{\gamma_{j}})\geq3$.
Since $g_{\alpha}\geq g^{D}_{\alpha}$ for every positive tangent
long root, we have $g_{\gamma_{j}}\geq3$. Hence $\operatorname{tr}K=\sum_{j}1/g_{\gamma_{j}}\leq1$,
and 
\[
2\operatorname{Id}+\mathcal{C}_{\tau}\geq2\operatorname{diag}\left(\frac{1}{g_{\gamma_{j}}}\right)>0.
\]

\subsubsection{Completion of the classical proof}

The coefficient-pattern list exhausts sums of two positive roots in
types $B$ and $C$. Every corresponding block satisfies $\mathcal{C}_{\tau}>-2\operatorname{Id}$.
There are finitely many blocks for a fixed flag, so $\lambda_{\min}(\mathcal{C})>-2$.
This proves Theorem~\ref{theorem:4.8}.

\subsection{The lower spectral bound in type $F_{4}$}

\begin{theorem}\label{theorem:4.9} For every generalized flag of
type $F_{4}$ with its canonical metric normalized by $\Ric=2g$,
the Calabi operator satisfies $\mathcal{C}>-2\operatorname{Id}$.
\end{theorem}

In an orthonormal basis of $\mathbb{R}^{4}$, the roots of $F_{4}$
are 
\[
\Delta_{\mathrm{long}}=\{\pm e_{i}\pm e_{j}:i<j\},\qquad\Delta_{\mathrm{short}}=\{\pm e_{i}\}\cup\left\{ \frac{1}{2}(\pm e_{1}\pm e_{2}\pm e_{3}\pm e_{4})\right\} .
\]
All sixteen sign choices occur in the last set. Put $H=\frac{\delta}{2}$
in the Euclidean root space, so $g_{\alpha}=(H,\alpha)$; the roots
with zero metric coefficient are precisely the Levi roots. Choose
Chevalley root vectors normalized by $z_{\alpha}=B(E_{\alpha},E_{-\alpha})=\frac{2}{|\alpha|^{2}}$,
and write $w_{\alpha}=g_{\alpha}z_{\alpha}$ as in \eqref{eq:2.4}.
Write $\mathcal{C}_{\tau}$ for the block indexed by $\tau$; it consists
of the unordered pairs $\alpha+\gamma=\tau$ with $\alpha,\gamma\in\Delta^{+}_{\Phi_{T}}$.

A Weyl-group transformation applied to both the roots and $H$ preserves
inner products and metric coefficients. We may therefore choose a
representative of each root sum up to the Weyl group. In the transformed
positive system, a root is positive and tangent precisely when $(H,\alpha)>0$.
We also compute brackets in closed root subsystems: a subsystem is
closed if the sum of two of its roots, whenever a root of $F_{4}$,
belongs to the subsystem. Thus the curvature formula has the same
bracket terms in the subsystem and in $F_{4}$.

The possible values of $|\tau|^{2}$ are $1,2,3,4,5,6,8$. Indeed,
a sum of two long roots has squared norm $2,4,6$, or $8$; a sum
of two short roots has squared norm $1,2,3$, or $4$; and a sum of
a long root and a short root has squared norm $1,3$, or $5$. A sum
of squared norm eight is a repeated long root. A sum of squared norm
six is a unique long-long pair, by the simply-laced argument. A sum
of squared norm five is a unique mixed pair: it can be sent to $2e_{1}+e_{2}$,
for which the only short root $\beta$ satisfying $(2e_{1}+e_{2},\beta)=2$
is $e_{1}$. The corresponding scalar blocks are positive by \eqref{eq:2.8}.
It remains to consider $|\tau|^{2}=1,2,3$, and $4$.

\subsubsection{Blocks whose sum is a long root}

Let $\tau$ be a long tangent root and put $s=g_{\tau}$. Write $m_{\tau}$
for the number of its unordered tangent decompositions. Every such
decomposition consists of two long roots or two short roots, as follows
from the coordinate description above.

As in Lemma~\ref{lemma:4.5}, the root count gives 
\[
s=1+m_{\tau}+\frac{h_{\tau}}{2},
\]
where $h_{\tau}$ counts tangent roots $\beta$ with $(\tau,\beta)=1$
and $\tau-\beta$ a Levi root. Indeed, every root other than $\pm\tau$
has pairing $0,1$, or $-1$ with $\tau$. Reflection in $\tau$ pairs
each tangent root with pairing $-1$ with $\beta+\tau$, cancelling
their contributions to $2s$. The remaining roots of pairing $1$
either form a tangent decomposition of $\tau$ or have a Levi difference.
In particular, $m_{\tau}<s$.

For any two tangent pairs $\alpha+\gamma=\beta+\epsilon=\tau$, we
claim 
\begin{equation}
\left|(\mathcal{C}_{\tau})_{\alpha\gamma,\beta\epsilon}\right|\le\frac{2}{s}.\label{eq:4.5}
\end{equation}
For a pair of long roots, the diagonal entry is $-2/s$. For a pair
of short roots, orient $\alpha,\gamma$ so that $g_{\alpha}\le g_{\gamma}$;
its diagonal entry is $2g_{\alpha}/(sg_{\gamma})$, which lies in
$[0,2/s]$. For two long-long pairs, the simply-laced estimate in
Subsection~\ref{subsec:4.2} gives \eqref{eq:4.5}.

For a long-long pair $(\alpha,\gamma)$ and a short-short pair $(\beta,\epsilon)$,
interchange the short roots if necessary so that $(\alpha,\beta)=1$.
Their off-diagonal entry satisfies 
\[
\left|(\mathcal{C}_{\tau})_{\alpha\gamma,\beta\epsilon}\right|=\frac{2}{s}\sqrt{\min\left\{ \frac{g_{\beta}g_{\gamma}}{g_{\epsilon}g_{\alpha}},\frac{g_{\epsilon}g_{\alpha}}{g_{\beta}g_{\gamma}}\right\} }\le\frac{2}{s}.
\]
This follows from the curvature formula in the closed $B_{3}$ subsystem
containing the two pairs. The two expressions according to the ordering
of the metric coefficients give the minimum above and agree when the
corresponding difference root is Levi.

For distinct short-short pairs, write $\alpha_{j}+\gamma_{j}=\tau$
and orient each pair so that $g_{\alpha_{j}}\le g_{\gamma_{j}}$.
The roots in the span of two such decompositions form a closed $C_{3}$
subsystem: the vectors $\tau/2$, $(\gamma_{j}-\alpha_{j})/2$, and
$(\gamma_{k}-\alpha_{k})/2$ are orthogonal and have squared length
$1/2$. These vectors are not roots. The squared root lengths in this
subsystem are half those in Lemma~\ref{lemma:4.7}. For the same
metric coefficients, the type-$C$ formula for the block with sum
$2e_{i}$ therefore gives half the curvature entries in the unit pair
basis. Thus, for $j\ne k$, 
\begin{equation}
\left|(\mathcal{C}_{\tau})_{\alpha_{j}\gamma_{j},\alpha_{k}\gamma_{k}}\right|=\frac{\min(g_{\alpha_{j}},g_{\alpha_{k}})\bigl(s-2\max(g_{\alpha_{j}},g_{\alpha_{k}})\bigr)}{s\sqrt{g_{\alpha_{j}}g_{\gamma_{j}}g_{\alpha_{k}}g_{\gamma_{k}}}}\le\frac{2}{s}.\label{eq:4.6}
\end{equation}
For the last inequality, suppose $g_{\alpha_{j}}\le g_{\alpha_{k}}$.
Since $g_{\gamma_{j}},g_{\gamma_{k}}\ge s/2$, the square-root factor
in the denominator is at least $(s/2)g_{\alpha_{j}}$, whereas the
numerator is at most $sg_{\alpha_{j}}$. Formula \eqref{eq:4.6} also
holds when $g_{\alpha_{j}}=g_{\gamma_{j}}$ or $g_{\alpha_{k}}=g_{\gamma_{k}}$,
so it includes the cases in which these difference roots are Levi
roots.

Thus every absolute row sum is at most $2m_{\tau}/s<2$, proving $\mathcal{C}_{\tau}>-2\operatorname{Id}$
when $\tau$ is a long root.

\subsubsection{Blocks with $|\tau|^{2}=4$}

Here $\tau=2u$ for a short root $u$. Indeed, two orthogonal long
roots have sum either $\pm2e_{i}$ or $\pm e_{1}\pm e_{2}\pm e_{3}\pm e_{4}$,
each twice a short root. The only other possible decomposition is
the repeated pair $(u,u)$.

Let $N$ count the long-long tangent decompositions of $2u$. Write
their roots as 
\[
\alpha_{j}=u-v_{j},\qquad\gamma_{j}=u+v_{j},
\]
Here each $v_{j}$ is a short root orthogonal to $u$. Orient it so
that $g_{v_{j}}\ge0$. Thus $0<g_{\alpha_{j}}\le g_{u}\le g_{\gamma_{j}}$
and $s=2g_{u}$. The vectors $u,v_{1},\ldots,v_{N}$ are orthonormal,
so $N\le3$. These roots and $u$ belong to a closed $B_{N+1}$ subsystem.
The type-$B$ curvature formulas give 
\[
\begin{aligned}(\mathcal{C}_{2u})_{uu,uu} & =\frac{1}{g_{u}}, & (\mathcal{C}_{2u})_{\alpha_{j}\gamma_{j},\alpha_{j}\gamma_{j}} & =0,\\
(\mathcal{C}_{2u})_{\alpha_{j}\gamma_{j},\alpha_{k}\gamma_{k}} & =-\frac{2\min(g_{\alpha_{j}},g_{\alpha_{k}})}{\sqrt{g_{\alpha_{j}}g_{\gamma_{j}}g_{\alpha_{k}}g_{\gamma_{k}}}} &  & (j\ne k),\\
(\mathcal{C}_{2u})_{\alpha_{j}\gamma_{j},uu} & =-\frac{\sqrt{2}\,g_{\alpha_{j}}}{g_{u}\sqrt{g_{\alpha_{j}}g_{\gamma_{j}}}},
\end{aligned}
\]
up to a diagonal change of root phases.

Here $g_{u}\ge N+1/2$. To prove this bound, evaluate $2g_{u}=\sum_{\beta\in\Delta^{+}_{\Phi_{T}}}(u,\beta)$.
Every long root with pairing $-1$ pairs with $\beta+2u$, a long
tangent root with pairing $1$, and cancels. The remaining positive
long-root contributions either come in pairs $\beta+(2u-\beta)=2u$,
contributing two per long-long decomposition, or have a Levi complement
and contribute nonnegatively. Thus long roots contribute at least
$2N$. Among short roots, the terms with pairing $-1/2$ cancel with
$\beta+u$; all remaining contributions are nonnegative. The root
$u$ itself contributes one. Hence $2g_{u}\ge2N+1$.

Use the positive scaling factors 
\[
U_{uu}=g_{u},\qquad U_{\alpha_{j}\gamma_{j}}=\sqrt{\frac{g_{\alpha_{j}}g_{\gamma_{j}}}{2}}
\]
in \eqref{eq:4.3}. The scaled off-diagonal sum in row $(\alpha_{j},\gamma_{j})$
is at most 
\[
\frac{2}{g_{\alpha_{j}}g_{\gamma_{j}}}\left(\sum_{k\ne j}\min(g_{\alpha_{j}},g_{\alpha_{k}})+g_{\alpha_{j}}\right)\le\frac{2N}{g_{\gamma_{j}}}\le\frac{2N}{g_{u}}<2.
\]
In row $(u,u)$, it is $\sum_{j}g_{\alpha_{j}}/g^{2}_{u}\le N/g_{u}<1$,
whereas the diagonal is positive. Strict scaled diagonal dominance
proves $\mathcal{C}_{2u}>-2\operatorname{Id}$. The case $N=0$ is
a positive one-dimensional block.

\subsubsection{Blocks with $|\tau|^{2}=3$}

Every root sum of squared norm three is Weyl-conjugate to $\tau=e_{1}+e_{2}+e_{3}$.
Its complete list of unordered root decompositions is 
\[
L_{i}=(e_{i},\tau-e_{i})\quad(i=1,2,3),\qquad S=(h_{-},h_{+}),\qquad h_{\pm}=\frac{1}{2}(\tau\pm e_{4}).
\]
Indeed, a long-long sum has even squared norm; an integer short-short
sum cannot have squared norm three; and the half-root coordinates
force the displayed unique short-short pair. The integer short-long
pairs are exactly the three displayed ones. To verify the Weyl-group
assertion, the coordinate possibilities are signed permutations of
$(1,1,1,0)$ and $\frac{1}{2}(3,1,1,1)$. Reflection in $\frac{1}{2}(e_{1}-e_{2}-e_{3}-e_{4})$
sends the first to $\frac{1}{2}(3,1,1,-1)$; signed permutations finish
the reduction.

Choose the sign of $e_{4}$ so that $g_{e_{4}}\ge0$. This fixes $\tau$
and gives $g_{h_{-}}\le g_{h_{+}}$. Put 
\[
s=g_{e_{1}}+g_{e_{2}}+g_{e_{3}}=(H,\tau).
\]
A pair belongs to $\mathcal{C}_{\tau}$ precisely when both its metric
coefficients are positive. In particular, $L_{i}$ is tangent when
$0<g_{e_{i}}<s$, and $S$ is tangent when $g_{e_{4}}<s$. The diagonal
and off-diagonal entries are 
\begin{align}
(\mathcal{C}_{\tau})_{L_{i},L_{i}} & =0,\qquad(\mathcal{C}_{\tau})_{S,S}=\frac{1}{g_{h_{+}}},\nonumber \\
|(\mathcal{C}_{\tau})_{L_{i},L_{j}}| & =\frac{2\min(g_{e_{i}},g_{\tau-e_{i}},g_{e_{j}},g_{\tau-e_{j}})}{\sqrt{g_{e_{i}}g_{\tau-e_{i}}g_{e_{j}}g_{\tau-e_{j}}}},\nonumber \\
|(\mathcal{C}_{\tau})_{L_{i},S}| & =\frac{\sqrt{2}\min(g_{e_{i}},g_{\tau-e_{i}},g_{h_{-}})}{\sqrt{g_{e_{i}}g_{\tau-e_{i}}g_{h_{-}}g_{h_{+}}}}.\label{eq:4.7}
\end{align}
Here $i\ne j$ in the off-diagonal formula, and only tangent pairs
are included. The order of the roots in $L_{i}$ is fixed by their
short and long lengths, not by the sizes of their coefficients.

To verify the constants, recall that $z_{\alpha}=2/(\alpha,\alpha)$,
so $z_{\mathrm{short}}=2$ and $z_{\mathrm{long}}=1$. The sum $\tau$
is not a root, so only difference-root terms contribute to an off-diagonal
entry. Between two $L$ pairs the difference root is long, and the
absolute bracket constants are $2$ and $1$. Between an $L$ pair
and $S$ the difference root is short, and the absolute bracket constants
are $1$ and $1$. Multiplication by $z$ and by the factor $2$ from
the unit pair basis, followed by division by the four root-vector
norms, gives respectively 
\[
\frac{2}{\sqrt{g_{e_{i}}g_{\tau-e_{i}}g_{e_{j}}g_{\tau-e_{j}}}},\qquad\frac{\sqrt{2}}{\sqrt{g_{e_{i}}g_{\tau-e_{i}}g_{h_{-}}g_{h_{+}}}}.
\]
The remaining coefficient is the least of the four metric coefficients,
giving the minima in \eqref{eq:4.7}. Finally, $(h_{+},h_{-})=1/2$,
and their difference is the short root $e_{4}$, which gives $(\mathcal{C}_{\tau})_{S,S}=1/g_{h_{+}}$.
The formulas extend continuously when a difference root is Levi, provided
that the roots indexing the entry remain tangent.

We apply the scaled diagonal-dominance criterion \eqref{eq:4.3}.
Since it involves only absolute values, the bracket signs need not
be chosen simultaneously. Use 
\[
U_{L_{i}}=\sqrt{g_{e_{i}}g_{\tau-e_{i}}},\qquad U_{S}=\sqrt{2g_{h_{-}}g_{h_{+}}}
\]
for the tangent pairs. We first bound the canonical metric coefficients.

The identity 
\[
H=\frac{1}{2}\sum_{\alpha:\,(H,\alpha)>0}\alpha=\frac{1}{4}\sum_{\alpha\in\Delta}\operatorname{sgn}(H,\alpha)\alpha
\]
follows from $H=\frac{\delta}{2}$, where $\operatorname{sgn}(0)=0$.
Write $H=b+f$, where $b$ is the vector contributed by the integer
$B_{4}$ roots and $f$ is the vector contributed by the half roots.
In a regular chamber, the absolute values of the nonzero coordinates
of $b$, in increasing order of the absolute coordinates of $H$,
are 
\[
\frac{1}{2},\quad\frac{3}{2},\quad\frac{5}{2},\quad\frac{7}{2},
\]
and their signs are those of $H$. For example, a positive coordinate
$(H,e_{i})$ contributes $\frac{1}{2}$ from $e_{i}$, and one further
unit for every other coordinate of smaller absolute value. On a wall,
tied values are replaced by their averages, and a zero coordinate
has contribution zero.

The half-root contribution has the same coordinate signs as $H$ or
is zero. To list its possible values, apply a signed permutation and,
for this coordinate calculation only, write the components of the
transformed $H$ as $H_{1}>H_{2}>H_{3}>H_{4}>0$. In the same coordinates,
$f$ is one of 
\begin{equation}
\begingroup
\global\long\def\arraystretch{1.25}%
\begin{array}{c|c}
\text{inequality} & f\\
\hline H_{1}>H_{2}+H_{3}+H_{4} & (2,0,0,0)\\
\begin{gathered}H_{1}<H_{2}+H_{3}+H_{4},\\
H_{1}+H_{4}>H_{2}+H_{3}
\end{gathered}
 & (\frac{3}{2},\frac{1}{2},\frac{1}{2},\frac{1}{2})\\
H_{1}+H_{4}<H_{2}+H_{3} & (1,1,1,0).
\end{array}\endgroup\label{eq:4.8}
\end{equation}
This table follows by grouping opposite half roots: for $H_{i}>0$,
with each $\epsilon_{j}\in\{-1,1\}$, 
\[
\begin{aligned}f_{i} & =\frac{1}{4}\#\left\{ (\epsilon_{j})_{j\ne i}:\left|\sum_{j\ne i}\epsilon_{j}H_{j}\right|<H_{i}\right\} \\
 & \quad+\frac{1}{8}\#\left\{ (\epsilon_{j})_{j\ne i}:\left|\sum_{j\ne i}\epsilon_{j}H_{j}\right|=H_{i}\right\} .
\end{aligned}
\]
Here $f_{i}$ denotes the corresponding component of the vector $f$.
In particular, it is nonnegative when the component of $H$ is positive.
We now return to the original coordinates, with $\tau=e_{1}+e_{2}+e_{3}$.
Boundary values of both $b$ and $f$ are averages of the values in
neighboring chambers: each root whose pairing vanishes contributes
zero, the average of its two signs. This also follows by averaging
sufficiently small generic perturbations with a distribution symmetric
about $H$. The bounds below apply to the root sums selected by each
such chamber; the perturbing vector itself need not be canonical.
Averaging those root-sum bounds therefore gives the asserted bounds
at the canonical vector on a wall.

If all three $L$ pairs belong to the block, then $g_{e_{1}},g_{e_{2}},g_{e_{3}}>0$.
Their integer-root contributions have sum at least $9/2$. If $S$
also belongs to the block, the half-root contribution to this sum
equals 
\[
\frac{3}{2}+\frac{1}{4}\sum^{3}_{i=1}\bigl[\operatorname{sgn}(s-2g_{e_{i}}+g_{e_{4}})+\operatorname{sgn}(s-2g_{e_{i}}-g_{e_{4}})\bigr].
\]
At most one $g_{e_{i}}$ exceeds $s/2$; hence at most one summand
can be negative, and it is at least $-2$. Consequently 
\begin{equation}
\begin{aligned}s & \ge\frac{11}{2} &  & \text{if all three \ensuremath{L} pairs and \ensuremath{S} occur},\\
s & \ge\frac{9}{2} &  & \text{if all three \ensuremath{L} pairs occur}.
\end{aligned}
\label{eq:4.9}
\end{equation}
If exactly two $L$ pairs and $S$ occur, relabel the first three
coordinates so that $g_{e_{1}}\ge g_{e_{2}}>0$ and $g_{e_{3}}\le0$.
The two tangent $L$ pairs imply $g_{e_{1}},g_{e_{2}}>-g_{e_{3}}$,
while $S$ implies $g_{e_{4}}-g_{e_{3}}<g_{e_{1}}+g_{e_{2}}$. First
take a regular chamber with $g_{e_{3}}<0$. If $g_{e_{4}}<g_{e_{1}}$,
then the first coordinate has largest absolute value. Its integer-root
contribution is $7/2$, and the contributions at $g_{e_{2}}$ and
$-g_{e_{3}}$ differ by at least one. Thus the integer contribution
to $s$ is at least $9/2$. In each row of \eqref{eq:4.8}, the half-root
contribution at the largest absolute coordinate is at least one, and
the contribution at $g_{e_{2}}$ is at least that at $-g_{e_{3}}$.
Their signed sum is therefore at least one.

If $g_{e_{4}}>g_{e_{1}}$, the coordinate order is $g_{e_{4}}>g_{e_{1}}>g_{e_{2}}>-g_{e_{3}}$,
so the integer contribution to $s$ is $5/2+3/2-1/2=7/2$. The inequality
$g_{e_{4}}-g_{e_{3}}<g_{e_{1}}+g_{e_{2}}$ selects the last row of
\eqref{eq:4.8}; the half-root contribution is $1+1-0=2$. Both cases
give $s\ge11/2$. When coordinates have equal absolute values, average
over neighboring chambers as above. At $g_{e_{3}}=0$, nearby chambers
have either this sign pattern or three positive first coordinates;
the strict inequalities defining the two tangent $L$ pairs and $S$
persist, and both cases give the same lower bound. It follows that
\[
s\ge\frac{11}{2}\qquad\text{if two \ensuremath{L} pairs and \ensuremath{S} occur}.
\]

Suppose first that two $L$ pairs and $S$ occur. The scaled off-diagonal
sum in row $L_{i}$ is 
\[
\begin{aligned} & \frac{2\min(g_{e_{i}},g_{\tau-e_{i}},g_{e_{j}},g_{\tau-e_{j}})+2\min(g_{e_{i}},g_{\tau-e_{i}},g_{h_{-}})}{g_{e_{i}}g_{\tau-e_{i}}}\\
 & \hspace{2cm}\le\frac{4}{\max(g_{e_{i}},g_{\tau-e_{i}})}\le\frac{8}{s}<2.
\end{aligned}
\]
The last inequality uses $s\ge11/2$. The $S$ row has scaled off-diagonal
sum at most $2/g_{h_{+}}$ and diagonal $1/g_{h_{+}}$. The diagonal
entry of $2\operatorname{Id}+\mathcal{C}_{\tau}$ therefore exceeds
the scaled off-diagonal sum by at least $2-1/g_{h_{+}}\ge2-2/s>0$.
Thus \eqref{eq:4.3} proves $2\operatorname{Id}+\mathcal{C}_{\tau}>0$.

Suppose next that all three $L$ pairs and $S$ occur. For $i\ne j$,
\[
\min(g_{e_{i}},g_{\tau-e_{i}},g_{e_{j}},g_{\tau-e_{j}})=\min(g_{e_{i}},g_{e_{j}}),
\]
because the three positive coefficients $g_{e_{1}},g_{e_{2}},g_{e_{3}}$
have sum $s$. We claim that the scaled off-diagonal sum in row $L_{i}$
and the corresponding sum minus the diagonal in row $S$ satisfy 
\[
\begin{aligned} & \frac{2\left(\sum_{j\ne i}\min(g_{e_{i}},g_{e_{j}})+\min(g_{h_{-}},g_{e_{i}},g_{\tau-e_{i}})\right)}{g_{e_{i}}g_{\tau-e_{i}}}\le\frac{9}{s},\\
 & \frac{\sum_{i}\min(g_{h_{-}},g_{e_{i}},g_{\tau-e_{i}})-g_{h_{-}}}{g_{h_{-}}g_{h_{+}}}\le\frac{3}{s}.
\end{aligned}
\]
For the first bound, relabel so that $0<g_{e_{1}}\le g_{e_{2}}\le g_{e_{3}}$.
The first row sum is at most $6/g_{\tau-e_{1}}\le9/s$. If $g_{e_{3}}\le s/2$,
the third row sum is at most $2s/(g_{e_{3}}g_{\tau-e_{3}})\le9/s$,
since $s/3\le g_{e_{3}}\le s/2$. If $g_{e_{3}}\ge s/2$, it is at
most $4/g_{e_{3}}\le8/s$. For the second row it suffices to check
\[
2s(g_{e_{1}}+2g_{e_{2}})\le9g_{e_{2}}(s-g_{e_{2}}).
\]
The right side minus the left side is $5sg_{e_{2}}-9g^{2}_{e_{2}}-2sg_{e_{1}}$.
If $g_{e_{2}}\le s/3$, use $g_{e_{1}}\le g_{e_{2}}$ to bound this
below by $3g_{e_{2}}(s-3g_{e_{2}})\ge0$. If $g_{e_{2}}\ge s/3$,
use $g_{e_{1}}\le s-2g_{e_{2}}$ to bound it below by 
\[
9g_{e_{2}}(s-g_{e_{2}})-2s^{2}=(3g_{e_{2}}-s)(2s-3g_{e_{2}})\ge0,
\]
where $g_{e_{2}}\le s/2$.

For the $S$ row, use 
\[
\sum_{i}\min(g_{h_{-}},g_{e_{i}},g_{\tau-e_{i}})\le\min(3g_{h_{-}},s).
\]
If $g_{h_{-}}\le s/3$, the desired quantity is at most $2/g_{h_{+}}\le3/s$;
if $g_{h_{-}}\ge s/3$, it is at most $1/g_{h_{-}}\le3/s$. The scaled
row estimate and \eqref{eq:4.9} now give 
\[
\mathcal{C}_{\tau}\ge-\frac{9}{s}\operatorname{Id}\ge-\frac{18}{11}\operatorname{Id}>-2\operatorname{Id}.
\]
If $S$ occurs with only one $L$ pair, the block has positive trace
$1/g_{h_{+}}$. The upper bound $\mathcal{C}\le2\operatorname{Id}$
from Theorem~\ref{theorem:4.1} therefore implies $\lambda_{\min}(\mathcal{C}_{\tau})>-2$
for this two-dimensional block. If no $L$ pair occurs, its single
entry is positive.

Finally suppose $S$ is absent. With at most two $L$ pairs the block
has trace zero. Proposition~\ref{proposition:4.2} gives the strict
upper bound $\lambda_{\max}(\mathcal{C}_{\tau})<2$ here: an equality
vector has torus weight twice a root, whereas $|\tau|^{2}=3$ is neither
$4$ nor $8$. Thus a nonempty block has spectrum $\{0\}$ or $\{-r,r\}$
with $0\le r<2$.

With all three $L$ pairs, the same scaling factors $U_{L_{i}}$ give
off-diagonal row sums at most $6/s$. To see this, relabel so that
$0<g_{e_{1}}\le g_{e_{2}}\le g_{e_{3}}$. The first and third row
sums are $4/g_{\tau-e_{1}}\le6/s$ and $2/g_{e_{3}}\le6/s$. For the
second row, the required inequality is 
\[
s(g_{e_{1}}+g_{e_{2}})\le3g_{e_{2}}(s-g_{e_{2}}).
\]
If $g_{e_{2}}\le s/3$, the difference between the right and left
sides is at least $g_{e_{2}}(s-3g_{e_{2}})\ge0$, using $g_{e_{1}}\le g_{e_{2}}$.
If $g_{e_{2}}\ge s/3$, it is at least $(s-g_{e_{2}})(3g_{e_{2}}-s)\ge0$,
using $g_{e_{1}}\le s-2g_{e_{2}}$. By \eqref{eq:4.9}, 
\[
\mathcal{C}_{\tau}\ge-\frac{6}{s}\operatorname{Id}\ge-\frac{4}{3}\operatorname{Id}>-2\operatorname{Id}.
\]
Thus $\mathcal{C}_{\tau}>-2\operatorname{Id}$ whenever $|\tau|^{2}=3$.

\subsubsection{Blocks whose sum is a short root}

Let $\tau$ be a short root and put $s=g_{\tau}>0$. Apply a Weyl-group
transformation so that $\tau=e_{1}$. Choose representatives $v_{1},v_{2},v_{3}$
of the three orthogonal short-root axes perpendicular to $\tau$,
oriented so that $g_{v_{i}}\ge0$. Put 
\[
M=\#\{i:g_{v_{i}}<s\},\qquad N=\#\{\text{short-short tangent decompositions of }\tau\}.
\]
Axes with zero metric coefficient are included in $M$. For an axis
with $0<g_{v_{i}}<s$, there is exactly one long-short tangent pair
$(v_{i},\tau-v_{i})$; if $g_{v_{i}}=0$, this pair is absent. With
the Chevalley normalization $z_{\alpha}=2/|\alpha|^{2}$, choose the
positive scaling factors 
\[
U_{\alpha\gamma}=\sqrt{\frac{w_{\alpha}w_{\gamma}}{2}}
\]
for the tangent pairs. Subscripts on $U$ denote these unordered root
pairs. Every pair in this block has distinct roots.

The canonical root count gives 
\begin{equation}
s\ge M+\frac{1}{2}+\frac{N}{2}.\label{eq:4.10}
\end{equation}
Indeed, long-root contributions to $2s$ cancel by reflection $\beta\mapsto\beta+2\tau$,
except for the $M$ tangent decompositions of $2\tau$, which contribute
$2M$, and nonnegative contributions from roots with Levi complements.
Among short roots, contributions with pairing $-1/2$ cancel with
$\beta+\tau$. The root $\tau$ contributes one, and each short-short
tangent decomposition of $\tau$ contributes one. This proves \eqref{eq:4.10}.
The inequality is strict if a short-short decomposition of $\tau$
has one summand in the Levi system: the other summand is tangent and
has inner product $1/2$ with $\tau$, giving an additional contribution
$1/4$ to $s$. We will use this strict inequality when exactly three
short-short pairs occur.

Fix a short-short row. Each other short-short pair corresponds to
a unique short-root axis orthogonal to $\tau$. Indeed, the summands
are half roots, and the other pair is obtained by changing one of
the last three signs, up to interchanging the summands. If both pairs
are tangent, then the coefficient of the corresponding axis root is
less than $s$, since it is the absolute difference of two metric
coefficients in $(0,s)$. The $N-1$ other short-short pairs therefore
correspond to distinct axes among the $M$ counted above, so 
\begin{equation}
N-1\le M.\label{eq:4.11}
\end{equation}
For each axis, group the corresponding short-short pair with the long-short
pair, when these pairs are tangent. We show that an axis with a second
tangent short-short pair contributes at most $3/s$ to the scaled
off-diagonal row sum. An axis without such a pair contributes at most
$1/s$.

To prove these estimates, fix one axis root $v$ with $0\le g_{v}<s$.
The two short-short pairs associated with this axis and the long-short
pair lie in a $C_{3}$ subsystem. Choose orthogonal vectors $f_{1},f_{2},f_{3}$
of squared length $1/2$ so that 
\[
\tau=f_{1}+f_{2},\qquad v=f_{1}-f_{2},
\]
and name the actual roots 
\[
\alpha=f_{2}-f_{3},\quad\gamma=f_{1}+f_{3},\quad\beta=f_{2}+f_{3},\quad\epsilon=f_{1}-f_{3}.
\]
Thus the three pairs are $(\alpha,\gamma)$, $(\beta,\epsilon)$,
and $(v,\tau-v)$. In particular, 
\[
\gamma=\beta+v,\qquad\epsilon=\alpha+v,\qquad\tau-v=\alpha+\beta.
\]
The vectors $f_{i}$ themselves are not roots. Change the sign of
$f_{3}$ if necessary so that $g_{\beta}\ge g_{\alpha}$ and $g_{\gamma}\ge g_{\epsilon}$.

When both short-short pairs are tangent, the $C_{3}$ curvature formula
gives 
\[
\begin{aligned} & |(\mathcal{C}_{\tau})_{\alpha\gamma,\beta\epsilon}|\frac{U_{\beta\epsilon}}{U_{\alpha\gamma}}+|(\mathcal{C}_{\tau})_{\alpha\gamma,v(\tau-v)}|\frac{U_{v(\tau-v)}}{U_{\alpha\gamma}}\\
 & \quad=\frac{s+g_{\epsilon}}{sg_{\gamma}}+\frac{g_{v}}{sg_{\gamma}}=\frac{3}{s}-\frac{2(g_{\gamma}-g_{\epsilon})}{sg_{\gamma}}\le\frac{3}{s},\\[3pt]
 & |(\mathcal{C}_{\tau})_{\beta\epsilon,\alpha\gamma}|\frac{U_{\alpha\gamma}}{U_{\beta\epsilon}}+|(\mathcal{C}_{\tau})_{\beta\epsilon,v(\tau-v)}|\frac{U_{v(\tau-v)}}{U_{\beta\epsilon}}\\
 & \quad=\frac{g_{\alpha}(s+g_{\epsilon})}{sg_{\beta}g_{\epsilon}}+\frac{g_{v}}{sg_{\epsilon}}=\frac{3}{s}-\frac{2(g_{\beta}-g_{\alpha})}{sg_{\beta}}\le\frac{3}{s}.
\end{aligned}
\]
If $g_{v}=0$, omit the terms indexed by $(v,\tau-v)$; their rational
expressions are zero, and this pair is not a matrix row.

When only one short-short pair and the long-short pair are tangent,
the short-short pair is $(\beta,\epsilon)$ and $g_{\alpha}\le0<g_{\epsilon}$.
The invariant curvature formula gives 
\[
|(\mathcal{C}_{\tau})_{\beta\epsilon,v(\tau-v)}|\frac{U_{v(\tau-v)}}{U_{\beta\epsilon}}=\frac{g_{\tau-v}}{sg_{\beta}}\le\frac{1}{s},
\]
since $g_{\tau-v}=g_{\alpha}+g_{\beta}\le g_{\beta}$. At $g_{\alpha}=0$
this equals the boundary value of the preceding grouped expression.
If $g_{v}\ge s$, neither a second tangent short-short pair nor a
tangent long-short pair occurs, so the axis contributes zero. These
formulas follow from the $C_{3}$ root matrices in Lemma~\ref{lemma:4.7}.
The squared root lengths are half those in that normalization, so
the curvature entries in the unit pair basis are halved for the same
metric coefficients. The factors $z_{\mathrm{short}}=2$ and $z_{\mathrm{long}}=1$
give the scaling factors $U$ above.

The diagonal of a short-short row is $-1/s$. Summing over the axes,
its negative diagonal plus scaled off-diagonal absolute sum is at
most 
\[
\frac{1+3(N-1)+(M-N+1)}{s}=\frac{M+2N-1}{s}.
\]
By \eqref{eq:4.10} and \eqref{eq:4.11}, the diagonal entry of $2\operatorname{Id}+\mathcal{C}_{\tau}$
exceeds the scaled off-diagonal sum by at least 
\[
2-\frac{M+2N-1}{s}\ge\frac{M-N+2}{s}\ge\frac{1}{s}>0.
\]
Thus every short-short row satisfies strict scaled diagonal dominance.

We next consider the long-short rows, with the same scaling factors
$U$. Write 
\[
L_{i}=(v_{i},\tau-v_{i})\qquad(0<g_{v_{i}}<s),
\]
and let $L$ be their number, so $L\le M$. Their diagonal entries
are $-2/s$. The label $i$ records the axis root, even when two axes
have the same metric coefficient. For distinct such rows, the scaled
entry is 
\begin{equation}
s|(\mathcal{C}_{\tau})_{L_{i},L_{j}}|\frac{U_{L_{j}}}{U_{L_{i}}}=2\min\left\{ \frac{g_{v_{j}}}{g_{v_{i}}},\frac{g_{\tau-v_{j}}}{g_{\tau-v_{i}}}\right\} \le2.\label{eq:4.12}
\end{equation}
Compute in the $B_{3}$ subsystem containing the two pairs. The absolute
value of the curvature component on the root vectors is 
\[
\frac{2}{s}\min(g_{v_{j}}g_{\tau-v_{i}},g_{v_{i}}g_{\tau-v_{j}});
\]
normalizing the pair vectors and multiplying by the ratio of scaling
factors gives \eqref{eq:4.12}.

Index the short-short pairs by sign triples: 
\[
\alpha_{\boldsymbol{\varepsilon}}=\frac{1}{2}\left(\tau+\sum^{3}_{j=1}\varepsilon_{j}v_{j}\right),\qquad\gamma_{\boldsymbol{\varepsilon}}=\tau-\alpha_{\boldsymbol{\varepsilon}},\qquad\varepsilon_{j}\in\{-1,1\}.
\]
Write $S_{[\boldsymbol{\varepsilon}]}$ for the unordered pair $(\alpha_{\boldsymbol{\varepsilon}},\gamma_{\boldsymbol{\varepsilon}})$,
where simultaneous reversal of all three signs gives the same pair.
These are four distinct geometric pairs, before omitting those that
are not tangent.

For a fixed row $L_{i}$, choose the representative with $\varepsilon_{i}=1$.
Then 
\[
\begin{aligned}2g_{\alpha_{\boldsymbol{\varepsilon}}} & =s+g_{v_{i}}+\sum_{j\ne i}\varepsilon_{j}g_{v_{j}},\\
2g_{\gamma_{\boldsymbol{\varepsilon}}} & =s-g_{v_{i}}-\sum_{j\ne i}\varepsilon_{j}g_{v_{j}}.
\end{aligned}
\]
The pair is tangent precisely when 
\[
-s-g_{v_{i}}<\sum_{j\ne i}\varepsilon_{j}g_{v_{j}}<g_{\tau-v_{i}}.
\]
For each existing row $L_{i}$, define the continuous nonnegative
function 
\begin{equation}
F_{i}(r)=\max\left\{ 0,\min\left\{ \frac{s+g_{v_{i}}+r}{g_{v_{i}}},\frac{s-g_{v_{i}}-r}{g_{\tau-v_{i}}}\right\} \right\} .\label{eq:4.13}
\end{equation}
In \eqref{eq:4.13}, the argument $r$ is a real number, not an index
for a root pair. Substitution in the local $C_{3}$ curvature formula
gives 
\[
s|(\mathcal{C}_{\tau})_{L_{i},S_{[\boldsymbol{\varepsilon}]}}|\frac{U_{S_{[\boldsymbol{\varepsilon}]}}}{U_{L_{i}}}=F_{i}\left(\sum_{j\ne i}\varepsilon_{j}g_{v_{j}}\right)
\]
when the pair is tangent. For an absent pair its contribution is zero,
as given by $F_{i}$; no matrix entry or scaling factor is assigned
to it. The two affine functions in \eqref{eq:4.13} meet at the value
$2$, so $0\le F_{i}\le2$. The two expressions according to the ordering
of the metric coefficients give the absolute curvature component 
\[
\frac{2}{s}\min\left(g_{v_{i}}g_{\gamma_{\boldsymbol{\varepsilon}}},g_{\tau-v_{i}}g_{\alpha_{\boldsymbol{\varepsilon}}}\right).
\]
Normalization and scaling give $2\min(g_{\alpha_{\boldsymbol{\varepsilon}}}/g_{v_{i}},g_{\gamma_{\boldsymbol{\varepsilon}}}/g_{\tau-v_{i}})$,
which is the value of $F_{i}$ above for a tangent pair.

For $r\ge0$, the paired-argument identity is 
\begin{equation}
F_{i}(r)+F_{i}(-r)=\begin{cases}
2, & 0\le r\le g_{\tau-v_{i}},\\
\dfrac{s+g_{v_{i}}-r}{g_{v_{i}}}, & g_{\tau-v_{i}}<r<s+g_{v_{i}},\\
0, & r\ge s+g_{v_{i}}.
\end{cases}\label{eq:4.14}
\end{equation}
For $\{i,j,k\}=\{1,2,3\}$, the total short-short contribution to
row $L_{i}$, multiplied by $s$, is therefore 
\[
\sum_{\varepsilon_{j}=\pm1}\sum_{\varepsilon_{k}=\pm1}F_{i}(\varepsilon_{j}g_{v_{j}}+\varepsilon_{k}g_{v_{k}})\le4.
\]
The sum is over sign choices, so equal arguments retain their multiplicity,
including when one or both coefficients vanish. If exactly three short-short
pairs are tangent, one of the sign groups $\{(+,+),(-,-)\}$ and $\{(+,-),(-,+)\}$
has two tangent pairs and the other has exactly one. The latter group
contributes strictly less than $2$ unless the absent pair has a root
with zero metric coefficient: this is the equality boundary $r=g_{\tau-v_{i}}$
in \eqref{eq:4.14}. Thus the short-short contribution to row $L_{i}$
is strictly less than $4/s$ except in that boundary case.

We now apply \eqref{eq:4.3}, adjusting the scaling factors when $N\le2$.

\smallskip{}
\noindent\emph{At most two short-short pairs.} If $N=0$, \eqref{eq:4.12}
bounds each row's negative diagonal plus absolute off-diagonal sum
by $2L/s\le2M/s<2$. If $1\le N\le2$, replace $U_{L_{i}}$ by $2U_{L_{i}}$
on every long-short pair and keep the original factors on the short-short
pairs. The scaled contributions of short-short columns to a long-short
row are halved; since $F_{i}\le2$, each contributes at most $1/s$.
Thus the negative diagonal plus scaled off-diagonal sum in each long-short
row is at most 
\[
\frac{2L+N}{s}\le\frac{2M+N}{s}<2.
\]
For a short-short row, an axis without a second tangent short-short
pair contributes at most $2/s$ after this change of scaling factors.
An axis with such a pair contributes strictly less than $4/s$: its
original grouped contribution is at most $3/s$, and the extra long-short
contribution is strictly less than $1/s$. Indeed, in the local root
notation above, both $g_{v}/g_{\gamma}$ and $g_{v}/g_{\epsilon}$
are strictly less than one, because $\gamma-v=\beta$ and $\epsilon-v=\alpha$
are tangent roots. If $g_{v}=0$, the long-short pair is absent and
its contribution is zero. There are $N-1$ axes with a second tangent
short-short pair and at most $M-N+1$ without one. For $N=1$, the
row bound is $(1+2M)/s<2$ by \eqref{eq:4.10}. For $N=2$, it is
strictly less than 
\[
\frac{1+4+2(M-1)}{s}=\frac{2M+3}{s}\le2.
\]
Hence all these blocks satisfy strict scaled diagonal dominance of
$\mathcal{C}_{\tau}+2\operatorname{Id}$.

\smallskip{}
\noindent\emph{Three or four short-short pairs.} Use the original
scaling factors $U$. The required strict diagonal dominance for the
short-short rows was proved above. Here $M=3$. To see this, relabel
the axes so that $g_{v_{1}}\ge g_{v_{2}}\ge g_{v_{3}}\ge0$. The two
largest of the four absolute signed sums defining the half-root pairs
are $g_{v_{1}}+g_{v_{2}}+g_{v_{3}}$ and $g_{v_{1}}+g_{v_{2}}-g_{v_{3}}$.
If three pairs are tangent, the latter is less than $s$, and hence
$g_{v_{1}}<s$; the assertion for four pairs follows as well. Every
long-short row is bounded by $(2L+4)/s\le10/s$. For $N=4$, \eqref{eq:4.10}
gives $2s\ge11$, so the dominance is strict. For $N=3$, it gives
$2s\ge10$. If a half root has zero metric coefficient, the strict
version of \eqref{eq:4.10} gives $2s>10$. Otherwise \eqref{eq:4.14}
makes the short-short contribution strictly less than $4/s$, and
again the dominance is strict.

Thus \eqref{eq:4.3} gives $\mathcal{C}_{\tau}>-2\operatorname{Id}$
when $\tau$ is a short root.

All blocks satisfy $\mathcal{C}_{\tau}>-2\operatorname{Id}$. Since
there are finitely many blocks, this proves Theorem~\ref{theorem:4.9}.

\subsection{The spectrum in type $G_{2}$}

Let $a$ be the short simple root and $b$ the long simple root, normalized
by 
\[
(a,a)=\frac{2}{3},\quad(b,b)=2,\quad(a,b)=-1.
\]
The positive roots are 
\[
a,\quad b,\quad c=a+b,\quad d=2a+b,\quad e=3a+b,\quad f=3a+2b.
\]
The three nonempty painted subsets are $\{a\}$, $\{b\}$, and $\{a,b\}$.
We compute the root-sum blocks for each of the corresponding flags.

\begin{proposition}\label{proposition:4.10} For the canonical metric
with $\Ric=2g$, the least eigenvalues of the Calabi operator on the
three type-$G_{2}$ flags are 
\[
\begin{array}{c|ccc}
\textit{painted roots} & \{a\} & \{b\} & \{a,b\}\\
\hline \lambda_{\min}(\mathcal{C}) & -\frac{6}{5} & -\frac{8}{9} & -\frac{3}{2}.
\end{array}
\]
In particular, all three least eigenvalues are strictly greater than
$-2$. \end{proposition} 
\begin{proof}
Choose a Chevalley basis with $B(E_{\alpha},E_{-\alpha})=z_{\alpha}=\frac{2}{(\alpha,\alpha)}$
and $[E_{\alpha},E_{-\alpha}]=z_{\alpha}H_{\alpha}$, where $[H_{\alpha},E_{\beta}]=(\alpha,\beta)E_{\beta}$.
Thus 
\[
z_{a}=z_{c}=z_{d}=3,\quad z_{b}=z_{e}=z_{f}=1.
\]
The only nonzero brackets between positive root vectors, apart from
antisymmetry, can be chosen as 
\begin{equation}
\begin{aligned}[E_a,E_b] & =E_{c}, & [E_{a},E_{c}] & =2E_{d}, & [E_{a},E_{d}] & =3E_{e},\\{}
[E_{b},E_{e}] & =E_{f}, & [E_{c},E_{d}] & =-3E_{f}.
\end{aligned}
\label{eq:4.15}
\end{equation}
The absolute values are determined by the root strings. The first
four brackets fix the signs, and the last follows from $0=[E_{a},[E_{b},E_{d}]]=[E_{c},E_{d}]+3[E_{b},E_{e}]$.
All mixed-root constants are fixed by antisymmetry, the opposite-root
brackets above, and invariance of $B$: 
\[
N_{-\alpha,-\beta}=-N_{\alpha\beta},\quad N_{\alpha\beta}z_{\alpha+\beta}=N_{\beta,-\alpha-\beta}z_{\alpha}\quad\text{when }\alpha+\beta\text{ is a root}.
\]
These relations determine the brackets used below.

Write $\alpha\gamma$ for the unit symmetric pair of roots $\alpha,\gamma$
defined in Section~\ref{sec:2}. All entries below are obtained from
\eqref{eq:2.6} and \eqref{eq:2.8}, with the normalization in \eqref{eq:2.5}.
The factor $\sqrt{p_{\alpha\gamma}p_{\beta\epsilon}}$ in that formula
is $1$ when both pairs are repeated, $2$ when both pairs consist
of distinct roots, and $\sqrt{2}$ when exactly one pair is repeated.

The canonical metric coefficients are listed below. A zero entry corresponds
to a Levi root and hence to an omitted tangent direction: 
\begin{equation}
\begin{array}{c|c|cccccc}
\text{painted roots} & \delta & g_{a} & g_{b} & g_{c} & g_{d} & g_{e} & g_{f}\\
\hline \{a\} & 10a+5b & \frac{5}{6} & 0 & \frac{5}{6} & \frac{5}{3} & \frac{5}{2} & \frac{5}{2}\\
\{b\} & 9a+6b & 0 & \frac{3}{2} & \frac{3}{2} & \frac{3}{2} & \frac{3}{2} & 3\\
\{a,b\} & 10a+6b & \frac{1}{3} & 1 & \frac{4}{3} & \frac{5}{3} & 2 & 3
\end{array}\label{eq:4.16}
\end{equation}
For each painted subset, $\delta$ is the sum of the positive tangent
roots, and $g_{\alpha}=(\delta,\alpha)/2$.

The sums with more than one unordered decomposition into positive
roots are 
\begin{equation}
b+d=2c,\quad b+e=c+d,\quad a+f=c+e=2d.\label{eq:4.17}
\end{equation}
For a partial flag, we omit the pairs containing a Levi root. Every
remaining pair whose sum occurs only once gives a scalar block.

We first compute the entries for arbitrary nonnegative simple-root
coefficients $g_{a},g_{b}$, omitting roots with zero metric coefficient.
The other coefficients are determined by additivity, and the squared
norms of the root vectors are 
\[
(w_{a},w_{b},w_{c},w_{d},w_{e},w_{f})=(3g_{a},g_{b},3g_{c},3g_{d},g_{e},g_{f}).
\]
Before normalization, the curvature entries in these blocks are 
\[
\begin{array}{c|c|c|c}
\text{entry} & \text{value} & \text{entry} & \text{value}\\
\hline R_{b\bar{b}d\bar{d}} & 0 & R_{c\bar{c}c\bar{c}} & 6g_{c}\\
R_{b\bar{b}e\bar{e}} & -\frac{g_{e}g_{b}}{g_{f}} & R_{c\bar{c}d\bar{d}} & 3g_{c}+\frac{9g^{2}_{c}}{g_{f}}\\
R_{a\bar{a}f\bar{f}} & 0 & R_{c\bar{c}e\bar{e}} & 0\\
R_{d\bar{d}d\bar{d}} & 6g_{d} & R_{b\bar{c}d\bar{c}} & 6g_{b}\\
R_{b\bar{c}e\bar{d}} & \frac{3g_{d}g_{b}}{g_{f}} & R_{a\bar{c}f\bar{e}} & -3g_{a}\\
R_{a\bar{d}f\bar{d}} & 6g_{a} & R_{c\bar{d}e\bar{d}} & 6g_{c}
\end{array}
\]
For the five off-diagonal entries, the bracket pairings in Proposition~\ref{proposition:2.1}
are 
\[
\begin{array}{c|ccccc}
\text{entry} & bd,cc & be,cd & af,ce & af,dd & ce,dd\\
\hline N_{1} & 6 & 3 & -3 & 6 & 6\\
N_{2} & 6 & 0 & -3 & 6 & 6.
\end{array}
\]
These values follow from \eqref{eq:4.15} and invariance of $B$.
Substitution into \eqref{eq:2.6} gives the off-diagonal entries in
the preceding table. The diagonal values use $(b,d)=(a,f)=(c,e)=0$,
$(b,e)=-1$, and $(c,d)=\frac{1}{3}$; in particular $R_{b\bar{b}e\bar{e}}=-g_{e}+\frac{g^{2}_{e}}{g_{f}}=-\frac{g_{e}g_{b}}{g_{f}}$.
When $g_{a}=0$ or $g_{b}=0$, the same formulas apply to the remaining
pairs. If a root difference becomes a Levi root, the corresponding
metric coefficients are equal, so the two cases of \eqref{eq:2.8}
agree. Substituting \eqref{eq:4.16} and applying \eqref{eq:2.5}
gives the matrices below.

The scalar entries obtained from \eqref{eq:2.8} and \eqref{eq:2.5}
are 
\[
\begin{array}{c|c@{\qquad}c|c}
aa & \frac{2}{3g_{a}} & bb & \frac{2}{g_{b}}\\
ab & -\frac{2}{g_{c}} & ac & \frac{2(2g_{a}-g_{b})}{3g_{c}g_{d}}\\
ad & \frac{2(6g_{a}+g_{b})}{3g_{d}g_{e}} & bc & \frac{2}{g_{c}}\\
ae,\ de,\ ee & \frac{2}{g_{e}} & bf,\ cf,\ df,\ ef,\ ff & \frac{2}{g_{f}}.
\end{array}
\]
Pairs involving a root with zero metric coefficient are omitted. When
$b$ is omitted, $cc$ and $cd$ also become scalar blocks; their
values follow from the preceding table of curvature components.

\smallskip{}
 \emph{The flag with painted short root $a$.} Here $\Delta^{+}_{\Phi_{T}}=\{a,c,d,e,f\}$.
The only nonscalar block has sum $a+f=c+e=2d$. In the ordered unit
pair basis $(af,ce,dd)$, its matrix is 
\[
\mathcal{C}_{a+f}=\frac{1}{5}\begin{pmatrix}0 & -4 & 2\sqrt{2}\\
-4 & 0 & 2\sqrt{2}\\
2\sqrt{2} & 2\sqrt{2} & 2
\end{pmatrix}.
\]
Every other pair, namely 
\[
aa,\ ac,\ cc,\ ad,\ cd,\ ae,\ cf,\ de,\ df,\ ee,\ ef,\ ff,
\]
gives the scalar $\frac{4}{5}$ by \eqref{eq:2.8} and \eqref{eq:2.5}.
Set $v=\frac{1}{\sqrt{5}}(\sqrt{2},\sqrt{2},-1)^{\mathsf{T}}$. The
displayed matrix equals $(\frac{4}{5})\operatorname{Id}-2vv^{\mathsf{T}}$,
and $\|v\|=1$. The symmetric square has dimension $15$, and the
complete spectrum is 
\[
\frac{4}{5}\quad\text{with multiplicity }14,\quad-\frac{6}{5}\quad\text{with multiplicity }1.
\]
Hence $\lambda_{\min}(\mathcal{C})=-\frac{6}{5}>-2$, attained on
the line spanned by $v$ in this block.

\smallskip{}
 \emph{The flag with painted long root $b$.} Here $\Delta^{+}_{\Phi_{T}}=\{b,c,d,e,f\}$.
The blocks with sums $b+d$, $b+e$, and $c+e$ have ordered unit
pair bases $(bd,cc)$, $(be,cd)$, and $(ce,dd)$, respectively. Their
matrices are 
\[
\begin{aligned}\mathcal{C}_{b+d}=\mathcal{C}_{c+e} & =\frac{1}{9}\begin{pmatrix}0 & 4\sqrt{6}\\
4\sqrt{6} & 4
\end{pmatrix},\\
\mathcal{C}_{b+e} & =\frac{1}{9}\begin{pmatrix}-6 & 6\\
6 & 10
\end{pmatrix}.
\end{aligned}
\]
The scalar blocks are 
\[
\begin{array}{c|l}
\frac{4}{3} & bb,\ bc,\ de,\ ee\\
\frac{2}{3} & bf,\ cf,\ df,\ ef,\ ff.
\end{array}
\]
Both displayed matrices have trace $\frac{4}{9}$ and determinant
$-\frac{32}{27}$. Their characteristic polynomial is therefore 
\[
(\lambda-\frac{4}{3})(\lambda+\frac{8}{9}).
\]
This gives the complete spectrum 
\[
\frac{4}{3}\ (7\text{ times}),\quad\frac{2}{3}\ (5\text{ times}),\quad-\frac{8}{9}\ (3\text{ times}).
\]
Consequently $\lambda_{\min}(\mathcal{C})=-\frac{8}{9}>-2$.

\smallskip{}
 \emph{The full flag.} All six roots are tangent. The blocks with
sums $b+d$, $b+e$, and $a+f$ have ordered unit pair bases $(bd,cc)$,
$(be,cd)$, and $(af,ce,dd)$, respectively. Their matrices are 
\[
\begin{aligned}\mathcal{C}_{b+d} & =\begin{pmatrix}0 & \frac{3\sqrt{10}}{10}\\
\frac{3\sqrt{10}}{10} & \frac{1}{2}
\end{pmatrix},\\
\mathcal{C}_{b+e} & =\begin{pmatrix}-\frac{2}{3} & \frac{\sqrt{10}}{6}\\
\frac{\sqrt{10}}{6} & \frac{14}{15}
\end{pmatrix},\\
\mathcal{C}_{a+f} & =\begin{pmatrix}0 & -\frac{\sqrt{6}}{6} & \frac{2\sqrt{6}}{15}\\
-\frac{\sqrt{6}}{6} & 0 & \frac{4}{5}\\
\frac{2\sqrt{6}}{15} & \frac{4}{5} & \frac{2}{5}
\end{pmatrix}.
\end{aligned}
\]
The remaining scalar blocks are 
\[
\begin{array}{c|l}
2 & aa,\ bb\\
-\frac{3}{2} & ab\\
-\frac{1}{10} & ac\\
\frac{3}{5} & ad\\
\frac{3}{2} & bc\\
1 & ae,\ de,\ ee\\
\frac{2}{3} & bf,\ cf,\ df,\ ef,\ ff.
\end{array}
\]
For example, $c-a=b$ is tangent, $(a,c)=-\frac{1}{3}$, and $N_{ac}=2$.
Thus \eqref{eq:2.8} gives the entry for $ac$: 
\[
R_{a\bar{a}c\bar{c}}=\frac{1}{3}\cdot9\left(-\frac{1}{3}\right)+\frac{(\frac{1}{3})^{2}}{\frac{5}{3}}\cdot3\cdot2^{2}=-\frac{1}{5}.
\]
Since $w_{a}=1$ and $w_{c}=4$, \eqref{eq:2.5} gives $(\mathcal{C}_{a+c})_{ac,ac}=\frac{2(-\frac{1}{5})}{4}=-\frac{1}{10}$.
Similarly, invariance of $B$ and \eqref{eq:4.15} give $N_{a,-c}=-3$
and $N_{f,-e}=1$. As $f-c=d\in\Delta^{+}_{\Phi_{T}}$ and $a+f$
is not a root, 
\[
R_{a\bar{c}f\bar{e}}=g_{a}z_{b}(-3)(1)=-1,\quad(\mathcal{C}_{a+f})_{af,ce}=\frac{-2}{\sqrt{1\cdot3\cdot4\cdot2}}=-\frac{\sqrt{6}}{6}.
\]
To compare the nonscalar eigenvalues with $-\frac{3}{2}$, add $(\frac{3}{2})\operatorname{Id}$
to each block. The resulting matrices are strictly diagonally dominant
with positive diagonal. For the first block this uses $\frac{3\sqrt{10}}{10}<\frac{6}{5}<\frac{3}{2}$.
For the second, its first diagonal is $\frac{5}{6}>\frac{\sqrt{10}}{6}$,
and the second is $\frac{73}{30}>\frac{\sqrt{10}}{6}$. For the last
block, the three sums of off-diagonal absolute values are 
\[
\frac{3\sqrt{6}}{10}<\frac{9}{10},\quad\frac{\sqrt{6}}{6}+\frac{4}{5}<\frac{13}{10},\quad\frac{2\sqrt{6}}{15}+\frac{4}{5}<\frac{6}{5},
\]
whereas its shifted diagonal is $(\frac{3}{2},\frac{3}{2},\frac{19}{10})$.
These estimates use $\sqrt{10}<4$ and $\sqrt{6}<3$. A Hermitian
matrix with positive diagonal and strict diagonal dominance is positive
definite. Thus every nonscalar eigenvalue is strictly greater than
$-\frac{3}{2}$. The scalar block $ab$ attains $-\frac{3}{2}$, so
\[
\lambda_{\min}(\mathcal{C})=-\frac{3}{2}>-2.
\]
\end{proof}

\subsection{Products and completion of the proof}
\begin{proof}[Proof of Theorem~\ref{theorem:1.2}]
Theorem~\ref{theorem:4.1} proves $\mathcal{C}\leq(4/k_{*})\operatorname{Id}\leq2\operatorname{Id}$
in every Lie type. For a simple group, the strict lower bound follows
from Theorem~\ref{theorem:4.4} in types $A$, $D$, and $E$, Theorem~\ref{theorem:4.8}
in types $B$ and $C$, Theorem~\ref{theorem:4.9} in type $F_{4}$,
and Proposition~\ref{proposition:4.10} in type $G_{2}$.

A generalized flag for a semisimple group is a product of the flags
of its simple factors, and its canonical metric is their product,
with Einstein constant two on each factor. On the symmetric square
of the product tangent space, $\mathcal{C}$ is the direct sum of
the operators on the factors and the zero operator on mixed-factor
tensors. Its spectrum therefore lies in $(-2,2]$. The sharp upper
bound and the description of the eigenvalue-$2$ space were proved
in Theorem~\ref{theorem:4.1} and Proposition~\ref{proposition:4.2}. 
\end{proof}

\section{Counterexamples to homogeneity}

\label{sec:5}

We prove Theorem~\ref{theorem:1.3} by constructing three metrics
on projective bundles. The first metric is Kähler-Einstein on a complex
13-fold. With $\Ric=g$, its Calabi operator has spectrum in $[-\frac{1}{2},\frac{11}{12}]$,
and both endpoints are included. Equivalently, after rescaling to
$\Ric=2g$, the spectral interval is $[-1,\frac{11}{6}]\subset(-2,2)$.
The other two metrics are restrictions of a product Fubini-Study metric
to an incidence submanifold. They are not Einstein, so the Calabi
spectrum must be combined with the individual Ricci eigenvalues to
estimate $\CQB$ and $\dCQB$.

The one-variable construction of the Einstein metric belongs to the
projective-bundle method of Koiso and Sakane \cite{KS88}. It also
fits the admissible setting of Apostolov, Calderbank, Gauduchon, and
Tønnesen-Friedman \cite{ACGTF08}, with two rank-two projectively
flat summands. We verify the potential, the Einstein equation, and
the extension across the two exceptional strata directly. The Calabi
spectral bounds below are the part needed to address the homogeneity
conjecture. For the 5-fold, the underlying nonhomogeneous complex
manifold and its positive orthogonal Ricci curvature were already
established by Ni, Wang, and Zheng \cite[Example~1.4 and p.~837]{NWZ21};
here we establish positivity of the full $\CQB$ and $\dCQB$ forms
for arbitrary maps.

\subsection{A Kähler-Einstein counterexample.}

Let 
\[
\begin{aligned}B & =\mathbb{P}^{5}\times\mathbb{P}^{5},\quad L_{1}=\mathcal{O}_{B}(-1,0),\quad L_{2}=\mathcal{O}_{B}(0,-1),\\
M & =\mathbb{P}_{\mathrm{lines}}(L^{\oplus2}_{1}\oplus L^{\oplus2}_{2}).
\end{aligned}
\]
Here $\mathcal{O}_{B}(-1,0)$ is the pullback of the tautological
line bundle from the first factor, and $\mathcal{O}_{B}(0,-1)$ is
its counterpart from the second factor. The fiber of $M\to B$ is
$\mathbb{P}^{3}$, so $\dim_{\mathbb{C}}M=13$. Since the base is
compact and the fibers are projective spaces, $M$ is compact.

The two rank-two summands distinguish two disjoint submanifolds, namely
their projectivizations. Between these submanifolds, a line has nonzero
components in both summands. The ratio of their squared norms supplies
the single real variable used in the construction. The main point
in extending the metric is that an entire complex two-dimensional
block vanishes at either end. Checking only a radial coordinate would
therefore be insufficient; we check the potential in the original
block coordinates.

\subsubsection{The profile and the global metric.}

Give both factors of $B=\mathbb{P}^{5}\times\mathbb{P}^{5}$ the Fubini-Study
potentials $\phi_{1},\phi_{2}$, normalized to have holomorphic sectional
curvature 2. Give $L_{1},L_{2}$ the induced tautological metrics:
their standard local frames have squared norms $e^{\phi_{1}}$ and
$e^{\phi_{2}}$.

For $0<x<4$ define 
\begin{equation}
p(x)=x(4-x)(4+x)^{5}(8-x)^{5},\quad\Theta(x)=\frac{{\displaystyle \int^{x}_{0}(2-t)p(t)\,dt}}{p(x)}.\label{eq:5.1}
\end{equation}
The symmetry $p(4-x)=p(x)$ implies $\int^{4}_{0}(2-t)p(t)\,dt=0$.
The numerator is positive on $(0,4)$: it increases on $(0,2)$ and
then decreases to zero. Thus $\Theta>0$ on this interval. Taylor
expansions at the endpoints are 
\[
\begin{aligned}\Theta(x) & =x+O(x^{2}) &  & \text{as }x\to0,\\
\Theta(x) & =(4-x)+O((4-x)^{2}) &  & \text{as }x\to4.
\end{aligned}
\]
Moreover, $\Theta$ extends real analytically to $[0,4]$. Indeed,
$p(x)=cx+O(x^{2})$ near zero for some $c>0$, while the numerator
is $cx^{2}+O(x^{3})$. The apparent singularity of their ratio is
thus removable. The same argument at the other end follows either
by direct expansion or by the symmetry $\Theta(4-x)=\Theta(x)$. Define
$x=x(s)$ and $f=f(s)$ by 
\[
\frac{dx}{ds}=\Theta(x),\quad x(0)=2,\quad f'(s)=x(s),\quad f(0)=0.
\]
Equivalently, $s(x)=\int^{x}_{2}\frac{dt}{\Theta(t)}$. This function
is strictly increasing, tends to $-\infty$ at zero, and tends to
$+\infty$ at four. Its inverse therefore defines $x(s)$ on all of
$\mathbb{R}$.

Write a nonzero vector in $L^{\oplus2}_{1}\oplus L^{\oplus2}_{2}$
as $(a,b)$, and put $U=\lVert a\rVert^{2}$, $V=\lVert b\rVert^{2}$.
On $UV>0$ use the potential 
\[
\widetilde{\Phi}=4(\phi_{1}+\phi_{2})+4\log V+f\left(\log\frac{U}{V}\right).
\]
A change of local representative $(a,b)\mapsto\lambda(a,b)$, with
$\lambda$ nowhere vanishing and holomorphic, adds the pluriharmonic
term $4\log|\lambda|^{2}$. Changes of the base potentials are pluriharmonic,
and $U,V$ are the intrinsic squared norms of the block vectors. Hence
the local forms $\sqrt{-1}\partial\bar{\partial}\widetilde{\Phi}$
patch on the projective bundle wherever $UV>0$.

Here is a direct check of extension across the two remaining strata.
Near $x=0$, the function $\frac{1}{\Theta(x)}-\frac{1}{x}$ is analytic.
Thus $s-\log x$ is analytic in $x$, and $y=e^{s}$ is analytic in
$x$ with a simple zero and positive derivative at zero. The analytic
inverse function theorem lets us regard $x$ as an analytic function
of $y$, again with a simple zero and positive derivative. Consequently
$\frac{df}{dy}=\frac{x}{y}$ extends analytically and is positive
at $y=0$. Integrating gives 
\[
f(s)=F(e^{s}),\quad F'(0)>0,
\]
where $F$ is analytic near zero. The same function describes the
other endpoint. Indeed, the symmetry $\Theta(4-x)=\Theta(x)$ and
the initial value $x(0)=2$ give $x(-s)=4-x(s)$ by uniqueness of
the differential equation. Integrating and using $f(0)=0$ gives $f(s)=4s+f(-s)$.
Hence, for large positive $s$, 
\[
f(s)=4s+F(e^{-s}).
\]
Thus the potential near $U=0$ is $4(\phi_{1}+\phi_{2})+4\log V+F(\frac{U}{V})$,
and near $V=0$ it is $4(\phi_{1}+\phi_{2})+4\log U+F(\frac{V}{U})$.
These expressions are smooth in the two normal block coordinates,
including where their projective directions cease to be defined. For
example, near $U=0$ choose a nonzero component of $b$ as the projective
affine denominator. Then $V$ is a smooth positive function and $U$
is a positive Hermitian quadratic form in the two components of $a$.
The normal Hessian of $F(\frac{U}{V})$ at $a=0$ is $F'(0)$ times
that quadratic form divided by $V$; it is positive definite. Its
mixed Hessian with directions tangent to the stratum vanishes at $a=0$.
The restriction to the stratum has coefficients $4,8,4$ on the two
base factors and the surviving projective line. The argument near
$V=0$ is identical, with the two base factors interchanged. This
proves positive definiteness at both strata in genuine smooth coordinates.

For the positivity check and subsequent computations, use the dense
chart $a=w(1,u)$, $b=(1,v)$, with $w\ne0$. Write $\psi_{u}=\log(1+|u|^{2})$
and $\psi_{v}=\log(1+|v|^{2})$. The potential becomes 
\begin{align}
\Phi & =4\phi_{1}+8\phi_{2}+4\psi_{v}+f(s),\nonumber \\
s & =\log|w|^{2}+\phi_{1}-\phi_{2}+\psi_{u}-\psi_{v}.\label{eq:5.2}
\end{align}
The connection-adapted orthogonal decomposition has four Fubini-Study
blocks of dimensions $5,5,1,1$ and a radial complex line. The first
two blocks are tangent to the base factors; the last two describe
the directions $[a]$ and $[b]$ inside their rank-two summands. The
remaining coordinate $w$ describes their relative magnitude and phase.
In this construction, $h_{j}$ denotes a horizontal metric coefficient;
it is distinct from the root coefficient $g_{\alpha}$ used above.
Write these coefficients and their derivatives with respect to $x$
as 
\[
\begin{aligned}(h_{1},h_{2},h_{3},h_{4}) & =(4+x,8-x,x,4-x),\\
(\varepsilon_{1},\varepsilon_{2},\varepsilon_{3},\varepsilon_{4}) & =(1,-1,1,-1).
\end{aligned}
\]
If $g_{\mathrm{FS},j}$ is the Fubini-Study metric on the $j$th block,
direct differentiation gives 
\[
g=\sum^{4}_{j=1}h_{j}g_{\mathrm{FS},j}+\Theta\,\partial s\otimes\bar{\partial}s.
\]
In the resulting adapted coframe, the radial coefficient is $\Theta$,
and the coefficient of $d\log w$ in $\partial s$ is one. Replacing
$d\log w$ by $\partial s$ absorbs the mixed horizontal-radial entries
of the metric. This is the only change of coframe used here. All these
coefficients are positive for $0<x<4$. At either endpoint, the noncollapsing
coefficients are $4,8,4$, and the normal block is positive because
$F'(0)>0$. The patched form consequently defines a smooth Kähler
metric on all of $M$.

\subsubsection{Verification of the Einstein equation.}

The determinant in the coordinates of \eqref{eq:5.2} is 
\[
\det(g_{I\bar{J}})=p(x)\Theta(x)\exp(-6\phi_{1}-6\phi_{2}-2\psi_{u}-2\psi_{v})\,|w|^{-2}.
\]
This follows either by a Schur complement or by the preceding orthogonal
decomposition; the connection change of frame is triangular with determinant
one in the coframe containing $d\log w$. More explicitly, the product
of the horizontal scale factors is $h^{5}_{1}h^{5}_{2}h_{3}h_{4}=p(x)$.
A Fubini-Study metric in complex dimension $d$ has determinant $e^{-(d+1)\phi}$
in the standard affine chart, and replacing $d\log w$ by $dw$ contributes
$|w|^{-2}$ to the determinant. These three observations give the
displayed formula. Equation \eqref{eq:5.1} says 
\[
(p\Theta)'=(2-x)p.
\]
Since $\frac{dx}{ds}=\Theta$ and $f'=x$, it follows that 
\[
\frac{d}{ds}\log(p\Theta)=2-x,\quad\log(p\Theta)=2s-f(s)+\text{constant}.
\]
Substituting in the determinant gives 
\[
\log\det(g_{I\bar{J}})=-\Phi+\log|w|^{2}+\text{constant}.
\]
The last logarithm is pluriharmonic on this chart. Hence $\operatorname{Ric}(g)=g$
there, and the identity extends to all of $M$ by smoothness. In particular,
the constants $4$ and $8$ in the potential are consistent with Einstein
normalization one. Indeed, after substituting $2s-f(s)$, the logarithm
of the determinant has coefficients $-4,-8,0,-4$ on $\phi_{1},\phi_{2},\psi_{u},\psi_{v}$,
respectively, exactly the negatives of their coefficients in $\Phi$.

\subsubsection{The complete curvature tensor.}

We compute the complete tensor in a unitary frame to determine the
Calabi operator on the entire symmetric square. Let $V_{j}$ denote
the four unitary horizontal blocks, and let $e_{0}$ span the unitary
radial line. Choose a unitary frame $e_{j,i}$ in each $V_{j}$, where
$1\le j\le4$ and $1\le i\le\dim_{\mathbb{C}}V_{j}$. In each $V_{j}$,
the curvature is the Fubini-Study tensor multiplied by $a_{j}$. Between
distinct blocks $V_{j},V_{l}$ it is the scalar mixed bisectional
tensor with coefficient $b_{jl}$; the radial mixed coefficient is
$b_{j0}$. The complete list is 
\begin{align}
a_{j} & =\frac{h_{j}-\Theta}{h^{2}_{j}},\nonumber \\
b_{jl} & =-\frac{\varepsilon_{j}\varepsilon_{l}\Theta}{h_{j}h_{l}}\quad\text{for }j\ne l,\nonumber \\
b_{j0} & =-\frac{\varepsilon_{j}\Theta'}{h_{j}}+\frac{\Theta}{h^{2}_{j}},\nonumber \\
R(e_{0},\overline{e_{0}},e_{0},\overline{e_{0}}) & =-\Theta''.\label{eq:5.3}
\end{align}
Here and below primes on $\Theta$ denote derivatives with respect
to $x$. In these frames the tensor identities are 
\[
\begin{aligned}R(e_{j,i},\overline{e_{j,k}},e_{j,p},\overline{e_{j,q}}) & =a_{j}(\delta_{ik}\delta_{pq}+\delta_{iq}\delta_{pk}),\\
R(e_{j,i},\overline{e_{j,k}},e_{l,p},\overline{e_{l,q}}) & =b_{jl}\delta_{ik}\delta_{pq}\quad(j\ne l),\\
R(e_{j,i},\overline{e_{j,k}},e_{0},\overline{e_{0}}) & =b_{j0}\delta_{ik}.
\end{aligned}
\]
All components not forced by these formulas and Kähler symmetry vanish.

For completeness, these formulas follow directly from \eqref{eq:5.2}.
Move the origins of all four projective factors to the point in question
and put $\zeta=\log w$. Write $z_{j,i}$ for the coordinates of the
vector $z_{j}$ in the $j$th block. The part of the potential needed
for purely horizontal fourth derivatives is 
\[
\sum_{j}h_{j}|z_{j}|^{2}-\frac{1}{2}\sum_{j}h_{j}|z_{j}|^{4}+\frac{\Theta}{2}\left(\sum_{j}\varepsilon_{j}|z_{j}|^{2}\right)^{2}.
\]
The horizontal connection corrections vanish at the origins. Taking
fourth derivatives gives the first two coefficients in \eqref{eq:5.3}:
the coefficient of $|z_{j}|^{4}$ is $\frac{1}{2}(\Theta-h_{j})$,
while that of $|z_{j}|^{2}|z_{l}|^{2}$ for $j\ne l$ is $\varepsilon_{j}\varepsilon_{l}\Theta$.
Dividing by $h^{2}_{j}$ or $h_{j}h_{l}$ converts coordinate vectors
to unit vectors. For radial mixed components, $\partial_{\zeta}h_{j}=\varepsilon_{j}\Theta$
and $\partial_{\zeta}\partial_{\bar{\zeta}}h_{j}=\varepsilon_{j}\Theta\Theta'$.
By \eqref{eq:2.1}, the corresponding coordinate curvature is 
\[
R(\partial_{z_{j,i}},\partial_{\bar{z}_{j,i}},\partial_{\zeta},\partial_{\bar{\zeta}})=-\varepsilon_{j}\Theta\Theta'+\frac{\Theta^{2}}{h_{j}}.
\]
Division by $h_{j}\Theta$ gives $b_{j0}$. For the radial term, 
\[
\begin{aligned}\partial_{\zeta}\Theta & =\Theta\Theta',\\
\partial_{\zeta}\partial_{\bar{\zeta}}\Theta & =\Theta(\Theta')^{2}+\Theta^{2}\Theta'',\\
R_{\zeta\bar{\zeta}\zeta\bar{\zeta}} & =-\Theta(\Theta')^{2}-\Theta^{2}\Theta''+\frac{(\Theta\Theta')^{2}}{\Theta}=-\Theta^{2}\Theta''.
\end{aligned}
\]
Dividing by $\Theta^{2}$ proves the last formula in \eqref{eq:5.3}.
The lifted unitary actions on the two base factors and the unitary
actions on the two rank-two summands preserve the metric. Their stabilizers
realize independent unitary rotations of all four horizontal blocks,
fixing the point and the radial direction. Invariance under their
scalar rotations forces the number of holomorphic and antiholomorphic
indices in each block to agree. This excludes all other index arrangements;
unitary invariance within a block then gives exactly the Kronecker
tensors displayed above.

The decomposition 
\[
\begin{aligned}\operatorname{Sym}^{2}T^{1,0}M & =\bigoplus_{j}\operatorname{Sym}^{2}V_{j}\oplus\bigoplus_{j<l}(V_{j}\odot V_{l})\\
 & \qquad\oplus\bigoplus_{j}(V_{j}\odot\mathbb{C}e_{0})\oplus\mathbb{C}(e_{0}\odot e_{0})
\end{aligned}
\]
is orthogonal. Here $\odot$ denotes symmetrization, and all summands
carry the inner product inherited from the tensor product, as in Section~\ref{sec:2}.
Formula \eqref{eq:5.3} shows that $\mathcal{C}$ is scalar on these
summands, with respective eigenvalues 
\[
2a_{j},\quad2b_{jl},\quad2b_{j0},\quad-\Theta''.
\]
There are no couplings between different summands. To see the factor
2 on a mixed summand, take a symmetric matrix supported on a pair
of distinct blocks $V_{j},V_{l}$, with entries $S_{(j,i),(l,p)}=S_{(l,p),(j,i)}=z$.
Its squared norm on this pair is $2|z|^{2}$, whereas the four Kähler
permutations contribute $4b_{jl}|z|^{2}$ to $\mathcal{R}(S)$. The
quotient is $2b_{jl}$. Within one horizontal block, contraction of
the two Kronecker tensors gives $2a_{j}\lVert S\rVert^{2}$. The radial-radial
summand has only one matrix entry and therefore has eigenvalue $-\Theta''$.

\subsubsection{Sharp bounds for the Calabi operator.}

Put $D=x(4-x)$. Since $(4+x)(8-x)=32+D$ and $\frac{dD}{dx}=2(2-x)$,
the numerator of \eqref{eq:5.1} is $\frac{1}{2}\int^{D}_{0}t(32+t)^{5}\,dt$.
This identity is valid on both sides of $x=2$: the two expressions
have the same derivative with respect to $x$ and both vanish at $x=0$.
Consequently 
\begin{align}
\Theta(x) & =G(D),\nonumber \\
G(D) & =\frac{D}{2}\int^{1}_{0}v\left(\frac{32+Dv}{32+D}\right)^{5}\,dv.\label{eq:5.4}
\end{align}
In particular, $0\le\Theta\le\frac{D}{4}\le1$.

We first bound the radial holomorphic sectional curvature. Set $r=\frac{D}{32+D}$,
so $0\le r\le\frac{1}{9}$. For a fixed $v$, put $k=\frac{32+Dv}{32+D}=1-r(1-v)$.
Here primes on $G$ and $k$ denote derivatives with respect to $D$.
Then $k'=-\frac{32(1-v)}{(32+D)^{2}}$ and $k''=\frac{64(1-v)}{(32+D)^{3}}$.
Differentiation under the integral in \eqref{eq:5.4} gives 
\[
\begin{aligned}G''(D) & =5\int^{1}_{0}vk^{3}\left(kk'+2D(k')^{2}+\frac{D}{2}kk''\right)\,dv\\
 & =-\frac{5(1-r)^{3}}{32}\int^{1}_{0}v(1-v)\bigl(1-r(1-v)\bigr)^{3}\bigl(1-3r(1-v)\bigr)\,dv.
\end{aligned}
\]
All factors following $v(1-v)$ belong to $(0,1]$. Since $\int^{1}_{0}v(1-v)\,dv=\frac{1}{6}$
and $G'(0)=\frac{1}{4}$, 
\[
\begin{aligned}-\frac{5}{192} & \le G''(D)<0,\\
\frac{7}{48} & \le G'(D)\le\frac{1}{4}.
\end{aligned}
\]
For the second line we integrated the first from zero and used $D\le4$,
giving $G'(D)\ge\frac{1}{4}-\frac{5}{192}D\ge\frac{7}{48}$. The chain
rule gives 
\begin{equation}
\frac{7}{24}\le-\Theta''(x)=2G'(D)-4(4-D)G''(D)\le\frac{11}{12}.\label{eq:5.5}
\end{equation}
The lower bound uses $-G''\ge0$; the upper bound is $2\cdot\frac{1}{4}+16\cdot\frac{5}{192}=\frac{11}{12}$.
Thus $\Theta$ is concave. Its endpoint slopes are $1$ and $-1$,
so $|\Theta'|\le1$.

This single estimate controls all curvature coefficients within the
projective fiber. Integrating $-\Theta''$ from $0$ and using $\Theta(0)=0$,
$\Theta'(0)=1$, we obtain 
\[
\begin{aligned}a_{3} & =\frac{x-\Theta(x)}{x^{2}}=\frac{1}{x^{2}}\int^{x}_{0}(x-t)(-\Theta''(t))\,dt,\\
b_{30} & =\frac{\Theta(x)-x\Theta'(x)}{x^{2}}=\frac{1}{x^{2}}\int^{x}_{0}t(-\Theta''(t))\,dt.
\end{aligned}
\]
It follows from \eqref{eq:5.5} that $0<2a_{3},2b_{30}\le\frac{11}{12}$.
Reflection under $x\mapsto4-x$ gives the same bounds for $2a_{4},2b_{40}$.
Also $2b_{34}=\frac{2\Theta}{D}\in(0,\frac{1}{2}]$.

For the base blocks, $h_{1},h_{2}\ge4$. The bounds for $\Theta$
and $\Theta'$ give, for $j=1,2$, 
\[
\begin{gathered}0<2a_{j}\le\frac{1}{2},\\
-\frac{1}{2}\le-\frac{2}{h_{j}}\le2b_{j0}\le\frac{2}{h_{j}}+\frac{2}{h^{2}_{j}}\le\frac{5}{8}.
\end{gathered}
\]
For a base-fiber coefficient, cancel $x$ or $4-x$ using $\Theta\le\frac{x(4-x)}{4}$.
More explicitly, for $j=1,2$, 
\[
\begin{aligned}|2b_{j3}| & =\frac{2\Theta}{h_{j}x}\le\frac{4-x}{2h_{j}}\le\frac{1}{2},\\
|2b_{j4}| & =\frac{2\Theta}{h_{j}(4-x)}\le\frac{x}{2h_{j}}\le\frac{1}{2}.
\end{aligned}
\]
For the two base blocks, $h_{1}h_{2}=32+D$ implies 
\[
|2b_{12}|\le\frac{D}{2(32+D)}\le\frac{1}{18}.
\]
These estimates cover every summand of the symmetric square. Smoothness
extends them to the endpoints, proving 
\begin{equation}
-\frac{1}{2}\operatorname{Id}\le\mathcal{C}\le\frac{11}{12}\operatorname{Id}\label{eq:5.6}
\end{equation}
on all of $M$. The two spectral endpoints are attained. At $D=0$,
$G'(0)=\frac{1}{4}$ and $G''(0)=-\frac{5}{192}$, so $-\Theta''(0)=\frac{11}{12}$.
Moreover, 
\[
\lim_{x\to0}2b_{13}=\lim_{x\to0}-\frac{2\Theta(x)}{(4+x)x}=-\frac{1}{2}.
\]
At points approaching the stratum $U=0$, the largest eigenvalue of
$\mathcal{C}$ is at least $-\Theta''(x)$ and the smallest is at
most $2b_{13}$. The eigenvalues of a smooth Hermitian operator vary
continuously, even though the adapted frame degenerates. Together
with \eqref{eq:5.6}, these limits show that the largest and smallest
eigenvalues on this stratum are exactly $\frac{11}{12}$ and $-\frac{1}{2}$.

Rescaling $g$ by $1/2$ gives Einstein constant two and doubles the
Calabi eigenvalues. The resulting spectrum lies in $[-1,\frac{11}{6}]\subset(-2,2)$,
with both endpoints attained. This gives the claimed spectral counterexample
to homogeneity.

The quadratic curvature inequalities follow from \eqref{eq:5.6}.
For an arbitrary matrix write $A=S+K$ with $S$ symmetric and $K$
skew-symmetric, as in Section~\ref{sec:2}. These parts are orthogonal,
and the curvature contraction depends only on $S$. Combining \eqref{eq:5.6}
with \eqref{eq:2.3}, with Einstein constant one, gives 
\[
\begin{aligned}\CQB(A) & \ge\lVert K\rVert^{2}+\frac{1}{12}\lVert S\rVert^{2}\ge\frac{1}{12}\lVert A\rVert^{2},\\
\dCQB(A) & \ge\lVert K\rVert^{2}+\frac{1}{2}\lVert S\rVert^{2}\ge\frac{1}{2}\lVert A\rVert^{2}.
\end{aligned}
\]
Symmetric eigentensors for the two extremal eigenvalues give equality,
so both constants are optimal. After rescaling to Einstein constant
two, the corresponding sharp constants are $\frac{1}{6}$ and $1$.

\subsubsection{Nonhomogeneity.}

Choose one copy of $L_{2}$ in $L^{\oplus2}_{1}\oplus L^{\oplus2}_{2}$.
Its lines define a section $S\cong B$. The projective-bundle tangent
formula from Section~\ref{sec:2} gives its normal bundle: 
\[
N_{S/M}=\operatorname{Hom}(L_{2},L^{\oplus2}_{1}\oplus L_{2})\cong\mathcal{O}_{B}(-1,1)^{\oplus2}\oplus\mathcal{O}_{B}.
\]
Let $C$ be a line in the first $\mathbb{P}^{5}$ factor of this section,
with the second factor fixed. The normal exact sequence gives a surjection
\[
TM|_{C}\rightarrow N_{S/M}|_{C}\rightarrow\mathcal{O}_{\mathbb{P}^{1}}(-1).
\]
Lemma~\ref{lemma:2.2} now proves that $M$ is not homogeneous. This
completes the proof of Theorem~\ref{theorem:1.3}(1).

\subsection{Lower-dimensional counterexamples.}

For $n=2,3$, let 
\[
E_{n}=\mathcal{O}^{\oplus2}_{\mathbb{P}^{n}}\oplus\mathcal{O}_{\mathbb{P}^{n}}(-1),\quad X_{n}=\mathbb{P}_{\mathrm{lines}}(E_{n}).
\]
The fiber is $\mathbb{P}^{2}$, so $X_{2}$ and $X_{3}$ have complex
dimensions 4 and 5. In the convention that projectivization parametrizes
one-dimensional quotients, the same manifold is written $\mathbb{P}(\mathcal{O}^{\oplus2}\oplus\mathcal{O}(1))$.
This explains the sign of the line bundle in the comparison with \cite{NWZ21}.

\subsubsection{The global metric.}

We define the metric through a projective embedding, which also makes
its global positivity immediate. Let $\pi:X_{n}\to\mathbb{P}^{n}$
be the bundle projection. The inclusion $E_{n}\subset\mathbb{P}^{n}\times\mathbb{C}^{n+3}$
induces a map $\beta:X_{n}\to\mathbb{P}^{n+2}$, and $(\pi,\beta)$
identifies $X_{n}$ with the closed incidence submanifold 
\[
\{([z],[a:b:v])\in\mathbb{P}^{n}\times\mathbb{P}^{n+2}:v\in\mathbb{C}z\}.
\]
The condition $v\in\mathbb{C}z$ is imposed by the equations $z_{i}v_{j}-z_{j}v_{i}=0$.
For a fixed $[z]$, the remaining data are a line in $\mathbb{C}^{2}\oplus\mathbb{C}z$,
precisely a line in the fiber of $E_{n}$. Local trivializations identify
this incidence variety smoothly and holomorphically with $\mathbb{P}_{\mathrm{lines}}(E_{n})$.
Thus $(\pi,\beta)$ is an embedding.

Give both projective factors their Fubini-Study forms, normalized
to have holomorphic sectional curvature two, and set 
\[
\omega_{n}=4\pi^{*}\omega_{\mathrm{FS}}+\beta^{*}\omega_{\mathrm{FS}}.
\]
This is the restriction of the product Kähler form to a complex submanifold,
so it defines a smooth Kähler metric on $X_{n}$. To compute it locally,
write $\phi(z)=\log(1+|z|^{2})$ on the base. The standard frame $(1,z)$
of $\mathcal{O}(-1)$ has squared norm $e^{\phi}=1+|z|^{2}$. In the
affine fiber chart represented by $(1,u,w)$, the second projective
coordinate of the embedding is $[1:u:w(1,z)]$. Pulling back the two
Fubini-Study potentials gives 
\[
\Phi(z,u,w)=4\log(1+|z|^{2})+\log\bigl(1+|u|^{2}+(1+|z|^{2})|w|^{2}\bigr).
\]
The $\operatorname{SU}(n+1)$ action on the base lifts to the tautological
line bundle, preserves its induced Hermitian metric, and acts trivially
on the two trivial summands. Hence it preserves $\omega_{n}$ and
allows us to move the base point to $z=0$. At that point put 
\[
t=\frac{|w|^{2}}{1+|u|^{2}+|w|^{2}},\quad h=4+t.
\]
The number $t\in[0,1]$ is globally the fraction of the squared norm
of a line representative contributed by the tautological summand.
The metric at $z=0$ has the orthogonal splitting 
\begin{equation}
g_{n}=h\sum^{n}_{i=1}dz_{i}\,d\bar{z}_{i}+g_{\mathrm{FS},\mathbb{P}^{2}}.\label{eq:5.7}
\end{equation}
The splitting follows because the first derivatives of $\phi$ vanish
at $z=0$. Its horizontal coefficient is $4+t$, and its vertical
restriction is the Fubini-Study metric on the fiber.

\subsubsection{Curvature from the incidence embedding.}

Write $H$ and $V$ for the horizontal and vertical spaces in \eqref{eq:5.7}.

\begin{lemma}\label{lemma:5.1} At every point there are unitary
frames $e_{1},\ldots,e_{n}$ of $H$ and $\xi_{0},\xi_{1}$ of $V$
such that, for $1\le i,j,k,l\le n$ and $p,q,r,s\in\{0,1\}$, 
\begin{align}
R(e_{i},\overline{e_{j}},e_{k},\overline{e_{l}}) & =a(\delta_{ij}\delta_{kl}+\delta_{il}\delta_{kj}),\nonumber \\
R(\xi_{p},\overline{\xi_{q}},\xi_{r},\overline{\xi_{s}}) & =\delta_{pq}\delta_{rs}+\delta_{ps}\delta_{rq},\nonumber \\
R(e_{i},\overline{e_{j}},\xi_{p},\overline{\xi_{q}}) & =\delta_{ij}\delta_{pq}c_{p},\label{eq:5.8}
\end{align}
where 
\[
a=\frac{4+t^{2}}{(4+t)^{2}},\quad c_{0}=1-\frac{20}{(4+t)^{2}},\quad c_{1}=\frac{t}{4+t}.
\]
All components not obtained from \eqref{eq:5.8} by Kähler symmetries
vanish. In particular $R(e_{i},\overline{\xi_{p}},e_{j},\overline{\xi_{q}})=0$.
\end{lemma} 
\begin{proof}
We apply the Gauss equation to the incidence embedding. By the unitary
symmetries, the first projective coordinate can be taken to be $[1:0:\cdots:0]$.
Choose orthonormal representatives $a_{0},a_{1}$ in the two trivial
summands and $c$ in the tautological summand over this point. Put
$r=\sqrt{t}$, $s=\sqrt{1-t}$. The point is represented by $v=sa_{0}+rc$,
and a unitary vertical frame is 
\[
\begin{aligned}\xi_{0} & =ra_{0}-sc,\\
\xi_{1} & =a_{1}.
\end{aligned}
\]
We identify $T_{[v]}\mathbb{P}^{n+2}$ with $v^{\perp}\subset\mathbb{C}^{n+3}$.
Let $f^{(1)}_{i}$ be ambient unit vectors in the first projective
factor with metric $4g_{\mathrm{FS}}$, and let $f^{(2)}_{i}$ be
the unit vectors in the second factor corresponding to the remaining
base coordinates. The differential of the incidence embedding sends
$\partial_{z_{i}}$ to $2f^{(1)}_{i}+rf^{(2)}_{i}$. An orthonormal
horizontal frame and an orthonormal normal frame are therefore 
\[
\begin{aligned}e_{i} & =\frac{2f^{(1)}_{i}+rf^{(2)}_{i}}{\sqrt{h}},\\
\nu_{i} & =\frac{-rf^{(1)}_{i}+2f^{(2)}_{i}}{\sqrt{h}}.
\end{aligned}
\]
Here vectors in the two ambient factors are identified with their
images in the product tangent space.

In centered affine projective coordinates, the second component of
the incidence embedding has representative 
\[
v+\zeta_{0}\xi_{0}+\zeta_{1}\xi_{1}+(r-s\zeta_{0})\sum_{i}z_{i}f^{(2)}_{i}.
\]
The first component is $[1:z_{1}:\cdots:z_{n}]$. The Christoffel
symbols of both ambient metrics vanish at the origins of these charts.
The only quadratic terms are $-s\zeta_{0}z_{i}f^{(2)}_{i}$. Projecting
their second derivatives onto the normal frame gives 
\[
\mathrm{II}(e_{i},\xi_{0})=\mathrm{II}(\xi_{0},e_{i})=-\frac{2\sqrt{1-t}}{h}\nu_{i},
\]
and all other second fundamental coefficients vanish. Here the unnormalized
mixed derivative is $-sf^{(2)}_{i}$; its normal projection is $-\frac{2s}{\sqrt{h}}\nu_{i}$,
and normalizing the horizontal input contributes the additional factor
$\frac{1}{\sqrt{h}}$. The vertical vectors already have unit length.
This accounts for the coefficient in the second fundamental form.

The ambient horizontal curvature coefficient is $\frac{4+t^{2}}{h^{2}}$.
Indeed, a unit horizontal vector has components of lengths $\frac{2}{\sqrt{h}}$
and $\frac{\sqrt{t}}{\sqrt{h}}$ in the two ambient factors. Their
Fubini-Study curvature coefficients are $\frac{1}{4}$ and $1$, respectively,
giving $\frac{1}{4}\frac{16}{h^{2}}+\frac{t^{2}}{h^{2}}$. The vertical
curvature is the Fubini-Study tensor, and the ambient mixed coefficient
is $\frac{t}{h}$ in both vertical directions. The Gauss equation
therefore changes only the mixed coefficient in the $\xi_{0}$ direction:
\[
c_{0}=\frac{t}{h}-\frac{4(1-t)}{h^{2}}=1-\frac{20}{h^{2}}.
\]
In our convention the Kähler-Gauss equation subtracts $\langle\mathrm{II}(X,Z),\mathrm{II}(Y,W)\rangle$
from the ambient component $R_{X\bar{Y}Z\bar{W}}$. Since the only
nonzero second fundamental coefficients pair a horizontal vector with
$\xi_{0}$, there is no correction to four horizontal or four vertical
indices. There is also no correction to $R(e_{i},\overline{\xi_{p}},e_{j},\overline{\xi_{q}})$,
whose holomorphic pair has two horizontal indices. Its ambient component
vanishes because the horizontal and vertical vectors in the second
factor are orthogonal. The coefficient in the $\xi_{1}$ direction
remains $c_{1}=\frac{t}{h}$. The same calculation gives every vanishing
component asserted in \eqref{eq:5.8}. The formulas extend to the
endpoints by smoothness. 
\end{proof}

Contracting \eqref{eq:5.8} gives the Ricci eigenvalues 
\begin{align}
\rho_{n} & =(n+1)a+c_{0}+c_{1}=\frac{4n+12t+(n+3)t^{2}}{(4+t)^{2}} &  & \text{on }H,\nonumber \\
\sigma_{p} & =3+nc_{p} &  & \text{on }\mathbb{C}\xi_{p}.\label{eq:5.9}
\end{align}
For $n=2,3$, the following elementary bounds will suffice: 
\begin{gather}
0<a\le\frac{1}{4},\quad-\frac{1}{4}\le c_{0}\le\frac{1}{5},\quad0\le c_{1}\le\frac{1}{5},\nonumber \\
\rho_{n}\ge\frac{n}{4},\quad\sigma_{p}\ge3-\frac{n}{4}.\label{eq:5.10}
\end{gather}
Indeed, $a\le\frac{1}{4}$ is equivalent to $8t-3t^{2}\ge0$. The
function $c_{0}=1-\frac{20}{(4+t)^{2}}$ is increasing, and $c_{1}$
is increasing as well. For the Ricci bound, subtracting $\frac{n}{4}(4+t)^{2}$
from the numerator in \eqref{eq:5.9} leaves 
\[
(12-2n)t+\left(3+\frac{3n}{4}\right)t^{2}\ge0.
\]

\subsubsection{The Calabi operator and inequalities involving Ricci eigenvalues.}

Lemma~\ref{lemma:5.1} makes the Calabi operator scalar on the orthogonal
summands 
\[
\operatorname{Sym}^{2}H,\qquad\operatorname{Sym}^{2}V,\qquad H\odot\mathbb{C}\xi_{0},\qquad H\odot\mathbb{C}\xi_{1},
\]
with respective eigenvalues $2a$, $2$, $2c_{0}$, and $2c_{1}$.
In particular, 
\[
-\frac{1}{2}\operatorname{Id}\le\mathcal{C}\le2\operatorname{Id}.
\]
Since these metrics are not Einstein, this spectral interval alone
does not give the desired quadratic curvature inequalities. We use
\eqref{eq:2.2}, retaining the individual Ricci eigenvalues in \eqref{eq:5.9},
and estimate both dimensions at once. Decompose $A$ into its $HH$,
$VV$, and two mixed blocks. The first block subscript specifies the
output space: $(A_{HV})_{ip}$ has horizontal row $i$ and vertical
column $p$, while $(A_{VH})_{pi}$ has vertical row $p$ and horizontal
column $i$. Put $S_{H}=\frac{A_{HH}+A^{\mathsf{T}}_{HH}}{2}$ and
$S_{V}=\frac{A_{VV}+A^{\mathsf{T}}_{VV}}{2}$. Lemma~\ref{lemma:5.1}
gives 
\[
\begin{aligned}\mathcal{R}(A) & =2a\lVert S_{H}\rVert^{2}+2\lVert S_{V}\rVert^{2}\\
 & \quad+\sum_{i,p}c_{p}|(A_{HV})_{ip}+(A_{VH})_{pi}|^{2}.
\end{aligned}
\]
In particular there are no couplings between these three types of
matrix block. The factor in the last term accounts for all four Kähler
permutations of a mixed component. For fixed horizontal and vertical
indices $i,p$, the contribution of $(A_{HV})_{ip}$ and $(A_{VH})_{pi}$
is 
\[
\begin{aligned} & c_{p}\bigl(|(A_{HV})_{ip}|^{2}+|(A_{VH})_{pi}|^{2}\\
 & \qquad+(A_{HV})_{ip}\overline{(A_{VH})_{pi}}+(A_{VH})_{pi}\overline{(A_{HV})_{ip}}\bigr)\\
 & \quad=c_{p}|(A_{HV})_{ip}+(A_{VH})_{pi}|^{2}.
\end{aligned}
\]
The corresponding Ricci eigenvalues are $\rho_{n}$ and $\sigma_{p}$,
respectively. Keeping these entries independent ensures that the estimates
apply to all linear maps.

On the horizontal block $\mathcal{M}=\rho_{n}\lVert A_{HH}\rVert^{2}$.
Since $a>0$ and $\lVert S_{H}\rVert\le\lVert A_{HH}\rVert$, we obtain
\[
\begin{aligned}\mathcal{M}-\mathcal{R} & \ge(\rho_{n}-2a)\lVert A_{HH}\rVert^{2},\\
\mathcal{M}+\mathcal{R} & \ge\rho_{n}\lVert A_{HH}\rVert^{2}.
\end{aligned}
\]
On the vertical block, the Ricci bound in \eqref{eq:5.10} gives 
\[
\mathcal{M}=\sum^{1}_{p,q=0}\sigma_{p}|(A_{VV})_{pq}|^{2}\ge\left(3-\frac{n}{4}\right)\lVert A_{VV}\rVert^{2}.
\]
Using $\lVert S_{V}\rVert\le\lVert A_{VV}\rVert$, we therefore have
\[
\begin{aligned}\mathcal{M}-\mathcal{R} & \ge\left(3-\frac{n}{4}\right)\lVert A_{VV}\rVert^{2}-2\lVert S_{V}\rVert^{2}\ge\left(1-\frac{n}{4}\right)\lVert A_{VV}\rVert^{2},\\
\mathcal{M}+\mathcal{R} & \ge\left(3-\frac{n}{4}\right)\lVert A_{VV}\rVert^{2}.
\end{aligned}
\]
For each mixed pair, the Cauchy-Schwarz inequality gives 
\[
|(A_{HV})_{ip}+(A_{VH})_{pi}|^{2}\le\frac{3}{2}|(A_{HV})_{ip}|^{2}+3|(A_{VH})_{pi}|^{2}.
\]
Indeed, the difference is $\left|\frac{(A_{HV})_{ip}}{\sqrt{2}}-\sqrt{2}(A_{VH})_{pi}\right|^{2}\ge0$.
Together with $|c_{p}|\le\frac{1}{4}$ and \eqref{eq:5.10}, this
yields, for either sign, 
\[
\begin{aligned} & \rho_{n}|(A_{HV})_{ip}|^{2}+\sigma_{p}|(A_{VH})_{pi}|^{2}\pm c_{p}|(A_{HV})_{ip}+(A_{VH})_{pi}|^{2}\\
 & \quad\ge\left(\frac{n}{4}-\frac{3}{8}\right)|(A_{HV})_{ip}|^{2}+\left(3-\frac{n}{4}-\frac{3}{4}\right)|(A_{VH})_{pi}|^{2}\\
 & \quad=\frac{2n-3}{8}|(A_{HV})_{ip}|^{2}+\frac{9-n}{4}|(A_{VH})_{pi}|^{2}.
\end{aligned}
\]
The right side is at least $\frac{1}{8}(|(A_{HV})_{ip}|^{2}+|(A_{VH})_{pi}|^{2})$
for $n=2$ and $\frac{3}{8}(|(A_{HV})_{ip}|^{2}+|(A_{VH})_{pi}|^{2})$
for $n=3$.

\noindent The 5-fold: For $n=3$, \eqref{eq:5.10} gives $\rho_{3}-2a\ge\frac{3}{4}-\frac{1}{2}=\frac{1}{4}$.
The preceding horizontal, vertical, and mixed estimates thus bound
each contribution to $\CQB$ and $\dCQB$ below by $\frac{1}{4}$
times the squared norm of that matrix block. Adding these bounds proves
the two inequalities in Theorem~\ref{theorem:1.3}(3).

\noindent The 4-fold: For $n=2$, the horizontal $\CQB$ margin has
the exact value 
\[
\rho_{2}-2a=\frac{3t}{4+t}\ge0.
\]
The vertical $\CQB$ contribution is at least $\frac{1}{2}\lVert A_{VV}\rVert^{2}$,
and each mixed pair contributes at least $\frac{1}{8}(|(A_{HV})_{ip}|^{2}+|(A_{VH})_{pi}|^{2})$.
Thus $\CQB$ is nonnegative. For $\dCQB$, the horizontal, vertical,
and mixed lower bounds are $\frac{1}{2}$, $\frac{5}{2}$, and $\frac{1}{8}$
times the respective squared norms. Their sum proves Theorem~\ref{theorem:1.3}(2).
The $\CQB$ bounds also show that $\CQB$ is strictly positive wherever
$t>0$. At $t=0$ its kernel consists precisely of symmetric matrices
supported on $H\times H$: the horizontal skew-symmetric contribution
and every other block remain positive. More explicitly, at $t=0$
we have $a=\frac{1}{4}$ and $\rho_{2}=\frac{1}{2}$, so the horizontal
$\CQB$ contribution is $\frac{1}{2}\lVert A_{HH}-S_{H}\rVert^{2}$.
The locus $t=0$ is precisely the divisor $\mathbb{P}_{\mathrm{lines}}(\mathcal{O}^{\oplus2})\subset X_{2}$.

\subsubsection{Nonhomogeneity and the Einstein condition.}

Choose a fixed line in the trivial rank-two summand of $E_{n}$. It
defines a section $S\cong\mathbb{P}^{n}$ of $\pi$. The tangent space
of a projective fiber at a line $L$ is $\operatorname{Hom}(L,E_{n}/L)$;
therefore the normal bundle of this section is 
\[
N_{S/X_{n}}\cong\mathcal{O}_{\mathbb{P}^{n}}\oplus\mathcal{O}_{\mathbb{P}^{n}}(-1).
\]
If $C\cong\mathbb{P}^{1}$ is a line in $S$, the normal exact sequence
gives a surjection 
\[
TX_{n}|_{C}\longrightarrow N_{S/X_{n}}|_{C}\longrightarrow\mathcal{O}_{C}(-1).
\]
Lemma~\ref{lemma:2.2} shows that neither $X_{2}$ nor $X_{3}$ is
homogeneous.

Finally, at $t=0$ the vertical Ricci eigenvalues are 
\[
\sigma_{0}=3-\frac{n}{4},\quad\sigma_{1}=3.
\]
They are distinct, so $g_{n}$ is not Einstein. This completes the
proof of Theorem~\ref{theorem:1.3}.


\begin{thebibliography}{ACGTF08}
\bibitem[AB03]{AB03} Hassan Azad and Indranil Biswas. Quasi-potentials
and Kähler-Einstein metrics on flag manifolds, II. \emph{J. Algebra},
269(2):480--491, 2003.

\bibitem[ACGTF08]{ACGTF08} Vestislav Apostolov, David M. J. Calderbank,
Paul Gauduchon, and Christina W. Tønnesen-Friedman. Hamiltonian 2-forms
in Kähler geometry. III. Extremal metrics and stability. \emph{Invent.
Math.}, 173(3):547--601, 2008.

\bibitem[AHZ18]{AHZ18} Angelynn Alvarez, Gordon Heier, and Fangyang
Zheng. On projectivized vector bundles and positive holomorphic sectional
curvature. \emph{Proc. Amer. Math. Soc.}, 146(7):2877--2882, 2018.

\bibitem[Bes87]{Bes87} Arthur L. Besse. \emph{Einstein Manifolds},
volume~10 of \emph{Ergebnisse der Mathematik und ihrer Grenzgebiete
(3)}. Springer-Verlag, Berlin, 1987. Chapter~8, pp.~208--234.

\bibitem[Bor53]{Bor53} Armand Borel. Sur la cohomologie des espaces
fibrés principaux et des espaces homogènes de groupes de Lie compacts.
\emph{Ann. of Math. (2)}, 57(1):115--207, 1953.

\bibitem[Bor54]{Bor54} Armand Borel. Kählerian coset spaces of semisimple
Lie groups. \emph{Proc. Nat. Acad. Sci. U.S.A.}, 40(12):1147--1151,
1954.

\bibitem[Bor60]{Bor60} Armand Borel. On the curvature tensor of the
Hermitian symmetric manifolds. \emph{Ann. of Math. (2)}, 71(3):508--521,
1960.

\bibitem[BR62]{BR62} Armand Borel and Reinhold Remmert. Über kompakte
homogene Kählersche Mannigfaltigkeiten. \emph{Math. Ann.}, 145:429--439,
1962.

\bibitem[CP91]{CP91} Frédéric Campana and Thomas Peternell. Projective
manifolds whose tangent bundles are numerically effective. \emph{Math.
Ann.}, 289:169--187, 1991.

\bibitem[CT13]{CT13} Albert Chau and Luen-Fai Tam. Kähler C-spaces
and quadratic bisectional curvature. \emph{J. Differential Geom.},
94(3):409--468, 2013.

\bibitem[CV60]{CV60} Eugenio Calabi and Edoardo Vesentini. On compact,
locally symmetric Kähler manifolds. \emph{Ann. of Math. (2)}, 71(3):472--507,
1960.

\bibitem[GO26]{GO26} Giovane Galindo and Ailton R. Oliveira. Curvature
positivity for Kähler and quasi-Kähler flag manifolds. \emph{Indag.
Math. (N.S.)}, 2026. Published online.

\bibitem[Hum72]{Hum72} James E. Humphreys. \emph{Introduction to
Lie Algebras and Representation Theory}, volume 9 of \emph{Graduate
Texts in Mathematics}. Springer-Verlag, New York, 1972.

\bibitem[Ito78]{Ito78} Mitsuhiro Itoh. On curvature properties of
Kähler C-spaces. \emph{J. Math. Soc. Japan}, 30(1):39--71, 1978.

\bibitem[JZ26]{JZ26} Run-Qiang Jian and Zhuhong Zhang. The cross
quadratic bisectional curvature of Kähler $C$-spaces of classical
type. arXiv:2610.39012v1, 2026. \url{https://arxiv.org/abs/2610.39012v1}.

\bibitem[Kad23]{Kad23} Fayadh Kadhem. A cluster structure on the
coordinate ring of partial flag varieties. \emph{J. Algebra}, 628:328--349,
2023. \url{https://doi.org/10.1016/j.jalgebra.2023.03.021}.

\bibitem[Kos55]{Kos55} Jean-Louis Koszul. Sur la forme hermitienne
canonique des espaces homogènes complexes. \emph{Canad. J. Math.},
7:562--576, 1955.

\bibitem[KS88]{KS88} Norihito Koiso and Yusuke Sakane. Nonhomogeneous
Kähler-Einstein metrics on compact complex manifolds. II. \emph{Osaka
J. Math.}, 25(4):933--959, 1988.

\bibitem[Loh19]{Loh19} Simon Peter Lohove. \emph{Holomorphic curvature
of Kähler Einstein metrics on generalised flag manifolds}. PhD thesis,
Università degli Studi di Firenze, 2019. \url{https://hdl.handle.net/2158/1151431}.

\bibitem[Mag22]{Mag22} Gunnar Þór Magnússon. Degenerate Hermitian
geometry and curvatures of holomorphic fibrations. arXiv:2210.01885,
2022.

\bibitem[Mat72]{Mat72} Yozo Matsushima. Remarks on Kähler--Einstein
manifolds. \emph{Nagoya Math. J.}, 46:161--173, 1972.

\bibitem[MZ86]{MZ86} Ngaiming Mok and Jia-Qing Zhong. Curvature characterization
of compact Hermitian symmetric spaces. \emph{J. Differential Geom.},
23:15--67, 1986.

\bibitem[Mok88]{Mok88} Ngaiming Mok. The uniformization theorem for
compact Kähler manifolds of nonnegative holomorphic bisectional curvature.
\emph{J. Differential Geom.}, 27(2):179--214, 1988.

\bibitem[Mor79]{Mor79} Shigefumi Mori. Projective manifolds with
ample tangent bundles. \emph{Ann. of Math. (2)}, 110(3):593--606,
1979.

\bibitem[Ni21]{Ni21} Lei Ni. The fundamental group, rational connectedness
and the positivity of Kähler manifolds. \emph{J. Reine Angew. Math.},
774:267--299, 2021.

\bibitem[Nom54]{Nom54} Katsumi Nomizu. Invariant affine connections
on homogeneous spaces. \emph{Amer. J. Math.}, 76(1):33--65, 1954.

\bibitem[NWZ21]{NWZ21} Lei Ni, Qingsong Wang, and Fangyang Zheng.
Manifolds with positive orthogonal Ricci curvature. \emph{Amer. J.
Math.}, 143(3):833--857, 2021.

\bibitem[NZ25]{NZ25} Lei Ni and Fangyang Zheng. Kähler manifolds
and cross quadratic bisectional curvature. \emph{J. Eur. Math. Soc.},
27(5):1715--1738, 2025.

\bibitem[PS25]{PS25} Peter Petersen and James Stanfield. Poincaré
polynomials and curvature operators of symmetric spaces. To appear
in \emph{MATRIX Annals}, 2025. \url{https://arxiv.org/abs/2508.11236}.

\bibitem[Sak86]{Sak86} Yusuke Sakane. Examples of compact Einstein
Kähler manifolds with positive Ricci tensor. \emph{Osaka J. Math.},
23:585--616, 1986.

\bibitem[SY80]{SY80} Yum-Tong Siu and Shing-Tung Yau. Compact Kähler
manifolds of positive bisectional curvature. \emph{Invent. Math.},
59(2):189--204, 1980.

\bibitem[Tian87]{Tian87} Gang Tian. On Kähler--Einstein metrics
on certain Kähler manifolds with $c_{1}(M)>0$. \emph{Invent. Math.},
89(2):225--246, 1987.

\bibitem[Yan18]{Yan18} Xiaokui Yang. RC-positivity, rational connectedness
and Yau's conjecture. \emph{Camb. J. Math.}, 6(2):183--212, 2018.

\bibitem[YZ19]{YZ19} Bo Yang and Fangyang Zheng. Hirzebruch manifolds
and positive holomorphic sectional curvature. \emph{Ann. Inst. Fourier
(Grenoble)}, 69(6):2589--2634, 2019.

\end{thebibliography}
\end{document}